\documentclass{amsart}

\date{}
\usepackage{amsmath}
\usepackage{amsfonts}
\usepackage{amssymb}
\usepackage{amsthm}
\usepackage[all]{xy}
\usepackage[CJKbookmarks]{hyperref}
\usepackage{fancyhdr}
\usepackage{tikz-cd}
\usepackage{geometry}

\counterwithin{equation}{section}

\newtheorem{definition}{Definition}[section]
\newtheorem{theorem}[definition]{Theorem}

\newtheorem{lemma}[definition]{Lemma}
\newtheorem{proposition}[definition]{Proposition}
\newtheorem{corollary}[definition]{Corollary}
\newtheorem{remark}[definition]{Remark}
\newtheorem{example}[definition]{Example}

\makeatletter\renewcommand\c@equation\c@definition\makeatother

\DeclareMathOperator{\Coker}{Coker}
 
\DeclareMathOperator{\add}{add}

\DeclareMathOperator{\depth}{depth}

\DeclareMathOperator{\End}{End}
\DeclareMathOperator{\gEnd}{\underline{End}}
\DeclareMathOperator{\Ext}{Ext}
\DeclareMathOperator{\gExt}{\underline{Ext}}
\DeclareMathOperator{\fd}{fd}
\DeclareMathOperator{\fg}{fg}

\DeclareMathOperator{\gr}{gr}

\DeclareMathOperator{\gldim}{gldim}

\DeclareMathOperator{\Gr}{Gr}

\DeclareMathOperator{\Hom}{Hom}
\DeclareMathOperator{\gHom}{\underline{Hom}}

\DeclareMathOperator{\id}{id}
\DeclareMathOperator{\idim}{idim}

\DeclareMathOperator{\im}{Im}
\DeclareMathOperator{\Ker}{Ker}

\DeclareMathOperator{\Lim}{Lim\,}

\DeclareMathOperator{\MCM}{MCM}

\DeclareMathOperator{\pdim}{pdim}

\DeclareMathOperator{\proj}{proj}

\DeclareMathOperator{\soc}{soc}

\DeclareMathOperator{\Tor}{Tor}

\DeclareMathOperator{\D}{\mathbf{D}} 
\DeclareMathOperator{\K}{\mathbf{K}} 

\DeclareMathOperator{\cA}{\mathcal{A}} 
\DeclareMathOperator{\cB}{\mathcal{B}} 
\DeclareMathOperator{\cC}{\mathcal{C}}

\begin{document}

\setcounter{tocdepth}{1}

\title{Commonly Graded Algebras and Their Homological Properties}
\author{Haonan Li}
\author{Quanshui Wu}
\address{School of Mathematical Sciences \\
Fudan University \\
Shanghai, 200433 \\
 China}
\email{lihn@fudan.edu.cn, qswu@fudan.edu.cn}
\thanks{This research has been supported by the NSFC (Grant No. 12471032) and the National Key Research and Development Program of China (Grant No. 2020YFA0713200).}


\newgeometry{left=3.18cm,right=3.18cm,top=2.54cm,bottom=2.54cm}

\begin{abstract}
    In this article, we study bounded-below locally finite $\mathbb{Z}$-graded algebras, which are referred to as commonly graded algebras in literature. Commonly graded algebras have almost similar theory as that of connected graded algebras, but sometimes the results need different methods of proof.  We give several characterizations of commonly graded AS-Gorenstein algebras, and show that any noetherian commonly graded AS-Gorenstein algebra admits a balanced dualizing complex.
    We then study (skew) Calabi-Yau properties of commonly graded algebras, and  give an example of graded algebra which is skew Calabi-Yau in ungraded sense but not in graded sense. We demonstrate that a noetherian commonly graded algebra is AS-regular if and only if the bounded derived category of its finite-dimensional graded modules constitutes a ``twisted" Calabi-Yau category. At the end of the article, we prove that the Auslander-Buchsbaum formula, along with the Bass theorem and the No-Hole theorem  hold for commonly graded algebras under appropriate conditions.
\end{abstract}

\maketitle

\tableofcontents

\section{Introduction}
This article serves as a preparation for the study of noncommutative resolutions of Artin-Schelter Gorenstein (AS-Gorenstein for short) isolated singularities. The concept of a noncommutative resolution for a normal Gorenstein domain was initially introduced by Van den Bergh \cite{VdB2}, characterized as the endomorphism ring of a specified reflexive module. Analogously, a noncommutative resolution $B$ of an AS-Gorenstein isolated singularity $A$ is defined as the graded endomorphism ring $B$ of some finitely generated graded $A$-module such that $B$ is an AS-regular algebra; see \cite{LSW,LW2} for further details.
However, for a finitely generated graded module $M$ over a noetherian connected graded algebra $A$,  the graded endomorphism ring $\gEnd_A(M)$ usually is not 
an $\mathbb{N}$-graded algebra, but only a bounded-below locally finite $\mathbb{Z}$-graded algebra. Algebras of this nature are called commonly graded algebras, as per the terminology used in \cite{CKWZ}.
To study the noncommutative resolutions of AS-Gorenstein isolated singularities, it is necessary to develop further the theory of commonly graded algebras. 

Locally finite $\mathbb{N}$-graded AS-regular algebras were investigated by Reyes and Rogalski \cite{RR}. 
In the context of commonly graded AS-Gorenstein algebras, 
we show that they share almost similar properties as the connected graded one. 
However, the original proofs applicable to connected graded algebras are not always available in commonly graded case, necessitating new proofs.
We first give some characterizations of commonly graded AS-Gorenstein algebras (Theorems \ref{equvalent definitions of AS-Gorenstein}, \ref{another equivalent definition of GAS}, \ref{supplement of equivalent definition of GAS} and \ref{GAS-Gorenstein and d-sCY}). 

\begin{theorem}[Theorem \ref{equvalent definitions of AS-Gorenstein}]
Let $A$ be a noetherian commonly graded algebra with finite left and right injective dimension $d$. Let $S=A/J$ where $J$ is the graded Jacobson radical of $A$. Suppose that $\gExt^i_A(S,A)=0$ for $i\neq d$. Then the following are equivalent.
\begin{itemize}
\item[(1)] $A$ is commonly graded AS-Gorenstein. 
\item[(2)] $\gExt^d_A(-,A)$ gives a contravariant equivalence between the categories  $\fd A$ and $\fd A^o$.
\item[(3)] $\gExt^d_A(-,A)$ gives a contravariant equivalence between the categories $\gr S$ and $\gr S^o$ where $\gr S$ (resp. $\gr S^o$) is regarded as the full subcategory of $\gr A$ (resp. $\gr A^o$) consisting of finitely generated graded semisimple objects via the algebra morphism $A\to A/J$.
\item[(4)] $\gExt^d_A(S,A)J=0$, and there is an invertible graded $(S,S)$-bimodule $V$ such that $V \cong \gExt^d_A(S,A)$ as graded $S^o$-modules.
\item[(5)] $\gExt^d_A(S,A)J=0$, and $\gExt^d_A(S,A)$ is invertible as a graded $(S,S)$-bimodule.
\end{itemize}
\end{theorem}

Balanced dualizing complexes are crucial tools to study noncommutative graded algebras.
The existence of balanced dualizing complexes for connected graded AS-Gorenstein algebras was proved by Yekutieli \cite{Y}.
A general existence theorem of balanced dualizing complexes over connected graded algebras was given by Van den Bergh \cite{VdB}. 
However, while considering a commonly graded algebra $A$, there are challenges: 
firstly, the isomorphism classes of graded simple modules are not unique and are not concentrated on a single degree; 
secondly, the invertible graded $(A,A)$-bimodules are not necessarily free. 
Consequently, the approach outlined in \cite{Y} cannot be directly applied to commonly graded algebras. 
In the proof of the existence of balanced dualizing complexes for noetherian commonly graded AS-Gorenstein algebras in Theorem \ref{existence of baalanced dualizing complex of GAS}, 
we utilize the structure of the minimal graded injective resolution of $A$ (Proposition \ref{gExt(M,A) cong gHom(M,I'd) as right module}).

\begin{theorem}[Theorem \ref{existence of baalanced dualizing complex of GAS}]
Any noetherian commonly graded AS-Gorenstein algebra admits a balanced dualizing complex.
\end{theorem}

As a consequence of the Van den Bergh existence theorem of balanced dualizing complexes (Theorem \ref{existence of balanced dualizing complex}), it can be deduced that any noetherian commonly graded AS-Gorenstein algebra satisfies the $\chi$-condition and has a finite local cohomological dimension (Corollary \ref{corollary of existence of balanced dualizing complex}).

There are several kinds of generalizations of AS-Gorenstein (regular) algebras in literature such as \cite{V,V2,RRZ,MM}. We discuss their relations and the differences. A special kind of graded algebras, called canonical commonly graded AS-Gorenstein (see Definition \ref{definition of canonical}), plays a critical role in this regard (Theorem \ref{GAS-Gorenstein and AS over A0}, \ref{GAS-regular and AS over A0} and Remark \ref{GAS-Gorenstein in RRZ}). 

We also concern the Calabi-Yau property for commonly graded ($\mathbb{N}$-graded) algebras and their derived categories. Calabi-Yau algebras \cite{G} (see Definition \ref{definition of CY of dimension d}) are closely related to AS-regular algebras. A connected graded algebra $A$ is a skew Calabi-Yau algebra if and only if $A$ is an AS-regular algebra \cite{RRZ}. Moreover, for any locally finite $\mathbb{N}$-graded algebra $A$, $A$ is a twisted Calabi-Yau algebra if and only if $A$ is an $\mathbb{N}$-graded AS-regular algebra (see Definition \ref{definition of GAS} and \ref{definition of CY of dimension d}) and $A/J$ is a separable algebra, where $J$ is the graded Jacobson radical of $A$ \cite{RR}. 
By the ``twisted" Calabi-Yau algebras $A$, we mean that the Nakayama bimodule $U$ of $A$ is invertible but not necessarily free as a left or right $A$-module; by the ``skew" Calabi-Yau algebras $A$, we mean that $U$ is free as a left or right $A$-module. Such differences vanish for connected graded algebras.

By employing canonical $\mathbb{N}$-graded AS-regular algebras, we investigate the conditions under which a graded twisted Calabi-Yau algebra becomes a graded or ungraded skew Calabi-Yau, as outlined in Theorem \ref{when twisted CY is skew CY}. Notably, some graded twisted Calabi-Yau algebras may be ungraded skew Calabi-Yau, but not be graded skew Calabi-Yau (see Example \ref{example: graded twisted CY is not graded skew CY}).

Recall that a Hom-finite (i.e. all Hom spaces are  finite-dimensional) $k$-linear triangulated category is called $d$-Calabi-Yau if the $d$-th shift functor is a Serre functor \cite{K}. It is well-known that for any noetherian Calabi-Yau algebra $A$, the bounded derived category of finite-dimensional $A$-modules is a Calabi-Yau category \cite{Kel}. Inspired by the work of Iyama and Reiten \cite{IR}, we study the converse of this result,
that is, if $\D^b(\fd A)$ is a Calabi-Yau category, does this imply that $A$ itself is Calabi-Yau?

A noetherian commonly graded algebra $A$ is called twisted derived-Calabi-Yau (twisted derived-CY for short) of dimension $d$, if the bounded derived category of finite-dimensional graded $A$-modules is a ``twisted" $d$-Calabi-Yau category, which means this triangulated category admits a Serre functor $U\otimes_A-[d]$, where $U$ is an invertible graded $(A,A)$-bimodule (see Definition \ref{definition of dCY algebras}). 
The following is proved in Theorem \ref{GAS-regular and d-CY}.

\begin{theorem}
Suppose $A$ is a noetherian commonly graded algebra. Then
\begin{itemize}
\item[(1)] $A$ is a commonly graded AS-regular algebra of dimension $d$ if and only if $A$ is a twisted derived-CY algebra of dimension $d$.
\item[(2)] If, further, $A$ is $\mathbb{N}$-graded, then $A$ is a twisted (resp. skew) Calabi-Yau algebra of dimension $d$ if and only if $A$ is a twisted (resp. skew) derived-CY algebra of dimension $d$ with $S=A/J$ being a separable algebra. 
\end{itemize}
\end{theorem}

In the last part of the article, we study incidentally the Auslander-Buchsbaum formula, the Bass theorem and the No-Hole theorem for commonly graded algebras, which are important tools in commutative ring theory. 
The investigation of these three formulae within the context of connected graded algebras is addressed in \cite{Jo1,Jo2}. Roughly speaking, for any noetherian connected graded algebra $A$, if $A$ satisfies $\chi$-condition, then the Auslander-Buchsbaum formula holds; if $A$ has a balanced dualizing complex, then the Bass theorem and the No-hole theorem holds.
Nevertheless, these results do not generally apply to non-connected graded algebras. The underlying issue is that the proofs in the connected graded case depend on the principle that the tensor product $M\otimes_A N$ is non-zero for any two non-zero bounded-below graded modules $M$ and $N$ because $k$ is the only graded simple $A$-module up to shifts and isomorphisms. In contrast, when $A$ is non-connected, it is possible for $M\otimes_A N$ to be zero.

For a noetherian commonly graded algebra $A$ with a balanced dualizing complex, we demonstrate that the Auslander-Buchsbaum formula holds if and only if the Bass theorem holds, which is also equivalent to a set of conditions being satisfied (see Theorem \ref{AB formula and Bass thm for locally finite alg}).
In cases where $A$ is a balanced Cohen-Macaulay ring, meaning that its balanced dualizing complex concentrates in a single degree (Definition \ref{def of balanced CM}), more equivalent conditions are provided (see Theorem \ref{AB formula and Bass thm for balanced CM alg}). These conditions extend those presented in \cite{Mo}.
For simplicity, we will refer to these equivalent conditions as $\mathbf{(C)}$, and when considering graded right modules, we will denote these conditions by $\mathbf{(C^o)}$. If $A$ is  balanced Cohen-Macaulay and meets the requirements $\mathbf{(C)}$ and $\mathbf{(C^o)}$, then the No-Hole theorem is applicable. The following theorem lists a portion of $\mathbf{(C)}$, while the complete conditions $\mathbf{(C)}$ are detailed in Section \ref{Auslander-Buchsbaum Formula and Bass Theorem}.

\begin{theorem}[Theorem \ref{AB formula and Bass thm for locally finite alg}, \ref{AB formula and Bass thm for balanced CM alg}, \ref{no-hole thm}]
    Let $A$ be a noetherian commonly graded algebra with a balanced dualizing complex $R$. Suppose $\depth_AA=d$. Then the following are equivalent.
    \begin{itemize}
        \item [C1] (Bass theorem) If $0 \neq X\in D^b_{\fg}(\Gr A)$ and $\idim_AX<\infty$, then $\idim_AX=d+\sup X$.
        \item [C2] (Auslander-Buchsbaum Formula) If $0 \neq Y\in D^b_{\fg}(\Gr A^o)$ and $\pdim_{A^o}Y<\infty$, then $\pdim_{A^o}Y+\depth_{A^o} Y=d$.
        \item [C4] $\gExt_A^d(M,A)\neq 0$ for every graded simple $A$-module $M$.
    \end{itemize}

    If $A$ is balanced Cohen-Macaulay and $R=\Omega[d]$ for some graded $(A,A)$-bimodule $\Omega$, then they are equivalent to each of the following conditions.
    \begin{itemize}
        \item [C12] For any $M\in\gr A$ with $\idim_AM<\infty$, $M$ has a resolution
        $$0\to W_m\to\cdots\to W_0\to M\to 0$$
        where $W_i\in \add_A\Omega$ and $m=d-\depth_AM$.
        \item [C16] Every maximal Cohen-Macaulay $A^o$-module with finite projective dimension is projective.
        \item [C17] Every maximal Cohen-Macaulay $A$-module with finite injective dimension is in $\add_A\Omega$.
    \end{itemize}
    In particular if $A$ is balanced Cohen-Macaulay and $A$ satisfies $\mathbf{(C)}$ and $\mathbf{(C^o)}$, then No-Hole Theorem holds:
    $$\gExt_A^i(S,M)\neq 0 \text{ if and only if } \depth_A M\leqslant i\leqslant \idim_AM$$
    for any $M\in\gr A$. 
\end{theorem}

It should be noted that condition C$4$ is inherently satisfied by any connected graded algebra with $\depth d$. Consequently, Theorem  \ref{AB formula and Bass thm for locally finite alg} generalizes the results established in \cite{Jo1, Jo2}. Furthermore, as per its definition, any noetherian commonly graded AS-Gorenstein algebra satisfies condition C$4$. Thus, the Auslander-Buchsbaum formula, the Bass theorem, and the No-Hole theorem hold in this context.

The article is organized as follows. In Section \ref{Preliminaries}, we introduce some notations and give some preliminaries. In Section \ref{Commonly graded Artin-Schelter Gorenstein Algebras}, we give some characterizations of Ext-finite commonly graded AS-Gorenstein algebras. In Section \ref{Existence of Balanced Dualizing Complex}, we prove the existence of the balanced dualizing complex over noetherian commonly graded AS-Gorenstein algebras. In particular we study the structure of the minimal graded injective resolution of $A$, which is important for the later sections. In Section \ref{Canonical commonly graded AS-Gorenstein Algebras}, we show the relations and the difference of kinds of generalizations of AS-Gorenstein algebra in literature. In Section \ref{Graded and Ungraded Calabi-Yau Algebras}, we study when $\mathbb{N}$-graded AS-regular algebras are graded and/or ungraded skew Calabi-Yau algebras. In Section \ref{commonly graded AS-Regular Algebras and Twisted Derived-Calabi-Yau Algebras}, we introduce the concept of twisted derived-CY algebras. Then we discuss the connections between derived-CY algebras and commonly graded AS-regular algebras. We show that noetherian commonly graded AS-regular algebras are the same as twisted derived-CY algebras. In Section \ref{Auslander-Buchsbaum Formula and Bass Theorem}, we investigate the Auslander-Buchsbaum formula, the Bass theorem and the No-Hole theorem for noetherian commonly graded algebras.

\section{Preliminaries}\label{Preliminaries}
\subsection{Notations and Conventions}\label{Notations and Conventions}
Let $k$ be a field. A $\mathbb{Z}$-graded $k$-space $X$ is called locally finite if $X_i$ is finite-dimensional for all $i\in\mathbb{Z}$; $X$ is called bounded-below (resp. bounded-above) if there is an integer $n$ such that $X_i=0$ for all $i<n$ (resp. $i>n$).  Let
$$b_l(X)=\inf\{i\mid X_i\neq 0\}\text{ and }b_u(X)=\sup\{i\mid X_i\neq 0\}.$$

A $k$-algebra $A$ with a $\mathbb{Z}$-graded vector space decomposition $A=\cdots\oplus A_{-1}\oplus A_0\oplus A_1\oplus A_2\oplus\cdots$ is called a $\mathbb{Z}$-graded algebra if $A_iA_j\subseteq A_{i+j}$ for all subspaces $A_i$ and $A_j$. If $A_i=0$ for all $i<0$, then $A$ is called an $\mathbb{N}$-graded algebra. An $\mathbb{N}$-graded algebra $A$  is called connected graded if $A_0=k$. A $\mathbb{Z}$-graded algebra $A$ is called commonly graded if $A$ is locally finite and bounded-below.

Let $A^e=A\otimes A^o$  where $A^o$ is the opposite ring of $A$. When we say a module, it means a left module. Usually, a right  $A$-module is viewed as a left $A^o$-module; and $(A,B)$-bimodules are viewed as $A\otimes B^o$-modules for $k$-algebras $A$ and $B$. 

For any graded $A$-modules $M, N$, a graded $A$-module morphism $f$ of degree $0$ is an  $A$-morphism $f: M \to N$ such that $f(M_i) \subseteq N_i$ for all $i \in \mathbb{Z}$.
The category of graded $A$-modules with morphisms of degree $0$ is denoted by $\Gr A$. The full subcategory of $\Gr A$ consisting of finitely generated (resp. finite-dimensional) graded $A$-modules is denoted by $\gr A$  (resp. $\fd A$). 

Let $M(n)$ be the $n$-th shift of a graded module $M$ with $M(n)_i =M_{n+i}$.
For any $M, N\in \Gr A$, 
$$\gHom_A(M,N):=\oplus_{n\in \mathbb{Z}}\Hom_{\Gr A}(M,N(n)).$$ 
Let $D(-):=\gHom_k(-,k)$ be the Matlis dual functor. 

Let $X[n]$ be the $n$-th shift of a complex $X$ with $X[n]^i=X^{n+i}$.
For any additive category or abelian category $\mathcal{C}$, let $\K(\mathcal{C})$, $\K^+(\mathcal{C}), \K^-(\mathcal{C})$ and $\K^b(\mathcal{C})$ denote the  homotopy category consisting of unbounded, left-bounded, right-bounded and bounded complexes of $\mathcal{C}$ respectively, and let $\D^*(\mathcal{C})$, where $*\in \{\emptyset, +, -, b\}$, denote the corresponding derived category of $\mathcal{C}$. 

For any $X,Y\in \D(\Gr A)$, let
$$\Ext_{\Gr A}^i(X,Y):=H^i(R\Hom_{\Gr A}(X,Y)) \text{ and } \gExt_A^i(X,Y):=H^i(R\gHom_A(X,Y)).$$

When $A$ is left noetherian, let $\D^*_{\fg}(\Gr A)$ be the full subcategory of $\D^*(\Gr A)$ consisting of the complexes with finitely generated cohomologies.

\subsection{Commonly Graded Algebras}\label{Graded Algebras}
In this subsection we study the basic properties of commonly graded algebras. 
Let $A$ be a commonly graded algebra and $J$ be its graded Jacobson radical.

An essential difference between commonly graded algebras and $\mathbb{N}$-graded algebras is that $M_{\geqslant i}:=\oplus_{n \geqslant i}M_n$ is not a submodule of $M$ in general for a graded $A$-module $M$. 
Fortunately, instead of $M_{\geqslant i}$, $J^iM$ works well under some finiteness conditions (see Lemma \ref{J^nM and M>n cofinal}).
Many proofs of the results concerning locally finite $\mathbb{N}$-graded (connected graded) algebras still work for commonly graded algebras without essential change but by replacing $M_{\geqslant i}$ with $J^iM$. So in the rest of this article, we will cite the results directly which are proved for $\mathbb{N}$-graded or connected graded algebras if their proofs can be modified in an aforementioned way.

\begin{lemma}\cite[Lemma 3.2]{CKWZ2}\label{A-J finite-dim}
\begin{itemize}
\item[(1)] If $L$ is a maximal graded left ideal of $A$, then $L\supseteq A_{> -b_l(A)}$. As a consequence, $J\supseteq A_{>-b_l(A)}$.
\item[(2)] For any maximal graded left ideal $L$  of $A$,
$$b_l(A) \leqslant b_l(A/L)\leqslant b_u(A/L)\leqslant -b_l(A).$$
As a consequence, $A/L$ is finite-dimensional.
\item[(3)] $b_l(A) \leqslant b_l(A/J)\leqslant b_u(A/J)\leqslant -b_l(A)$,  so $A/J$ is finite-dimensional and is graded semisimple, i.e. $A/J$ is a direct sum of graded simple modules.
\end{itemize}
\end{lemma}

 We fix some notations to study finite-dimensional semisimple graded algebras. Let $D$ be an algebra, and 
$n\in \mathbb{N}^+$, $\gamma_1,\cdots,\gamma_n\in \mathbb{Z}$. 
Regard $D$ as a $\mathbb{Z}$-graded algebra concentrated in degree $0$. Then
$$\begin{pmatrix}
D& D(\gamma_1-\gamma_2)&\cdots& D(\gamma_1-\gamma_n)\\
D(\gamma_2-\gamma_1)& D& \cdots & D(\gamma_2-\gamma_n)\\
\vdots &\vdots &\ddots& \vdots\\
D(\gamma_n-\gamma_1)& D(\gamma_n-\gamma_2)&\cdots &D
\end{pmatrix}$$
is a $\mathbb{Z}$-graded algebra with the usual matrix multiplication. This graded algebra is denoted by $M_n(D)(\bar{\gamma})$, where $\bar{\gamma}=(\gamma_1,\cdots,\gamma_n)\in\mathbb{Z}^n$.

\begin{lemma}\label{A/J is graded semisimple and its structure}
    Let $S=A/J$. Then 
    \begin{itemize}
        \item [(1)] $S$ is a graded symmetric algebra, that is, $D(S)\cong S$ as graded $(S,S)$-bimodules.
        \item [(2)] $S\cong \mathop{\oplus}\limits_{i=1}^n M_{r_i}(D_i)(\bar{\gamma^i})$ for some integer $n$, division algebras $D_i$,  and $\bar{\gamma^i}\in\mathbb{Z}^{r_i}\, (1 \leqslant i \leqslant n)$.
    \end{itemize}
\end{lemma}
\begin{proof}
(1) It follows from Lemma \ref{A-J finite-dim} (3) and \cite[Theorem A]{DNN} (or \cite[16F]{Lam2}).

    (2) It follows from Lemma \ref{A-J finite-dim} (3) and \cite[Theorem 2.10.10, Exercise 5.6.2]{NO2}. 
\end{proof}

\begin{definition}
The socle of a graded $A$-module $M$ is the largest graded semisimple submodule of $M$, which is denoted by $\soc M$.
\end{definition}

Obviously, $\soc M\cong \gHom_A(S,M)$ as graded $A$-modules (or $S$-modules).

\begin{lemma}\label{soc}
Let $M$ be a locally finite graded $A$-module. Then $D(M/JM)=\soc(D(M))$ and $D(JM) \cong D(M)/\soc D(M)$, where $D$ is the Matlis dual functor.
\end{lemma}
\begin{proof}
The exact sequence $0\to JM\to M\to M/JM\to 0$ induces the exact sequence $0\to D(M/JM)\to D(M)\to D(JM)\to 0$. Since $D(M/JM)$ is a semisimple graded $S^o$-module, $D(M/JM)\subseteq \soc D(M)$. 

The injection $\soc D(M)\to D(M)$ induces the surjection $M\to D(\soc D(M))$. Since $D(\soc D(M))$ is semisimple, $M\to D(\soc D(M))$ factors through $M/JM$. Hence there is a surjection $M/JM\to D(\soc D(M))$ whose dual $\soc D(M)\to D(M/JM)$ is an injection. It follows that $\soc D(M)\subseteq D(M/JM)$. Hence
\[\soc D(M) = D(M/JM) \, \textrm{ and } \, D(JM) \cong D(M)/D(M/JM)=D(M)/\soc D(M).
\qedhere
\]
\end{proof}

\begin{lemma}\label{J^i is contained in A>j} For any $i>0$, there is some integer $n_i$ such that $J^{n_i}\subseteq A_{\geqslant i}$. As a consequence, $\bigcap_nJ^n=0$.
\end{lemma}
\begin{proof}
Let $I=AA_{\geqslant 1-2b_l(A)}A$, which is an ideal of $A$ contained in $A_{\geqslant 1}$. Let $B=A/I$. Then $B$ is finite-dimensional. Hence $J(B)^k=0$ for some integer $k$, which means $J^k\subseteq I\subseteq A_{\geqslant 1}$. Therefore $J^{ik}\subseteq A_{\geqslant i}$ for all $i$. 
\end{proof}

The following is a version of Nakayama's Lemma for commonly graded algebras \cite[Proposition 3.3(2)]{CKWZ2} \cite[Corollary 2.9.2]{NO2}.
\begin{lemma}\label{Nakayama's Lemma}
If $M$ is a bounded-below graded $A$-module such that $MJ=M$, then $M=0$.
\end{lemma}
\begin{proof}
For any fixed $i > 0$, by Lemma \ref{J^i is contained in A>j}, there is some $n$ such that  $J^n\subseteq A_{\geqslant i}$. Then $M=MJ=MJ^n\subseteq MA_{\geqslant i}$. If $M \neq 0$, it follows from $M$ being bounded-below that  $b_l(M)\geqslant b_l(M) + i$, which is a contradiction.
\end{proof}

The following lemma is a version of \cite[Lemma 2.3]{RR} for commonly graded algebras.

\begin{lemma}\label{finitely generated algebra}
Suppose $A$ is a commonly graded algebra. Let $S=A/J$. Then the following are equivalent.
\begin{itemize}
    \item [(1)] $J/J^2$ is finite-dimensional.
    \item [(2)] $S$ is finitely presented as a graded $A$-module.
    \item [(3)] $S$ is finitely presented as a graded $A^o$-module.
    \item [(4)] $A$ is finitely generated as a $k$-algebra.
    \item [(5)] $\{J^n\}_n$ and $\{A_{\geqslant n}\}_n$ are cofinal.
\end{itemize}
\end{lemma}
\begin{proof}
(1) $\Leftrightarrow$ (2) $\Leftrightarrow$ (3) By Lemma \ref{Nakayama's Lemma}, $J/J^2$ is finite-dimensional if and only if $J$ is finitely generated as a left (right) graded $A$-module.
This is equivalent to that ${}_AS$ (or equivalently, $S_A$) is a finitely presented graded module. Hence (1)-(3) are equivalent. 

(1) $\Rightarrow$ (4)  By Lemma \ref{A-J finite-dim}, $A/J$ is finite-dimensional. 
Let $V$ be a finite-dimensional subspace of $J$ such that $V+J^2=J$, and $X$ be a finite-dimensional space of $A$ such that $X+J=A$. By induction, $V^n+J^{n+1}=J^n$. Then $A=X+V+V^2+\cdots+V^n+J^{n+1}$ for all $n$. By lemma \ref{J^i is contained in A>j}, for any $i > 0$ there is an integer $n_i$ such that $J^{n_i}\subseteq A_{\geqslant i}$. Thus $A_{i-1}\cap J^{n_i}=0$. This implies $A_{i-1}\subseteq X+V+V^2+\cdots+V^{n_i-1}$. Therefore $A$ is generated as a $k$-algebra by $A_{\leqslant 0}$, $V$ and $X$.

(4) $\Rightarrow$ (1)  By Lemma \ref{A-J finite-dim}, there is an integer $n\in\mathbb{N}$ such that $A_{\geqslant n}\subseteq J$.  
We claim that there is some $n'$ such that $A_{\geqslant n'}\subseteq A_{\geqslant n}A_{\geqslant n}\subseteq J^2$. Then it follows that $J/J^2$ is finite-dimensional.
Let $\{a_1,\cdots,a_t\}$ be a set of homogeneous generators of the graded algebra $A$, and $m=\max\{\deg a_1,\cdots,\deg a_t\}$.
Set $n'=2n+m$. Suppose $a=b_1\cdots b_s \in A_{\geqslant n'}$, where $b_i\in\{a_1,\cdots,a_t\}$ for all $1\leqslant i\leqslant s$. Then for some $1\leqslant s'\leqslant s$, $n\leqslant \deg (b_1\cdots b_{s'})\leqslant n+m$. It follows from $a\in A_{\geqslant n'}$ that $\deg(b_{s'+1}\cdots b_s)\geqslant n$. Hence $a\in A_{\geqslant n}A_{\geqslant n}$ which means $A_{\geqslant n'}\subseteq A_{\geqslant n}A_{\geqslant n}$.

(4) $\Rightarrow$ (5) By Lemma \ref{J^i is contained in A>j}, we only need to prove that for any $i>0$ there is an integer $j$ such that $A_{\geqslant j}\subseteq J^i$.
With the notations as before, $A=X+V+V^2+\cdots+V^{i-1}+J^i$.
Since $V$ and $X$ are finite-dimensional, $X+V+V^2+\cdots+V^{i-1}$ are both bounded-below and bounded-above. Thus for $j\gg 0$, $A_{\geqslant j}\subseteq J^i$.

(5) $\Rightarrow$ (1) By assumption, there is some integer $n$ such that $A_{\geqslant n}\subseteq J^2\subseteq J$.  Since $A$ is locally finite, $J/J^2$ is finite-dimensional.
\end{proof}

\begin{lemma}\label{J^nM and M>n cofinal}
If $A$ is a finitely generated $k$-algebra, then for any finitely generated graded $A$-module $M$, $\{M_{\geqslant n}\}_n$ and $\{MJ^n\}_n$ are cofinal.
\end{lemma}
\begin{proof}
By Lemma \ref{finitely generated algebra}, it suffices to show that $\{M_{\geqslant n}\}_n$  and  $\{MA_{\geqslant n}\}_n$ are cofinal.

Since $M$ is bounded-below, for any $n$, there is an integer $m$ such that $MA_{\geqslant m} \subseteq M_{\geqslant n}$.

Let $\{m_1,\cdots, m_s\}$ be a set of homogeneous generators of $M$ and  $r=\max\{r_1,\cdots, r_s\}$ where $r_i=\deg m_i$. Then, for any $n$ and  any homogeneous element  $m\in M$ with $\deg m \geqslant n+r$, $m=\Sigma m_ia_i$ for some $a_i \in A$ with $\deg a_i=\deg m-r_i$. So,  $m\in MA_{\geqslant n}$. Hence, $M_{\geqslant n+r}\subseteq MA_{\geqslant n}$.
\end{proof}

\subsection{Local Cohomology and Minimal Resolutions}
In this subsection, we introduce local cohomology theory for commonly graded algebras, and discuss homological dimensions of complexes in terms  of minimal graded projective  or injective resolutions. 

Let $A$ be a commonly graded algebra. For any graded $A$-module $M$, let 
$$\Gamma_A(M)=\{m\in M\mid J^i\cdot m=0, i\gg 0\}.$$
Clearly, $\Gamma_A(M)\cong  \underrightarrow{\Lim}\gHom_A(A/J^i,M)$.

If $A$ is a finitely generated $k$-algebra, then  $\Gamma_A(M)=\{m\in M\mid A_{\geqslant i}\cdot m=0, i\gg 0\}$ by Lemma \ref{finitely generated algebra}.

Let $R\Gamma_A$ be the right derived functor of the left exact functor $\Gamma_A: \Gr A \to \Gr A$. The $n$-th local cohomology of $M$ is defined to be $R^n\Gamma_A(M)$. In fact, for $n\geqslant 0$,
$$R^n\Gamma_A(M)\cong \underrightarrow{\Lim}\gExt_A^n(A/J^i,M).$$ 
We say that $\Gamma_A$ has finite cohomological dimension, or $A$ has finite local cohomological dimension, if there exists an integer $n$ such that for all $M\in \Gr A$ and $i>n$, $R^i\Gamma_A(M)=0$. If further there is some $M\in \Gr A$, such that $R^n\Gamma_A(M)\neq 0$, then we say $n$ is the cohomological dimension of $\Gamma_A$.

 If $\Gamma_A(M)=M$ (resp. $\Gamma_A(M)=0$) then $M$ is called torsion (resp. torsion-free). Any $x\in \Gamma_A(M)$ is called a torsion element of $M$.
 Note that a finitely generated graded $A$-module is torsion if and only if it is a finite-dimensional $A$-module as $A$ is assumed to be locally finite. 

$\Gamma_{A^o}$ and $R^n\Gamma_{A^o}$ are defined similarly.

\begin{lemma}\label{torsion elements and soc M}
For any $M\in \Gr A$, $\soc M\neq 0$ if and only if $M$ contains non-zero torsion elements.
\end{lemma}

For any $X\in \D(\Gr A)$,
the injective dimension of $X$ is defined by
$$\idim_A(X):=\sup\{i\mid \gExt_A^i(M,X)\neq 0 \text{ for some } M\in \Gr A\}.$$
When $X\in \D^+(\Gr A)$, by Baer's criterion,
$$\idim_A(X)=\sup\{i\mid \gExt_A^i(M,X)\neq 0 \text{ for some } M\in \gr A\}.$$
The projective dimension of $X$ is defined by
$$\pdim_A(X):=\sup\{i\mid \gExt_A^i(X,M)\neq 0 \text{ for some } M\in \Gr A\}.$$

 Next lemma shows the existence of minimal graded injective (projective) resolutions for left (right) bounded complexes over commonly graded algebras.

\begin{lemma}\label{minimal injective and projective complex}
    Let $A$ be a commonly graded algebra.
    \begin{itemize}
        \item [(1)]  If $X\in \D^+(\Gr A)$, then $X$ is quasi-isomorphic to a graded injective complex 
        $$\cdots \longrightarrow  I^{n-1}\xrightarrow[]{\partial^{n-1}} I^{n}\xrightarrow[]{\partial^n} I^{n+1} \longrightarrow  \cdots$$ such that $\Ker \partial^i$ is an essential submodule of $I^i$ for any $i$.
        \item[(2)] If $X\in \D^-(\Gr A)$ such that all modules in $X$ are bounded-below, then $X$ is quasi-isomorphic to a graded projective complex 
        $$\cdots \longrightarrow   P^{n-1}\xrightarrow[]{\partial^{n-1}} P^{n}\xrightarrow[]{\partial^n} P^{n+1} \longrightarrow  \cdots$$ 
        such that $\im \partial^i\subseteq JP^{i+1}$ for any $i$.
    \end{itemize}
\end{lemma}
\begin{proof}
    (1) The proof is similar to \cite[Lemma 4.2]{Y}.

    (2) By Nakayama's Lemma \ref{Nakayama's Lemma}, any bounded-below graded $A$-module admits a projective cover.
    Then it follows from the dual version of \cite[Lemma 4.2]{Y} that (2) holds.
\end{proof}

The graded injective (resp. projective) complexes satisfying the condition in Lemma \ref{minimal injective and projective complex} are called minimal. For any minimal graded injective complex $I^\bullet$ (resp. minimal graded projective complex $P^\bullet$), the differentials of $\gHom_A(A/J,I^\bullet)$ (resp. $\gHom_A(P^\bullet,A/J)$ and $A/J\otimes_A P^\bullet$) are zero.

The following lemmas are well-known for connected graded algebras, say, see \cite{Jo2}.
\begin{lemma}\label{pdim and idim and resolution}
    Let $A$ be a commonly graded algebra.
    \begin{itemize}
        \item [(1)] For any  $X\in \D^+(\Gr A)$,
        \begin{align*}
            \idim_A(X)&=\inf\limits_{I^\bullet}\{\sup\{i\mid I^i\neq 0, X \text{ is quasi-isomorphic to } I^\bullet\}\}\\
            &=\sup\{i\mid I^i\neq 0, I^\bullet \text{ is a minimal graded injective resolution of } X\}.
        \end{align*}
        \item [(2)] For any $X\in \D^-(\Gr A)$ such that all modules in $X$ are bounded-below,
        \begin{align*}
            \pdim_A(X)&=\inf\limits_{P^\bullet}\{-\inf\{i\mid P^i\neq 0, P^\bullet \text{ is quasi-isomorphic to } X\}\}\\
            &=-\inf\{i\mid P^i\neq 0, P^\bullet \text{ is a minimal graded projective resolution of } X\}.
        \end{align*}
    \end{itemize}
\end{lemma}

The depth of a complex $X\in \D^+(\Gr A)$ is defined by
$$\depth_A X:=\min\{i\mid \gExt_A^i(A/J,X)\neq 0\}.$$

\begin{lemma}\label{depth and local dimension}
Let $A$ be a commonly graded algebra and $S=A/J$. If $A$ is a finitely generated $k$-algebra, then, for any $X\in \D^+(\Gr A)$,
$$\depth_A X=\min\{i\mid R^i\Gamma_A(X)\neq 0\}.$$
\end{lemma}
\begin{proof}
Let $d_1=\min\{i\mid \gExt_A^i(S,X)\neq 0\}$ and $d_2=\min\{i\mid R^i\Gamma_A(X)\neq 0\}$. Let $E^\bullet$ be a minimal graded injective resolution of $X$. Then $\gExt^i_A(S,X)=\gHom_A(S,E^i)\cong \soc E^i$. Lemma \ref{torsion elements and soc M} implies that for any $i<d_1$, $E^i$ is torsion-free. 
Hence $R^i\Gamma_A(X)=H^i(\Gamma_A(E^\bullet))=0$ for all $i<d_1$ which yields that $d_2\geqslant d_1$.

Since $\gExt^{d_1-1}_A(S,X)=0$, $\gExt^{d_1-1}_A(M,X)=0$ for any graded simple module $M$. 
An induction on the length of a finite-dimensional graded $A$-module $Y$ implies that $\gExt^{d_1-1}_A(Y,X)=0$. 

Note that $J^j/J^i$ and $A/J^i$ are finite-dimensional for any $i>j$ by Lemma \ref{finitely generated algebra}.
The exact sequence $0\to J^j/J^i\to A/J^i\to A/J^j\to 0$ induces a long exact sequence
$$\cdots\to \gExt^{d_1-1}_A(J^j/J^i,X)\to \gExt_A^{d_1}(A/J^j,X)\to \gExt_A^{d_1}(A/J^i,X)\to\cdots.$$
So $\gExt_A^{d_1}(A/J^j,X)\to \gExt_A^{d_1}(A/J^i,X)$ is injective and $R^{d_1}\Gamma_A(X)$ is the union of $\gExt_A^{d_1}(A/J^i,X)$.  It follows from $\gExt_A^{d_1}(S,X)\neq 0$ that $R^{d_1}\Gamma_A(X) \neq 0$. Hence $d_1\geqslant d_2$ and in conclusion $d_1=d_2$.
\end{proof}

\subsection{Local Duality and Balanced Dualizing Complex}

The results in \cite[Section 3-8]{VdB} still hold for commonly graded algebras with similar proofs by changing $A_{\geqslant i}$ to $J^i$, and $M_{\geqslant i}$ to $J^iM$. 
In particular, the existence theorem of balanced dualizing complexes over connected graded algebras  \cite[Theorem 6.3]{VdB} still holds for commonly graded algebras.  We recall them here. 

\begin{definition}
A commonly graded algebra $A$ is called Ext-finite, if $\gExt_A^i(S,S)$ is finite-dimensional for all $i$, where $S=A/J$.
\end{definition}

The notion of ``Ext-finite" is left-right symmetric, as showed in the following.

\begin{lemma}\label{Ext-finite is left-right symmetric}
Let $A$ be a commonly graded algebra and $S=A/J$. Then the following are equivalent.
\begin{itemize}
\item [(1)] $A$ is Ext-finite.
\item [(2)] $\gExt_{A^o}^i(S,S)$ is finite-dimensional for all $i$.
\item [(3)] $\Tor_i^A(S,S)$ is finite-dimensional for all $i$.
\item [(4)] ${}_AS$ has a finitely generated graded projective resolution.
\item [(5)] $S_A$ has a finitely generated graded projective resolution.
\end{itemize}
\end{lemma}
\begin{proof}
(1) $\Leftrightarrow$ (3) $\Leftrightarrow$ (4). Let $P_\bullet\to S\to 0$ be a minimal graded projective resolution of ${}_AS$. Then $\gExt_A^i(S,S)\cong \gHom_A(P_i,S)$ and $\Tor_i^A(S,S)=S\otimes_A P_i$. Thus $A$ is Ext-finite if and only if $P_i$ is finitely generated for all $i$; if and only if $\Tor_i^A(S,S)$ is finite-dimensional for all $i$.

(2) $\Leftrightarrow$ (3) $\Leftrightarrow$ (5). Similar to the proof of (1) $\Leftrightarrow$ (3) $\Leftrightarrow$ (4).
\end{proof}

\begin{corollary}\label{Ext-finite implies f.d. module has f.g. proj. reso.}
Let $A$ be a commonly graded algebra. 
\begin{itemize}
    \item [(1)] If $A$ is left or right noetherian, then $A$ is Ext-finite.
    \item [(2)] If $A$ is Ext-finite, then every finite-dimensional graded $A$-module (resp. $A^o$-module) has a finitely generated graded projective resolution.
    \item [(3)] If $A$ is Ext-finite, then $A$ is a finitely generated $k$-algebra.
\end{itemize} 
\end{corollary}

\begin{theorem}[Local Duality]\label{local duality}\cite[Theorem 5.1]{VdB}
Suppose $A$ is an Ext-finite commonly graded algebra and $\Gamma_A$ has finite cohomological dimension. Then
\begin{itemize}
\item[(1)] The injective dimension of $D(R\Gamma_A(A))$ is $0$ in $\D(\Gr A)$.
\item[(2)] For any graded algebra $B$, and $M\in \D(\Gr A\otimes B^o)$,
$$D(R\Gamma_A(M))\cong R\gHom_A(M,D(R\Gamma_A(A)))$$
in $\D(\Gr B\otimes A^o)$.
\end{itemize}
\end{theorem}

The $\chi$-condition is defined in \cite{AZ} to establish cohomology theory for noncommutative projective schemes in the $\mathbb{N}$-graded case. 
It has a strong relation with the existence of balanced dualizing complexes. Here is the definition of the $\chi$-condition for commonly graded algebras.

\begin{definition}
A commonly graded algebra $A$ is called satisfying the $\chi$-condition if for any $M\in \gr A$ and $i\geqslant 0$, $\gExt_A^i(A/J,M)$ is bounded-above.
\end{definition}

If $A$ is right noetherian, then $A$ satisfies the $\chi$-condition 
if and only if $\gExt_A^i(A/J,M)$ is finite-dimensional for any $M\in \gr A$ and $i\geqslant 0$.

\begin{lemma}\label{depth A and A^o}
    Let $A$ be a noetherian commonly graded algebra. Suppose both $A$ and $A^o$ satisfy the $\chi$-condition. Then $\depth_AA=\depth_{A^o}A$.
\end{lemma}
\begin{proof}
    It follows from \cite[Corollay 4.8]{VdB} and Lemma \ref{depth and local dimension}.
\end{proof}

Now we recall the definition of dualizing complexes \cite{Y}.
\begin{definition}\label{Dualizing complex}
Let $A$ be a noetherian commonly graded algebra. A complex $R\in \D^b(\Gr A^e)$ is called a dualizing complex of $A$, if it satisfies the following conditions:
\begin{itemize}
\item[(1)] $R$ has finite injective dimension over $A$ and $A^o$ respectively.
\item[(2)] The cohomologies of $R$  are finitely generated both as $A$-module and $A^o$-module.
\item[(3)] The natural morphisms $\Phi: A\to R\gHom_A(R,R)$ and $\Phi^o:A\to R\gHom_{A^o}(R,R)$ are isomorphisms in $\D(\Gr A^e)$.
\end{itemize}
If moreover, $R\Gamma_A(R)\cong D(A)$ and $R\Gamma_{A^o}(R)\cong D(A)$ in $\D(\Gr A^e)$, then $R$ is called a balanced dualizing complex of $A$.
\end{definition}

Next lemma justifies the name ``dualizing complex".

\begin{lemma}\label{dualizing complex induces a duality}
    Let $A$ be a noetherian commonly graded algebra with a dualizing complex $R$. 
    \begin{itemize}
        \item [(1)] The functors $R\gHom_A(-,R)$ and $R\gHom_{A^o}(-,R)$ give a duality between $\D_{\fg}(\Gr A)$ and $\D_{\fg}(\Gr A^o)$, restricting to a duality between $\D^b_{\fg}(\Gr A)$ and $\D^b_{\fg}(\Gr A^o)$ (resp. $\D^-_{\fg}(\Gr A)$ and $\D^+_{\fg}(\Gr A^o)$, $\D^+_{\fg}(\Gr A)$ and $\D^-_{\fg}(\Gr A^o)$).
        \item [(2)] For any $X\in \D^b_{\fg}(\Gr A)$, if $\idim_AX<\infty$ (resp. $\pdim_AX<\infty$), then 
        $$\pdim_{A^o}R\gHom_A(X,R)<\infty\, (\textrm{ resp. }\, \idim_{A^o}R\gHom_A(X,R)<\infty).$$ 
        Similar results  hold for any $Y \in \D^b_{\fg}(\Gr A^o)$. 
    \end{itemize}
\end{lemma}
\begin{proof}
    (1) See \cite[Proposition 1.3]{YZ}.

    (2) It follows from (1).
\end{proof}

The following is the commonly graded version of \cite[Theorem 6.3]{VdB}.

\begin{theorem}\label{existence of balanced dualizing complex}
Let $A$ be a noetherian commonly graded algebra. Then $A$ admits a balanced dualizing complex if and only if $A$ satisfies the following two conditions:
\begin{itemize}
\item[(1)] both $A$ and $A^o$ satisfy the $\chi$-condition;
\item[(2)] both $\Gamma_A$ and $\Gamma_{A^o}$ have finite cohomological dimension.
\end{itemize}
If $A$ admits a balanced dualizing complex $R$, then $R \cong D(R\Gamma_A(A))\cong D(R\Gamma_{A^o}(A))$.
\end{theorem}

\subsection{Invertible Graded Bimodules}
The results in this subsection are known, whereas the proofs are given here because we cannot locate a suitable reference.

Let $A$ be a $k$-algebra. For any $A$-module $M$ and any  algebra endomorphism $\varphi$ of $A$, ${}^\varphi M$ is the $A$-module with the same $k$-space structure as $M$ but with a new $A$-action: $a \cdot m = \varphi(a)m$ for any $a\in A$ and $m\in {}^\varphi M$. If $M$ is an $(A,A)$-bimodule, $\varphi$ and $\psi$ are algebra endomorphisms of $A$, then ${}^\varphi M^\psi$ is an $(A,A)$-bimodule defined similarly. If $\varphi = \id$,  ${}^\varphi M^\psi$ is simply denoted as ${}^1 M^\psi$ or $ M^\psi$.

\begin{lemma}\label{bimod-summand-of-S}
Suppose $S=\mathop{\oplus}\limits_{i=1}^n M_{r_i}(D_i)$ is an artinian semisimple algebra where the $D_i$'s are division algebras. If $U = {}^1 S^\varphi$ as $(S,S)$-bimodules  
for some automorphism $\varphi$ of $S$, and $V$ is an $(S,S)$-bimodule direct summand of $U$, then there is a subset $I$ of $\{1,\cdots,n\}$ such that $$V \cong {}^1 (\mathop{\oplus}\limits_{i\in I}   M_{r_i}(D_i))^\varphi  = \mathop{\oplus}\limits_{i\in I}  {}^1 M_{r_i}(D_i)^\varphi.$$
\end{lemma}
\begin{proof}
Suppose $U=V\oplus W$. Then ${}^1U^{\varphi^{-1}}={}^1V^{\varphi^{-1}}\oplus {}^1W^{\varphi^{-1}}\cong S$. 
Thus ${}^1V^{\varphi^{-1}}$ can be regarded as an ideal of $S$, 
which implies that ${}^1V^{\varphi^{-1}}\cong \mathop{\oplus}\limits_{i\in I} M_{r_i}(D_i)$ for some subset $I$ of $\{1,\cdots,n\}$. Hence 
\[
V \cong {}^1 (\mathop{\oplus}\limits_{i\in I}  M_{r_i}(D_i))^\varphi = \mathop{\oplus}\limits_{i\in I}  {}^1 M_{r_i}(D_i)^\varphi.
\qedhere
\]
\end{proof}

\begin{lemma}\label{graded invertible S bimodule}
Suppose $S=\mathop{\oplus}\limits_{i=1}^n M_{r_i}(D_i)$ is an artinian semisimple algebra where  the $D_i$'s  are division algebras. Regard $S$ as an $\mathbb{N}$-graded algebra concentrated in degree $0$.
Let $U$ be an invertible graded $(S,S)$-bimodule which is isomorphic to $\mathop{\oplus}\limits_{i=1}^n {}^1 M_{r_i}(D_i)^\varphi$ as ungraded modules, where $\varphi$ is an automorphism of $S$. Then, as graded $(S,S)$-bimodules, $U\cong \mathop{\oplus}\limits_{i=1}^n {}^1 M_{r_i}(D_i)^\varphi(l_i)$  for some integers $l_i$.
\end{lemma}
\begin{proof}
Since $U$ is finitely generated, $U =\mathop{\oplus}_{l \in \mathbb{Z}} U_l$ is concentrated in finitely many degrees, that is, only finitely many $U_l \neq 0$. 
Note that each $U_{l} \neq 0$ is an $(S,S)$-bimodule direct summand of $U \cong \mathop{\oplus}\limits_{i=1}^n {}^1 M_{r_i}(D_i)^\varphi $.
By Lemma \ref{bimod-summand-of-S}, $U_{l}\cong  \mathop{\oplus}\limits_{j\in I} {}^1 M_{r_j}(D_j)^\varphi$ as ungraded $(S,S)$-bimodules for some subset $I$ of $\{1,\cdots,n\}$. 
It follows that $U\cong \mathop{\oplus}\limits_{i=1}^n {}^1 M_{r_i}(D_i)^\varphi(l_i)$ for some integers $l_i$ as graded $(S,S)$-bimodules.
\end{proof}

\begin{lemma}\label{invertible bimodule of basic semisimple algebra}
Suppose $S=D_1\oplus \cdots\oplus D_n$ where the $D_i$'s are division algebras. Then every invertible $(S,S)$-bimodule is isomorphic to $\mathop{\oplus}\limits_{i=1}^n {}^1 D_i^\varphi$, where $\varphi$ is an automorphism of $S$. If $S$ is considered as an $\mathbb{N}$-graded algebra concentrated in degree $0$, then every invertible  graded $(S,S)$-bimodule is isomorphic to $\mathop{\oplus}\limits_{i=1}^n\,{}^1 D_i ^\varphi(l_i)$ for some integers $l_i$.
\end{lemma} 
\begin{proof}
Let $U$ be an invertible $(S,S)$-bimodule. On one hand, since $U_S$ is a finitely generated projective generator, $U\cong D_1^{(s_1)}\oplus \cdots\oplus D_n^{(s_n)}$ as $S^o$-modules, where $D_i^{(s_i)}$ means the direct sum of $s_i$ copies of $D_i$. On the other hand, it follows from the ring isomorphism
$$S\cong \End_{S^o}(U)\cong \End_{S^o}(D_1^{(s_1)}\oplus \cdots\oplus D_n^{(s_n)})\cong M_{s_1}(D_1)\oplus\cdots\oplus M_{s_n}(D_n)$$
that $s_1=\cdots=s_n=1$. Therefore $U\cong S$ as $S^o$-modules. 

Similarly, $U\cong S$ as $S$-modules. Since $U$ is invertible, there is an automorphism $\varphi$ of $S$ such that $U\cong {}^1 S^\varphi\cong  \mathop{\oplus}\limits_{i=1}^n {}^1 D_i ^\varphi$. 
The last conclusion follows from Lemma \ref{graded invertible S bimodule}.
\end{proof}

\section{Commonly graded AS-Gorenstein Algebras}\label{Commonly graded Artin-Schelter Gorenstein Algebras}
In this section we study the commonly graded AS-Gorenstein algebras.

Let $A$ be a commonly graded algebra and $S=A/J$. By Lemma \ref{A/J is graded semisimple and its structure}, we may assume that $S \cong \mathop{\oplus}\limits_{i=1}^n M_{r_i}(D_i)(\bar{\gamma}^i)$ where the $D_i$'s are division algebras and $\bar{\gamma}^i\in\mathbb{Z}^{r_i}$. 
For any $1\leqslant i\leqslant n$, let 
$$S_i:=(D,D(\gamma^i_2-\gamma^i_1),\cdots,D(\gamma^i_{r_i}-\gamma^i_1))^T$$
which is a graded simple module of $M_{r_i}(D_i)(\bar{\gamma}^i)$. Then $\{S_1,\cdots,S_n\}$ represents all the isomorphic classes of graded simple $S$-modules,
and graded simple $A$-modules up to shifting. There is a  set of orthogonal primitive idempotents $e_1,\cdots,e_n$ of $A_0$ (it may happen $e_1+\cdots+e_n\neq 1$) such that $S_i\cong S\bar{e}_i$ where $\bar{e}_i$ is the image of $e_i$ in $A_0/{J(A_0)}$.  Let $Q_i=Ae_i$ and $Q'_i=e_iA$. Then clearly $Q_i$ is a graded projective cover of $S_i$, and $Q'_i$ is a graded projective cover of $S'_i=\bar{e}_iS$.

\begin{lemma}\label{D(Si) cong Si'}
$D(S_i) \cong S'_i$ for any $i$.
\end{lemma}
\begin{proof}
Since $D(S)\cong S$ as $(S,S)$-bimodules (see Lemma \ref{A/J is graded semisimple and its structure}), 
\[
    S'_i \cong \Hom_S(S_i, S) \cong  \Hom_S(S_i, D(S)) \cong \Hom_k(S_i, k) = D(S_i).
\qedhere
\]
\end{proof}

\begin{proposition}\label{injective hull of simple module}
$D(Q'_i)$ is the injective hull of $S_i$, $D(Q_i)$ is the injective hull of $S'_i$, and $D(A_A)$ is the injective hull of ${}_AS$. 
\end{proposition}
\begin{proof}
    It follows from \cite[Corollary 2.8.2]{NO2}.
\end{proof}

AS-Gorenstein algebras are generalized in several slightly different ways \cite{V,V2,MM,RRZ,HY,IKU}. The following definition is similar to the one given by Mart\'inez-Villa \cite{V,V2}, where the algebra $A$ considered is not required to have finite left and right injective dimension, but $A_0$ is a finite product of base field. 

\begin{definition}\label{definition of GAS}
A commonly graded algebra $A$ is called Artin-Schelter Gorenstein (for short, AS-Gorenstein) of dimension $d$ if the following conditions hold.
\begin{itemize}
\item[(1)] $A$ has  finite left and right injective dimension $d$.
\item[(2)] For every graded simple $A$-module $M$, $\gExt^i_A(M,A)=0$ if $i\neq d$; for every graded simple $A^o$-module $N$, $\gExt_{A^o}^i(N,A)=0$ if $i\neq d$.
\item[(3)] $\gExt_A^d(-,A)$ and $\gExt_{A^o}^d(-,A)$ give a bijection between the isomorphic classes of graded simple $A$-modules to the isomorphic classes of graded simple $A^o$-modules.
\end{itemize}
If moreover $A$ has finite global dimension $d$, then $A$ is called Artin-Schelter regular (for short, AS-regular).

If $A$ is $\mathbb{N}$-graded, then it is called an $\mathbb{N}$-graded AS-Gorenstein (regular) algebra.
\end{definition}

When $A$ is connected graded, this definition coincides with the classic one.
If $A$ is $\mathbb{N}$-graded, then $A$ is an $\mathbb{N}$-graded AS-regular algebra defined in Definition \ref{definition of GAS} if and only if that $A$ is a generalized AS-regular algebra defined in \cite{RR}, as the global dimension is equal to the graded global dimension for any $\mathbb{N}$-graded algebra (see \cite[Corollary \uppercase\expandafter{\romannumeral1.7.8}]{NO1}).

It is established that AS-regular algebras are Ext-finite, as demonstrated in \cite[Proposition 3.1]{SZ}. An analogous assertion is valid for commonly graded AS-regular algebras. The following proposition covers a more extensive scope, and its proof is similar to \cite[Proposition 3.4]{MM}.

\begin{proposition}\label{criterion for Ext-finite} 
Let $A$ be a commonly graded algebra. Suppose $M$ is a graded $A$-module with finite projective dimension $n$. If $\gExt_A^i(M,A)=0$ for any $i\neq n$ and $\gExt_A^n(M,A)$ is finite-dimensional, then $M$ has a finitely generated graded projective resolution. 
\end{proposition}

\begin{corollary}\label{GAS-regular is Ext-finite}
If $A$ is commonly graded AS-regular, then $A$ is Ext-finite.
\end{corollary}

The following is an example of commonly graded AS-Gorenstein (regular) algebras with nonzero negative part.

\begin{example}
    Let $A$ be an AS-Gorenstein (regular) algebra and $B=\gEnd_A(A\oplus A(1))$. Then by Proposition \ref{GAS-Gorenstein is preserved by Morita}, $B$ is a commonly graded AS-Gorenstein (regular) algebra. Obviously, $B$ has nonzero negative part. 
\end{example}

Suppose $A$ is a commonly graded AS-Gorenstein algebra of dimension $d$. With the notations as above, there is a permutation $\sigma \in \mathfrak{S}_n$ such that for any $1\leqslant i\leqslant n$, $\gExt^d_A(S_i,A) \cong S'_{\sigma(i)}(l_i)$ for some $l_i \in \mathbb{Z}$. Note that $\{l_1,\cdots,l_n\}$ relies on the choice of the isomorphic classes of graded simple modules $\{S_1,\cdots,S_n\}$. However, if $A/J$ is concentrated in degree $0$, then $\{l_1,\cdots,l_n\}$ is independent on the choices of $\{S_1,\cdots,S_n\}$.

\begin{definition} With the notations as above,
$\{l_1,\cdots,l_n\}$ is called the set of Gorenstein parameters of $A$ with respect to  $\{S_1,\cdots,S_n\}$.
\end{definition}

\begin{lemma}\label{Gorenstein parameters}
Suppose $A$ is a commonly graded AS-Gorenstein algebra of dimension $d$. With the notations as above, for any $1\leqslant i\leqslant n$, $\gExt^d_{A^{o}}(S_i',A)\cong S_{\sigma^{-1}(i)}(l_{\sigma^{-1}(i)})$.
\end{lemma}
\begin{proof}
Assume that $\gExt^d_{A^{o}}(S_i',A)\cong S_{\tilde{\sigma}(i)}(l'_i)$ for some $\tilde{\sigma} \in \mathfrak{S}_n$. Then
\begin{align*}
\gExt_A^d(\gExt_{A^{o}}^d(S'_i,A),A)&\cong\gExt_A^d(S_{\tilde{\sigma}(i)}(l'_i),A)\\
&\cong S'_{\sigma\tilde{\sigma}(i)}(-l'_i+l_{\tilde{\sigma}(i)}).
\end{align*}
Since $\gExt_A^d(\gExt_{A^{o}}^d(S'_i,A),A)\cong S'_i$ by  definition,  $\tilde{\sigma}=\sigma^{-1}$ and $l'_i=l_{\sigma^{-1}(i)}$.
\end{proof}

The following lemma is proved for $\mathbb{N}$-graded algebras $A$ where $A_0$ is a finite product of the base field in \cite{V,V2}. The proof there works in the more general case.
\begin{lemma} \label{ext(S,Q)} 
Suppose $A$ is a commonly graded AS-Gorenstein algebra of dimension $d$. With the notations as above, for any $1\leqslant i,j\leqslant n$,
$$\gExt_A^d(S_i,Q_j)\left\{
\begin{aligned}
&=0,&j\neq \sigma(i)\\
&\neq 0, &j=\sigma(i).
\end{aligned}
\right.
$$
\end{lemma}

\begin{proof}
 Since $Q_j$ is finitely generated graded projective, 
\begin{align*}
 D(\gExt^d_A(S_i,Q_j))&\cong D(\gExt_A^d(S_i,A)\otimes_A Q_j)\cong D(S'_{\sigma(i)}(l_i)\otimes_A Q_j)\\
 &\cong \gHom_{A^{o}}(S'_{\sigma(i)}(l_i),D(Q_j))\\
 &=\gHom_{A^{o}}(S'_{\sigma(i)}(l_i),\soc D(Q_j)).
\end{align*} 
By Lemmas \ref{soc} and \ref{D(Si) cong Si'}, $\soc D(Q_j)=D(Q_j/JQ_j)=D(S_j)\cong S'_j$.  
Hence $\gExt^d_A(S_i,Q_j)=0$ for any $j\neq \sigma(i)$,
which forces
$\gExt_A^d(S_i,Q_{\sigma(i)})\neq 0$.
\end{proof}

Next lemma is well-known (see also \cite[Theorem 1]{V2}).
\begin{lemma}\label{dual property}
Let $R$ be a ring such that $R_R$ has finite injective dimension and $M$ be an $R$-module with a finitely generated projective resolution. Assume that there is an integer $n$ such that $\Ext_R^i(M,R)=0$ for $i\neq n$. Then
\begin{itemize}
\item[(1)] $\Ext_{R^o}^i(\Ext_R^n(M,R),R)=0$ for $i\neq n$.
\item[(2)] $\Ext_{R^o}^n(\Ext_R^n(M,R),R)\cong M$.
\end{itemize}
In the case that $M$ is an $(R,S)$-bimodule, the isomorphism in (2) is an $(R,S)$-bimodule morphism.
\end{lemma}
This lemma also holds in graded case.
\begin{proof}
The conclusion follows from the cohomological spectral sequence of the double complex $$I^\bullet \otimes_R P_\bullet \cong \Hom_{R^o}(\Hom_R(P_{\bullet}, R), I^\bullet)$$ where $P_\bullet$ is a finitely generated projective resolution of ${}_RM$ and $I^\bullet$ is an injective resolution of $R_R$ of finite length. 
\end{proof}

Let $A$ be a commonly graded $k$-algebra with finite left and right injective dimension $d$. Let $\mathcal{X}$ be the full subcategory of $\gr A$ consisting of finite-dimensional modules $M$ such that
\begin{itemize}
\item [(1)] $M$ has a finitely generated graded projective resolution;
\item [(2)] $\gExt_A^i(M,A)=0$ for $i\neq d$;
\item [(3)] $\gExt_A^d(M,A)$ is finite-dimensional and has a finitely generated graded projective resolution as a graded $A^o$-module.
\end{itemize}
Dually, let $\mathcal{Y}$ be the full subcategory of $\gr A^o$  consisting of finite-dimensional modules $N$ such that
\begin{itemize}
\item [(1)] $N$ has a finitely generated graded projective resolution;
\item [(2)] $\gExt_{A^o}^i(N,A)=0$ for $i\neq d$;
\item [(3)] $\gExt_{A^o}^d(N,A)$ is finite-dimensional and has a finitely generated graded projective resolution as a graded $A$-module.
\end{itemize}

Then we have the following lemma, compared with \cite[Lemma 5.1]{RR}.
\begin{lemma}\label{preparation for equivalent defintions of GAS} 
Keep the assumptions as above.
Then
\begin{itemize}
\item[(1)] both $\mathcal{X}$ and $\mathcal{Y}$ are closed under taking extensions and direct summands.
\item[(2)] $\gExt_A^d(-,A)$ and $\gExt^d_{A^o}(-,A)$ induce a duality between $\mathcal{X}$ and $\mathcal{Y}$.
\item[(3)] If $M\in \mathcal{X}$ is a graded $A\otimes B^o$-module, where $B$ is a commonly graded algebra, then $M\cong \gExt^d_{A^o}(\gExt^d_A(M,A),A)$ as graded $A\otimes B^o$-modules.
\end{itemize}
\end{lemma}
\begin{proof}
(1) is obvious. 
(2) and (3) follows from Lemma \ref{dual property}.
\end{proof}

Next theorem characterizes Ext-finite commonly graded AS-Gorenstein algebras. For $\mathbb{N}$-graded AS-regular algebras, it is proved in \cite[Theorem 5.2]{RR}.
We add a proof here for the convenience of readers.
\begin{theorem}\label{equvalent definitions of AS-Gorenstein}
Let $A$ be an Ext-finite commonly graded algebra with finite graded left and right injective dimension $d$ and $S=A/J$. Suppose $\gExt^i_A(S,A)=0$ for $i\neq d$. Then the following are equivalent.
\begin{itemize}
\item[(1)] $A$ is commonly graded AS-Gorenstein. 
\item[(2)] $\gExt^d_A(-,A)$ gives a contravariant equivalence between the categories  $\fd A$ and $\fd A^o$.
\item[(3)] $\gExt^d_A(-,A)$ gives a contravariant equivalence between the categories $\gr S$ and $\gr S^o$ where $\gr S$ (resp. $\gr S^o$) is regarded as the full subcategory of $\gr A$ (resp. $\gr A^o$) consisting of finitely generated graded semisimple objects via the algebra morphism $A\to A/J$.
\item[(4)] $\gExt^d_A(S,A)J=0$, and there is an  invertible graded $(S,S)$-bimodule $V$ such that $\gExt^d_A(S,A)$, viewed as a graded $S^o$-module, is isomorphic to $V$.
\item[(5)] $\gExt^d_A(S,A)J=0$, and $\gExt^d_A(S,A)$, viewed as a graded $(S,S)$-bimodule, is invertible.
\end{itemize}
\end{theorem}
\begin{proof}
Since $A$ is Ext-finite, every finite-dimensional graded $A$-module or $A^o$-module has a finitely generated graded projective resolution by Corollary \ref{Ext-finite implies f.d. module has f.g. proj. reso.}. 

Let $\mathcal{X}$ and $\mathcal{Y}$ be the full subcategories as defined right before Lemma \ref{preparation for equivalent defintions of GAS}.

(1) $\Rightarrow$ (2). By 
the definition of commonly graded AS-Gorenstein algebras, every graded simple $A$-module (resp. $A^o$-module) is contained in $\mathcal{X}$ (resp. $\mathcal{Y}$).
Since $\mathcal{X}$ and $\mathcal{Y}$ are closed under extensions by Lemma \ref{preparation for equivalent defintions of GAS},  $\mathcal{X}$ is in fact the full subcategory of  finite-dimensional graded $A$-modules and $\mathcal{Y}$ is the full subcategory of  finite-dimensional  graded $A^o$-modules. Then (2) follows from  Lemma \ref{preparation for equivalent defintions of GAS} (2).

(2) $\Rightarrow$ (3). 
The contravariant equivalent functors preserve finitely generated graded semisimple objects.

(3) $\Rightarrow$ (1). It follows from Lemma \ref{dual property} that $\gExt_A^d(-,A)$ and $\gExt_{A^o}^d(-,A)$ give a duality between $\gr S$ and $\gr S^o$. Thus $\gExt_A^d(-,A)$ gives a bijection between the isomorphic classes of graded simple $A$-modules and the isomorphic classes of graded simple $A^o$-modules, with $\gExt_{A^o}^d(-,A)$ as its inverse.

Now for any graded simple $S^o$-module $M$, there is a graded simple $S$-module $M'$ such that $\gExt_A^d(M',A)\cong M$. Then by Lemma \ref{dual property}, $\gExt_{A^o}^i(M,A)=0$ for $i\neq d$. Therefore, $A$ is commonly graded AS-Gorenstein.

(3) $\Rightarrow$ (4). Since $\gExt_A^d(-,A)$ sends $\gr S$ to $\gr S^o$, $\gExt_A^d(S,A)J=0$. Note that $\gExt_A^d(D(-),A)$ is an auto-equivalent functor of $\gr S^o$. Hence there is an invertible graded $(S,S)$-bimodule $V$ such that $\gExt_A^d(D(-),A)\cong -\otimes_{S} V$. Therefore $\gExt_A^d(-,A)\cong D(-)\otimes_{S}V$ and $\gExt_A^d(S,A)\cong D(S)\otimes_S V\cong V$ as graded $S^o$-modules.

(4) $\Rightarrow$ (2). Clearly, as a left module, $S$ belongs to $\mathcal{X}$. Since $\mathcal{X}$ is closed under taking direct summands and extensions by Lemma \ref{preparation for equivalent defintions of GAS} (1),  $\mathcal{X}$ is the full subcategory of  finite-dimensional graded $A$-modules. By Lemma \ref{dual property}, as a right module, $V \cong  \gExt_A^d(S,A)$ belongs to $\mathcal{Y}$. Since $V$ is a graded projective generator of $\gr S^o$, each indecomposable graded projective  $S^o$-module is a direct summand of some shift $V$, which implies that every simple right $A_0$-module is a direct summand of $V$ up to shift. Thus $\mathcal{Y}$ is the full subcategory  of  finite-dimensional graded $A^o$-modules. Then (2) follows from Lemma \ref{preparation for equivalent defintions of GAS} (2). 

(4) $\Rightarrow$ (5). Since $S\cong \gExt_{A^o}^d(\gExt_A^d(S,A),A)$ as $(S,S)$-bimodules by Lemma \ref{dual property}, the left $S$-module structure of $\gExt_A^d(S,A)$ is faithful. This means the canonical morphism $S\to \gEnd_{S^o}(\gExt_A^d(S,A))$
is injective. 
On the other hand, the isomorphisms $\gEnd_{S^o}(\gExt_A^d(S,A))\cong\gEnd_{S^o}(V)\cong S$ requires the dimension of $S$ is equal to the one of $\gEnd_{S^o}(\gExt_A^d(S,A))$. 
Therefore the structure morphism $S\to \gEnd_{S^o}(\gExt_A^d(S,A))$
is an isomorphism which yields that $\gExt_A^d(S,A)$ is an invertible graded $(S,S)$-bimodule. 

(5) $\Rightarrow$ (4). Obvious.
\end{proof}

\begin{remark} A right-module version of Theorem \ref{equvalent definitions of AS-Gorenstein} holds
by interchanging the left graded module category with the right graded module category and replacing $\gExt_A^i(-,A)$ with $\gExt_{A^o}^i(-,A)$.
\end{remark}

Next we show that every commonly graded AS-Gorenstein algebra is graded Morita equivalent to a basic commonly graded AS-Gorenstein algebra. Similar result is proved in \cite{IKU}. Recall the concepts first.

\begin{definition}
A (graded) $k$-algebra $A$ is called a basic algebra if $A/J$ is a finite direct sum of division algebras where $J$ is the (graded) Jacobson radical of $A$.
\end{definition}

Note that if $A$ is a basic commonly graded algebra, then $S=A/J$ is concentrated in degree $0$, because finite-dimensional graded division algebras are concentrated in degree $0$ by \cite[Exercise 5.6.2]{NO2}.

Two  commonly graded algebras $A$ and $B$ are graded Morita equivalent if there exists an equivalent functor $F:\Gr A\to \Gr B$ commuting with the shifting functor. This condition is equivalent to that there exists a graded $(B,A)$-bimodule $U$ such that $U_A$ is a finitely generated graded projective generator and the structure morphism $B\to \gEnd_{A^o}(U)$ is an isomorphism \cite[Theorem 2.3.8]{H}.

\begin{proposition}\label{Morita to basic algebra}
Suppose $A$ is a commonly graded algebra. Then $A$ is graded Morita equivalent to a basic commonly graded algebra $B$. 
\end{proposition}
\begin{proof}
Suppose $S=A/J=\mathop{\oplus}\limits_{i=1}^n M_{r_i}(D_i)(\bar{\gamma}^i)$ where the $D_i$'s are division algebras and $\bar{\gamma}^i\in \mathbb{Z}^{r_i}$.  Let $e_i$ be a primitive idempotent of $A$ such that $Se_i$ is a graded  simple module of $M_{r_i}(D_i)(\bar{\gamma}^i)$ for $1 \leq i \leq n$.  Let $e=e_1+\cdots+e_n$ and $P=eA$. Then clearly $P_A$ is a projective generator in $\Gr A^o$ and $\gEnd_{A^o}(P)\cong eAe$. Thus $A$ is graded Morita equivalent to $eAe$. Let $B=eAe$. By \cite[Theorem 21.10]{Lam1}, $B/J(B)\cong \bar{e}(A/J)\bar{e}\cong D_1\oplus\cdots\oplus D_n$. Hence $B$ is a  basic commonly graded algebra.
\end{proof}

The next proposition shows that the commonly graded AS-Gorenstein property is a graded Morita invariant.

\begin{proposition}\label{GAS-Gorenstein is preserved by Morita}
Let $A$ and $B$ are two commonly graded algebras. If $A$ and $B$ are graded Morita equivalent, then
\begin{itemize}
\item [(1)] $A$ is Ext-finite if and only if so is $B$;
\item [(2)] $A$ is noetherian with a balanced dualizing complex if and only if so is $B$;
\item [(3)] $A$ is commonly graded AS-Gorenstein of dimension $d$ if and only if so is $B$.
\end{itemize}
\end{proposition}
\begin{proof}
Let $U$ be a graded $(B,A)$-bimodule such that $U\otimes_A -:\Gr A\to \Gr B$ is an equivalent functor. Let $U^{-1}=\gHom_B(U,B)$. Then $U^{-1}\otimes_B -$ is a quasi-inverse of $U\otimes_A -$.

(1) Suppose $A$ is Ext-finite, which means that every graded simple $A$-module has a finitely generated graded projective resolution. 
As $U\otimes_A-$ sends graded simple $A$-modules to graded simple $B$-modules and preserves the finiteness of projective resolution, $B$ is Ext-finite.

(2) Suppose $A$ is noetherian with a balanced dualizing complex $R$. Note that the equivalent functor $F=U\otimes_A - \otimes_A U^{-1}:\Gr A^e\to \Gr B^e$ extends to an equivalence $F:\D^b(\Gr A^e)\to\D^b(\Gr B^e)$. Direct verification shows that $F(R)$ is a balanced dualizing complex of $B$.

(3) Suppose $A$ is a commonly graded AS-Gorenstein algebra of dimension $d$. 
Then the injective dimension of ${}_AU^{-1}$ is the same as $A$. 
Thus the injective dimension of ${}_BB\cong U\otimes_A U^{-1}$ is $d$. Similarly, the injective dimension of $B_B$ is $d$.

For any graded simple $B$-module $X$, there is a graded simple $A$-module $X'$ such that $X\cong U\otimes_A X'$. Thus, for any $i\neq d$, 
$$\gExt_B^i(X,B)\cong \gExt_B^i(U\otimes_A X',B)\cong \gExt_A^i(X',U^{-1})\cong \gExt_A^i(X',A)\otimes_A U^{-1}=0.$$

Let $S_A=A/J(A)$ and $S_B=B/J(B)$. Regard $\gr S_A$ and $\gr S_B$ as the full subcategory of $\gr A$ and $\gr B$ consisting of graded semisimple modules respectively. 
Since every equivalent functor preserves semisimple objects, 
$U^{-1}\otimes_B -$ can be restricted to an equivalence $\gr S_B\to \gr S_A$.
Similarly, $-\otimes_A U^{-1}:
\gr S_A^o\to \gr S_B^o$ is also an equivalence. 
Therefore 
$$\gExt_B^d(-,B)\cong \gExt_A^d(U^{-1}\otimes_B -,U^{-1})\cong \gExt_A^d(U^{-1}\otimes_B -,A)\otimes_A U^{-1}$$ 
gives a contravariant equivalence from $\gr S_B$ to $\gr S_B^o$. It follows from Theorem \ref{equvalent definitions of AS-Gorenstein} that $B$ is a commonly graded AS-Gorenstein algebra of dimension $d$.
\end{proof}

\section{Existence of Balanced Dualizing Complexes}\label{Existence of Balanced Dualizing Complex}
This section is first dedicated to examining the structure of the minimal graded injective resolution of a commonly graded AS-Gorenstein algebra. Then we prove the existence of balanced dualizing complex for a noetherian commonly graded AS-Gorenstein algebra.

\begin{lemma}\label{local cohomology and injective}
Suppose $A$ is a commonly graded algebra and $I$ is a graded injective $A$-module. If $\soc I$ is finite-dimensional, then $\Gamma_{A}(I)$ is the injective hull of $\soc I$ and $I=\Gamma_A(I)\oplus I'$ where $I'$ is a torsion free graded injective $A$-module.

Moreover, if $A$ is left noetherian, then for any graded injective module $I$, $\Gamma_A(I)$ is the injective hull of $\soc I$ and $I=\Gamma_A(I)\oplus I'$ where $I'$ is a torsion free graded injective $A$-module.
\end{lemma}
\begin{proof}
Assume $I=E\oplus I'$ where $E$ is the graded injective hull of $\soc I$.
If $\soc I$ is finite-dimensional,  then 
$E$ is a finite direct sum of the graded injective hulls of graded simple modules. 
By Proposition \ref{injective hull of simple module},
$E$ is a torsion module. We claim that $I'$ is torsion-free. Otherwise, there is  $0 \neq x\in I'$ such that $Jx=0$. 
Then $Ax$ is a finite-dimensional submodule of $I'$ and $\soc I'\neq 0$. 
Hence $\soc I = \soc E \oplus \soc I' \neq \soc E$, which is a contradiction. Thus $I'$ is torsion-free and $\Gamma_A(I)=E$.

If $A$ is left noetherian, then any direct sum of injective modules is still injective.  
Therefore the injective hull of an arbitrary direct sum of graded simple modules is the direct sum of the injective hulls of those graded simple modules. The assertion follows from a similar argument as above.
\end{proof}

In the rest of this section, we assume that $A$ is an Ext-finite  commonly graded AS-Gorenstein algebra of dimension $d$ and $S=A/J$.

Let $I^\bullet$ be a minimal graded injective resolution of ${}_AA$. Then $\gExt_A^i(S,A)=\gHom_A(S,I^i)\cong \soc I^i$ as graded left $A$-modules. Since $\gExt_A^i(S,A)=0$ for all $i\neq d$, $I^i$ is torsion-free by Lemma \ref{torsion elements and soc M}. Since $\gExt_A^d(S,A) \cong \soc I^d$ is finite-dimensional, by Lemma \ref{local cohomology and injective}, $I^d=I'^d\oplus I''^d$ where $I'^d$ is the injective hull of $\gExt_A^d(S,A)$ 
and $I''^d$ is torsion-free. 

\begin{proposition}\label{D(I'd) is right progenerator}
$D(I'^d)$ is a finitely generated graded projective generator of $\Gr A^o$.
\end{proposition}
\begin{proof}
Since $\soc I'^d\cong\gExt_A^d(S,A)\cong V$ as left $A$-modules for some invertible graded $(S,S)$-bimodule $V$ by Theorem \ref{equvalent definitions of AS-Gorenstein}, $\soc I'^d$ is a graded projective generator in $\gr S$.
It follows from Proposition \ref{injective hull of simple module} that
$D(I'^d)$ is a finitely generated graded projective generator in $\Gr A^o$.
\end{proof}

Propositions \ref{gExt(M,A) cong gHom(M,I'd) as right module} is proved in \cite{V2} in the case where $A_0$ is a direct product of the base field $k$. The proof still works when $A$ is commonly graded.

\begin{proposition}\label{gExt(M,A) cong gHom(M,I'd) as right module} The left $A$-module
$I'^d$ can be endowed with a graded $(A,A)$-bimodule structure, so that
for any $M\in \fd A$, as graded $A^o$-modules
$$\gExt_A^i(M,A)\cong \left\{
\begin{aligned}
&0,&i\neq d,\\
&\gHom_A(M,I'^d) &i=d.
\end{aligned}
\right.$$
If moreover $M$ is a graded bimodule, then the isomorphism is a bimodule isomorphism.
\end{proposition}
\begin{proof}
Since for all $p\leqslant d-1$, $I^p$ is torsion-free,
$\gHom_A(M,I^p)=0$ 
and thus $\gExt_A^d(M,A)=\gHom_A(M,I'^d\oplus I''^d)=\gHom_A(M,I'^d)$ as vector spaces.

Let $\phi(a): {}_AA\to {}_AA$ be the morphism induced by the right multiplication of $a\in A$. By lifting $\phi(a)$ to the minimal graded injective resolution of ${}_AA$, there is a     commutative diagram:
\begin{center}
\begin{tikzcd}
0 \arrow[r] & A \arrow[r, "\partial^{-1}"] \arrow[d, "\phi(a)"] & I^0 \arrow[r, "\partial^0"] \arrow[d, "\phi^0(a)"] & \cdots \arrow[r, "\partial^{d-1}"] & I^d \arrow[r] \arrow[d, "\phi^d(a)"] & 0 \\
0 \arrow[r] & A \arrow[r, "\partial^{-1}"]                      & I^0 \arrow[r, "\partial^0"]                        & \cdots \arrow[r, "\partial^{d-1}"] & I^d \arrow[r]                        & 0.
\end{tikzcd}
\end{center}
Since $I^d=I'^d\oplus I''^d$ and $\gHom(I'^d,I''^d)=0$, the morphism $\phi^d(a)$ can be as expressed  
\begin{center}
$\phi^d(a)=
\begin{pmatrix}
\lambda_1(a) & \lambda_2(a)\\
0 & \lambda_3(a)
\end{pmatrix}$
\end{center}
where $\lambda_1(a):I'^d\to I'^d,\lambda_2(a):I''^d\to I'^d$ and $\lambda_3(a):I''^d\to I''^d$ are graded $A$-module morphisms.

We claim that $\lambda_1(a)$ is independent of the choices of the lifting of $\phi(a)$.

Suppose $\{\phi'^i(a)\}$ is another lifting of $\phi(a)$. Then there is a sequence of morphisms $s_k:I^k\to I^{k-1}$ such that $\phi^k(a)-\phi'^k(a)=\partial^{k-1}s_k+s_{k+1}\partial^k$. 
In particular, $\phi^d(a)-\phi'^d(a)=\partial^{d-1}s_d$. We may assume 
\begin{center}
$\partial^{d-1}=
\begin{pmatrix}
\partial'^{d-1}\\
\partial''^{d-1}
\end{pmatrix}$
and
$s_d=
\begin{pmatrix}
s'_d &s''_d
\end{pmatrix}$.
\end{center}
Note that $s'_d:I'^d\to I^{d-1}$ is a morphism from torsion module to torsion-free module. Hence $s'_d=0$
and
\begin{center}
$\phi^d(a)-\phi'^d(a)=
\begin{pmatrix}
\lambda_1(a)-\lambda'_1(a) & \lambda_2(a)-\lambda'_2(a)\\
0 & \lambda_3(a)-\lambda'_3(a)
\end{pmatrix}
=
\begin{pmatrix}
0 & \partial'^{d-1}s''_d\\
0 & \partial''^{d-1}s''_d
\end{pmatrix}$,
\end{center}
It follows that $\lambda_1(a)=\lambda'_1(a)$.

Clearly the map $\lambda_1:A\to \gEnd_A(I'^d)^{o}, \, a\mapsto \lambda_1(a)$ is an algebra morphism. Thus the left $A$-module $I'^d$ is an $(A,A)$-bimodule with the right action given by the algebra morphism $\lambda_1$.

Meanwhile, the canonical right $A$-module structure of $\gExt^d_A(M,A)$
can be interpreted as in the following way. Any $x\in\gExt_A^d(M,A)$ can be seen as $f+\im(\partial^{d-1})_*$, where $f$ is a morphism from $M$ to $I^d$. Since $M$ is finite-dimensional and $I^{d-1}$ is torsion-free, $(\partial^{d-1})_*=0$. Hence $x$ can be identified with $f$. Then for any $a\in A$, $x\cdot a$ is the composition of $f$ and the $d$-th lifting of $\phi(a)$, i.e. $x\cdot a=\phi^d(a)\circ f$. Now for $f\in \gHom_A(M,I^d)=\gHom_A(M,I'^d)$, which is a right $A$-module induced by the $(A,A)$-bimodule structure of $I'^d$ just given,  $f\cdot a=\lambda_1(a)\circ f$. By abuse of notation, $\phi^d(a)\circ f=\lambda_1(a)\circ f\in \gHom_A(M,I'^d)$. Thus $\gExt^d_A(M,A)\cong \gHom_A(M,I'^d)$ as right $A$-modules.
\end{proof}

Let $0\to A_A\to J^\bullet$ be a minimal graded injective resolution of $A_A$. Similarly, $J^i$ is torsion-free for all $i\neq d$ and $J^d=J'^d\oplus J''^d$ where
$J'^d$ is the injective hull of $\gExt_{A^o}^d(S,A)$ and
$J''^d$ is torsion-free.

\begin{proposition}\label{D(J') is left progenerator} 
$D(J'^d)$ is a finitely generated graded projective generator  of $\Gr A$.
\end{proposition} 

\begin{proposition}
\label{gExt(M,A) cong gHom(M,J'd) as left module} The right $A$-module
$J'^d$ can be endowed with a graded $(A,A)$-bimodule structure, so that for any $N\in \fd A^o$, as graded $A$-modules,
$$\gExt_{A^{o}}^i(N,A)\cong \left\{
\begin{aligned}
&0,&i\neq d\\
&\gHom_{A^{o}}(N,J'^d),&i=d.
\end{aligned}
\right.$$
 If moreover $N$ is a graded bimodule, then the isomorphism is a bimodule isomorphism.
\end{proposition}

Next we study the relation between $I'^d$ and $J'^d$. Similar result is proved in \cite{V2} in $\mathbb{N}$-graded case where $A_0$ is a direct product of the base field $k$.  Recall that $\Gamma_A$ and $\Gamma_{A^o}$ are the local cohomology functors.

\begin{proposition}\label{I and J iso as bimodule}
With the above notations, if $A$ is an Ext-finite commonly graded AS-Gorenstein algebra of dimension $d$, then $I'^d\cong J'^d$ as graded $(A,A)$-bimodules.
\end{proposition}
\begin{proof}
Note that in $\D(\Gr A^e)$, $R\Gamma_A(A)\cong I'^d[-d]$ and $R\Gamma_{A^o}(A)\cong J'^d[-d]$, and  both of them have torsion cohomologies. Then
$$R\Gamma_{A^o}\circ R\Gamma_A(A)\cong R\Gamma_{A^o}(I'^d[-d])\cong I'^d[-d]$$ and 
$$R\Gamma_A\circ R\Gamma_{A^o}(A)\cong
R\Gamma_A(J'^d[-d])\cong J'^d[-d].$$

By \cite[Lemma 4.5]{VdB}, $R\Gamma_{A^e}(A)=R\Gamma_A\circ R\Gamma_{A^o}(A)=R\Gamma_{A^o}\circ R\Gamma_A(A)$ in $\D(\Gr A^e)$.
It follows that $I'^d\cong H^d(R\Gamma_{A^e}(A))\cong J'^d$ as graded $(A,A)$-bimodules.
\end{proof}

Next lemma is a well-known fact.
\begin{lemma}\label{Hom iso}
Let $A$ be a $k$-algebra and $M$ a finitely presented left $A$-module. Then for any left $A$-module $X$, there is a natural vector space isomorphism
$$\Hom_k(\Hom_A(M,X),k))\cong \Hom_k(X,k)\otimes_A M.$$
The result holds also in the graded case.
\end{lemma}

Based on Proposition \ref{I and J iso as bimodule}, we prove an important fact that $D(I'^d)$ and $D(J'^d)$ are invertible graded $(A,A)$-bimodules. When $A$ is basic noetherian $\mathbb{N}$-graded, this is also proved in \cite[Proposition 2.17]{IKU}.

\begin{proposition}\label{D(I'd) is invertible}
Let $A$ be an Ext-finite commonly graded AS-Gorenstein algebra of dimension $d$. Then $D(R^d\Gamma_A(A))$ and $D(R^d\Gamma_{A^o}(A))$ are invertible graded $(A,A)$-bimodules.
\end{proposition}
\begin{proof}
With the above notations, $R^d\Gamma_A(A)\cong I'^d$ and $R^d\Gamma_{A^o}(A)\cong J'^d$. By Propositions \ref{D(I'd) is right progenerator}, \ref{D(J') is left progenerator} and \ref{I and J iso as bimodule},
$D(I'^d)$ is a finitely generated graded projective generator as a left or a right graded $A$-module.

Proposition \ref{gExt(M,A) cong gHom(M,I'd) as right module} yields that 
$\gHom_A(S,I'^d)=\gExt_A^d(S,A)\cong V$ as $(S,S)$-bimodules for some invertible graded $(S,S)$-bimodule $V$.
By Lemma \ref{Hom iso},
$$D(I'^d)\otimes_A S\cong D(\gHom_A(S,I'^d))\cong D(V),$$
where $D(V)=V^{-1}$ is the inverse of $V$ as an invertible $(S,S)$-bimodule. Thus we have 
$$J(D(I'^d)\otimes_A S)=0.$$

Similarly $(S\otimes_A D(J'^d))J=0$. Since $I'^d\cong J'^d$,  
$(S\otimes_A D(I'^d))J=0.$
By \cite[Lemma 5.13]{RR}, $D(I'^d)$ is an invertible graded $(A,A)$-bimodule, and so is $D(J'^d)$.
\end{proof}

\begin{corollary}\label{gExt_A(S,A) cong gExt_Ao(S,A)}
Let $A$ be an Ext-finite commonly graded AS-Gorenstein algebra of dimension $d$ and $\omega=\gHom_{A^o}(D(R^d\Gamma_A(A)),A)$. Then $\gExt_A^d(S,A)\cong \gExt_{A^o}^d(S,A)\cong \omega/J\omega$ as graded $(A,A)$-bimodules.
\end{corollary}
\begin{proof}
Let $V_1=\gExt_A^d(S,A)$ and $V_2=\gExt_{A^o}^d(S,A)$. 
The proof of Proposition \ref{D(I'd) is invertible} shows that
$$D(I'^d)\otimes_A S\cong V_1^{-1} \, \textrm{ and } \, S\otimes_A D(J'^d)\cong V_2^{-1}.$$
Since $D(I'^d)\cong D(J'^d)$, by \cite[Lemma 2.11]{RR},
$$D(I'^d)\otimes_A S\cong S\otimes_A D(J'^d).$$
Therefore $\gExt_A^d(S,A)\cong \gExt_{A^o}^d(S,A)$ as graded $(A,A)$-bimodules.

Note that the inverse $(S,S)$-bimodule of $D(I'^d)\otimes_A S$ is $\gHom_{A^o}(D(I'^d),A)\otimes_A S$. Thus 
$$\gExt_A^d(S,A)\cong \gHom_{A^o}(D(I'^d),A)\otimes_A S=\omega \otimes_A S\cong \omega/J\omega$$
as graded $(S,S)$-modules, and as graded $(A,A)$-modules. 
\end{proof}

Besides Theorem \ref{equvalent definitions of AS-Gorenstein}, we give more characterizations of Ext-finite commonly graded AS-Gorenstein algebras.

\begin{theorem}\label{another equivalent definition of GAS}
Let $A$ be an Ext-finite commonly graded algebra. Then the following are equivalent:
\begin{itemize}
\item[(1)] $A$ is a commonly graded AS-Gorenstein algebra of dimension $d$.
\item[(2)] $A$ has finite injective dimension $d$ as a left and right $A$-module and there is an invertible graded $(A,A)$-bimodule $U$ such that
$R\Gamma_A(A)\cong D(U)[-d]$
in $\D(\Gr A^e)$.
\end{itemize}
\end{theorem}
\begin{proof}
(1) $\Rightarrow$ (2). It follows from Proposition \ref{D(I'd) is invertible} that $R^d\Gamma_A(A) \cong I'^d\cong D(U)$ for some invertible $(A,A)$-bimodule $U$.

(2) $\Rightarrow$ (1). 
Let $I^\bullet$ be a minimal graded injective resolution of $A$. By Lemma \ref{depth and local dimension},  $\gExt_A^i(S,A)=0$ for all $i\neq d$. 
Since $\gExt_A^i(S,A)\cong\gHom_A(S,I^i)\cong \soc I^i$ as graded $A$-modules, $I^i$ is torsion-free for all $i\neq d$ and $R^d\Gamma_A(A)=\Gamma_A(I^d)=D(U)$ as graded $A$-modules. 
Then $I^d$ can be written as $I^d=I'^d\oplus I''^d$ where $I'^d\cong D(U)$ is torsion and $I''^d$ is torsion-free.

By a similar argument as in the proof of Proposition \ref{gExt(M,A) cong gHom(M,I'd) as right module},  $I'^d$ can be endowed with a graded $(A,A)$-bimodule structure  such that $I'^d\cong R^d\Gamma_A(A)\cong D(U)$ as graded $(A,A)$-bimodules. Then
\begin{align*}
\gExt_A^d(S,A)&\cong \gHom_A(S,I'^d)\\
&\cong \gHom_A(S,D(U))\\
&\cong \gHom_k(U\otimes_A S,k)\\
&\cong \gHom_S(U\otimes_A S,S)
\end{align*}
as graded $(A,A)$-bimodules.

Let $V=U\otimes_A S$. 
Then $V$ is an invertible graded $(S,S)$-bimodule by \cite[Lemma 2.11]{RR}
and thus so is $\gExt_A^d(S,A)$.
By Theorem \ref{equvalent definitions of AS-Gorenstein}, $A$ is a commonly graded AS-Gorenstein algebra.
\end{proof}

When $A$ is $\mathbb{N}$-graded, we can say more than Theorem \ref{equvalent definitions of AS-Gorenstein}.

\begin{theorem}\label{supplement of equivalent definition of GAS}
Let $A$ be an Ext-finite $\mathbb{N}$-graded algebra with finite left and right injective dimension $d$. Suppose $\gExt^i_A(S,A)=0$ for $i\neq d$. Consider the following statements:
\begin{itemize}
\item [(1)] $\gExt^d_A(-,A)$ gives a contravariant equivalence from $\gr A_0$ to $\gr A_0^o$ where $\gr A_0$ (resp. $\gr A_0^o$) is regarded as the full subcategory of $\gr A$ (resp. $\gr A^o$) via the algebra map $A\to A/A_{\geqslant 1}$.
\item [(2)] $\gExt_A^d(A_0,A)\,A_{\geqslant 1}=0$, and when $\gExt^d_A(A_0,A)$ is regarded as a graded $A_0^o$-module, it is isomorphic to $D(A_0)\otimes_{A_0} W$ for some  invertible graded $(A_0,A_0)$-bimodule $W$.
\item [(3)] $\gExt_A^d(A_0,A)\,A_{\geqslant 1}=0$, and when $\gExt^d_A(A_0,A)$ is regarded as a graded $(A_0,A_0)$-bimodule, it is isomorphic to $D(A_0)\otimes_{A_0} W$ for some  invertible graded $(A_0,A_0)$-bimodule $W$.
\item [(4)] $A$ is $\mathbb{N}$-graded AS-Gorenstein.
\end{itemize}
Then (1) $\Rightarrow$ (2) $\Rightarrow$ (3) $\Rightarrow$ (4).

If $A$ is an $\mathbb{N}$-graded AS-Gorenstein algebra such that the Gorenstein parameters are all equal, then (1), (2) and (3) hold.
\end{theorem}
\begin{proof}
(1) $\Rightarrow$ (2). Since $\gExt_A^d(-,A)$ sends $\gr A_0$ to $\gr A_0^o$ and $A_0\in \gr A_0$, $\gExt_A^d(A_0,A)\,A_{\geqslant 1}=0$. By hypothesis $\gExt_A^d(D(-),A)$ is an auto-equivalent functor of $\gr A_0$, and thus there is an  invertible  graded $(A_0,A_0)$-bimodule $W$ such that $\gExt_A^d(D(-),A)\cong -\otimes_{A_0} W$. Hence $\gExt_A^d(A_0,A)\cong D(A_0)\otimes_{A_0} W$ as graded $A_0^o$-modules.

(2) $\Rightarrow$ (3). Let $U=\gExt_A^d(A_0,A)$. By Lemma \ref{dual property}, 
$A_0\cong \gExt_{A^o}^d(\gExt_A^d(A_0,A),A)$ 
as graded $(A_0,A_0)$-bimodules. Thus $U$ is a faithful $A_0$-module. Let $\varphi_1: A_0 \to \gEnd_{A^o}(U)$ and $\varphi_2: A_0\to \gEnd_{A^o}(D(A_0)\otimes_{A_0} W)$ be the morphisms giving by the left $A_0$-module structures. Then $\varphi_1$ is injective as $U$ is a faithful $A_0$-module. 
In fact, $\varphi_2$ is the composition $A_0 \cong (\gEnd_A(A_0))^o \cong \gEnd_{A^o}(D(A_0))\cong \gEnd_{A^o}(D(A_0)\otimes_{A_0} W)$ which shows that
$\varphi_2$ is an isomorphism. 
Since $\gEnd_{A_0}(U)\cong \gEnd_{A_0}(D(A_0)\otimes_{A_0} W)\cong A_0$ is finite-dimensional, $\varphi_1$ is an isomorphism. Then there is a commutative diagram:
\[\xymatrix{
A_0\ar[r]^{\varphi_1}\ar@{.>}[d]^{\psi} & \gEnd_{A^o}(U)\ar[d]^{\cong} \\
A_0\ar[r]^(0.25){\varphi_2}& \gEnd_{A^o}(D(A_0)\otimes_{A_0} W)
}\] 
where $\psi$ is the composition of 
$$A_0\xrightarrow{\varphi_1}\gEnd_{A_0}(U)\cong \gEnd_{A_0}(D(A_0)\otimes_{A_0} W)\xrightarrow{\varphi_2^{-1}} A_0.$$
Thus $U\cong {}^\psi (D(A_0)\otimes_{A_0} W)^1$ as graded $(A_0,A_0)$-bimodules. Since $${}^\psi (D(A_0)\otimes_{A_0} W)^1
\cong D(A_0^\psi) \otimes_{A_0} W  \cong D({}^{\psi^{-1}}\!A_0) \otimes_{A_0} W) \cong D(A_0)^{\psi^{-1}}\otimes_{A_0}  W \cong D(A_0)\otimes_{A_0} {}^\psi W^1,$$ 
we have $U\cong D(A_0)\otimes_{A_0} {}^\psi W^1$, where ${}^\psi W^1$ is an invertible graded $(A_0,A_0)$-bimodule. 

(3) $\Rightarrow$ (4). 
Since $\gExt_A^i(S,A)=0$ for all $i\neq d$, 
by an induction on the length of finite-dimensional graded $A$-module $M$,
$\gExt_A^i(M,A)=0$ for all $i\neq d$. 
In particular,  $$R\gHom_A(A_0,A)\cong D(A_0)\otimes_{A_0} W[-d].$$ Then
\begin{align*}
R\gHom_A(S,A)&\cong R\gHom_A(A_0\otimes_{A_0} S,A)\\
&\cong R\gHom_{A_0}(S,R\gHom_A(A_0,A))\\
&\cong R\gHom_{A_0}(S,D(A_0)\otimes_{A_0} W[-d])\\
&\cong R\gHom_{A_0}(S,D(A_0))\otimes_{A_0} W[-d]\\
&\cong D(S)\otimes_{A_0} W[-d]\\
&\cong S\otimes_{A_0} W[-d].
\end{align*}
It follows from \cite[Lemma 2.11]{RR} that $S\otimes_{A_0} W$ is an invertible  graded $(S,S)$-bimodule which implies that
(4) is true
by Theorem \ref{equvalent definitions of AS-Gorenstein}.

Next we prove (4) $\Rightarrow$ (1) under the condition that all the Gorenstein parameters of $A$ are equal, that is, $l_1=\cdots =l_n$ (for the definition and notations  see Section \ref{Commonly graded Artin-Schelter Gorenstein Algebras}). 
Let $l=l_1=\cdots=l_n$.

We claim that $\gExt_A^d(-,A)$ sends $\gr A_0$ to $\gr A_0^o$.
Obviously, $\gExt_A^d(S,A)$ is concentrated on degree $-l$ and so is $\gExt_{A^o}^d(S,A)$.
With the notations above, 
$D(I'^d)$ is generated in degree $l$ as a graded $A^o$-module because $I'^d$ is the injective hull of $\gExt_A^d(S,A)$. Similarly $D(J'^d)$ is also generated in degree $l$ as a graded $A$-module.
Let $P=D(I'^d)$. By Proposition \ref{I and J iso as bimodule}, $D(I'^d)\cong D(J'^d)$ as graded $(A,A)$-bimodules. So $P$ is generated in degree $l$ on both sides. Therefore $A_{\geqslant 1}P=P_{\geqslant l+1}=PA_{\geqslant 1}$.

Suppose $M\in \gr A_0$. Then by Proposition \ref{gExt(M,A) cong gHom(M,I'd) as right module}, 
$$\gExt_A^d(M,A)\cong \gHom_A(M,I'^d)\cong \gHom_A(M,D(P))\cong D(P\otimes_A M)$$
as graded $A^o$-modules.
Since $A_{\geqslant 1}(P\otimes_A M)=A_{\geqslant 1}P\otimes_A M=PA_{\geqslant 1}\otimes_A M=P\otimes_A A_{\geqslant 1}M=0$, 
$D(P\otimes_A M)A_{\geqslant 1}=0$ which implies that $\gExt_A^d(M,A)$ is in $\gr A_0^o$. Similarly $\gExt_{A^o}^d(-,A)$ sends $\gr A_0^o$ to $\gr A_0$.

It follows from Lemma \ref{dual property} that
$$\gExt_{A^o}^d(\gExt_A^d(M,A),A)\cong M \textrm{ and }  \gExt_A^d(\gExt_{A^o}^d(N,A),A)\cong N$$ 
for any $M\in \gr A_0$ and  any $N\in \gr A_0^o$.
Hence $\gExt^d_A(-,A)$ gives a contravariant equivalence from $\gr A_0$ to $\gr A_0^o$, that is, (1) holds.
\end{proof}

We now prove the main result of this section.
\begin{theorem}\label{existence of baalanced dualizing complex of GAS}
Let $A$ be a noetherian commonly graded AS-Gorenstein algebra of dimension $d$. Then $A$ admits a balanced dualizing complex given by $D(R^d\Gamma_A(A))[d]$.
\end{theorem}
\begin{proof}
Recall that if $I^\bullet$ is a minimal graded injective resolution of ${}_AA$, then $I^i$ is torsion-free for all $i\neq d$ and $R^d\Gamma_A(A) \cong I'^d$. Thus $R\Gamma_A(A)\cong R^d\Gamma_A(A)[-d] \cong I'^d[-d]$ in $\D(\Gr A^e)$.

It follows from Proposition \ref{D(I'd) is invertible} that $D(R^d\Gamma_{A}(A))\cong D(I'^d)$ is an invertible graded $(A,A)$-bimodule. 
Hence $D(I'^d)$ has finite injective dimension on both sides.
These imply that $ D(R^d\Gamma_A(A))[d]$ is a dualizing complex of $A$ according to Definition \ref{Dualizing complex}. 

It remains to show that the dualizing complex $D(R^d\Gamma_A(A))[d]$ is balanced.
For convenience, let $U=D(I'^d)$.
Then $R^d\Gamma_{A}(A)\cong D(U)$ and thus
 $$R\Gamma_{A}(U)\cong R\Gamma_{A}(A)\otimes_A U\cong D(U)\otimes_A U[-d]\cong D(A)[-d]$$
in $\D(\Gr A^e)$. 
So,
$R\Gamma_{A}(D(I'^d)[d])\cong D(A)$
in $\D(\Gr A^e)$. 

Similarly, $R\Gamma_{A^o}(D(J'^d)[d])\cong D(A)$
in $\D(\Gr A^e)$.
Since $I'^d\cong J'^d$ as graded $(A,A)$-bimodules, $D(R^d\Gamma_A(A))[d]$ is a balanced dualizing complex of $A$.
\end{proof}

The following corollary is direct from Theorems \ref{existence of baalanced dualizing complex of GAS} and \ref{existence of balanced dualizing complex}.
 
\begin{corollary}\label{corollary of existence of balanced dualizing complex}
If $A$ is a noetherian commonly graded AS-Gorenstein algebra, then
\begin{itemize}
\item[(1)] $A$ satisfies the $\chi$-condition.
\item[(2)] both $\Gamma_A$ and $\Gamma_{A^o}$ have finite cohomological dimension.
\end{itemize}
\end{corollary}

Next proposition is originally proved in \cite[Proposition 8.4]{VdB} and \cite[Lemma 3.5]{RRZ}.

\begin{proposition}\label{Ext(A,Ae) and balanced dualizing complex}
If $A$ is a noetherian commonly graded AS-Gorenstein algebra of dimension $d$, then for any $i\neq d$, $\gExt^i_{A^e}(A,A^e)=0$, and as graded $(A,A)$-bimodules
$$\gExt^d_{A^e}(A,A^e)\cong \gHom_{A^o}(D(R^d\Gamma_A(A)),A).$$

\end{proposition}
\begin{proof}
Let $R=D(R\Gamma_A(A))$ and $R'=R\gHom_{A^o}(R,A)$.
Then $R$ is a rigid dualizing complex (see \cite[Section 8]{VdB}). 
By Theorem \ref{existence of baalanced dualizing complex of GAS}, $R\cong D(I'^d)[d]$. 
Since $D(I'^d)$ is an invertible graded $(A,A)$-bimodule,
$R\otimes_A^L R'\cong A$ and $R'\otimes_A^L R\cong A$
in $\D(\Gr A^e)$. Therefore
\begin{align*}
R\gHom_{A^e}(A,A^e)&\cong R\gHom_{A^e}(A,(R\otimes_A^L R')\otimes (R'\otimes_A^L R))\\
&\cong R\gHom_{A^e}(A,R\otimes R)\otimes_{A^e} (R'\otimes R')\\
&\cong R\otimes_{A^e} (R'\otimes R')\\
&\cong R'.
\end{align*}
So, for any $i\neq d$, $\gExt^i_{A^e}(A,A^e)=0$, and
$\gExt^d_{A^e}(A,A^e)\cong \gHom_{A^o}(D(R^d\Gamma_A(A)),A).$
\end{proof}

\section{Canonical commonly graded AS-Gorenstein Algebras}\label{Canonical commonly graded AS-Gorenstein Algebras}
In this section we introduces a notion of canonical commonly graded AS-Gorenstein algebras to explore the relationships between different definitions of generalized AS-Gorenstein algebras in literature, in particular, the one named by ``AS-Gorenstein over $A_0$" in \cite{MM}.

To simplify the notations, in this section \textbf{we always assume that $A$ is a commonly graded AS-Gorenstein algebra of dimension $d$ such that $S=A/J$ is concentrated on degree 0} unless otherwise stated. By Propositions \ref{Morita to basic algebra} and \ref{GAS-Gorenstein is preserved by Morita}, every commonly graded AS-Gorenstein algebra is graded Morita equivalent to such an algebra.
If $A$ is $\mathbb{N}$-graded, then $S$ is concentrated in degree $0$ automatically.

\begin{definition}\label{definition of canonical}
Let $A$ be an Ext-finite commonly graded AS-Gorenstein algebra of dimension $d$ with $S=A/J$ concentrated in degree $0$. If $\gExt_A^d(S,A)\cong {}^1S^\delta$ as ungraded $(S,S)$-bimodules where $\delta$ is an automorphism of $S$, then $A$ is called canonical.
\end{definition}

Keep the notations as in Section \ref{Commonly graded Artin-Schelter Gorenstein Algebras}.
Since $S$ is assumed to be concentrated in degree $0$, $S=A/J\cong \mathop{\oplus}\limits_{i=1}^n M_{r_i}(D_i)$ for some division algebras $D_i$, and the Gorenstein parameters are independent on the choice of $\{S_1,\cdots,S_n\}$. Let $\{l_1,\cdots, l_n\}$ be the Gorenstein parameters of $A$ and $\sigma:\{1,\cdots,n\}\to \{1,\cdots,n\}$ the permutation such that $\gExt^d_A(S_i,A)=S'_{\sigma(i)}(l_i)$.
For convenience, let $\tau=\sigma^{-1}$.

Recall that $Q_i$ is a graded projective cover of $S_i$, $Q'_i$ is a graded projective cover of $S'_i$, and $I^\bullet$ is a minimal graded injective resolution of ${}_AA$. For any $1\leqslant i \leqslant n$, let $I^\bullet_i$ be a minimal graded injective resolution of $Q_i$. Then $I^p=\mathop{\oplus}\limits_i (I^p_i)^{(r_i)}$ for any $p\geqslant 0$, 
as $S\cong \mathop{\oplus}\limits_i S_i^{(r_i)}$,  where $(-)^{(r_i)}$ means the direct sum of $r_i$ copies of $(-)$. For any $p\neq d$, $I^p$ is torsion-free, and so is $I^p_i$.

\begin{lemma}\label{soc I}
Keep the assumptions as above. 
Then $\soc I^d_i=S_{\tau(i)}(l_{\tau(i)})$.
\end{lemma}
\begin{proof} 

Since $\gExt_{A^o}^d(-,A)$ gives a bijection between graded simple modules, the length of $\gExt_{A^o}^d(S,A)$ as a graded $A$-module is equal to the one of $S_A$ which is $r_1+\cdots +r_n$. By Corollary \ref{gExt_A(S,A) cong gExt_Ao(S,A)}, as graded $A$-modules the length of $\gExt_A^d(S,A)$ is also equal to the one of $\gExt_{A^o}^d(S,A)$.
So the length of $\mathop{\oplus}\limits_{i=1}^n (\soc I^d_i)^{(r_i)}$ is $r_1+\cdots+r_n$, 
as $\gExt_A^d(S,A)\cong\gHom_A(S,I^d)\cong\soc I^d=\mathop{\oplus}\limits_{i=1}^n (\soc I^d_i)^{(r_i)}$. 
Then, 
each $\soc I^d_i$ has length one.
Note that
$\gExt_A^d(S_j,Q_i)\cong\gHom_A(S_j,I_i^d)\cong\gHom_S(S_j,\soc I_i^d).$
Hence $\gExt_A^d(S_{\tau(i)},A)=
S'_i(l_{\tau(i)})$.
It follows from
Lemma \ref{ext(S,Q)} that
$\gHom_S(S_{\tau(i)},\soc I^d_i)$ is concentrated on degree $-l_{\tau(i)}$
which forces $\soc I^d_i \cong S_{\tau(i)}(l_{\tau(i)})$.
\end{proof}

By Lemma \ref{local cohomology and injective}, suppose $I^d_i=I'^d_i\oplus I''^d_i$ where $I'^d_i$ is torsion and $I''^d_i$ is torsion-free. Then $I'^d_i$ is the injective hull of $\soc I^d_i$ and $I'^d_i \cong D(Q'_{\tau(i)})(l_{\tau(i)})$
by Lemma \ref{soc I} and Proposition \ref{injective hull of simple module}.

\begin{lemma}\label{ri=r_sigma i}
Keep the assumptions as above. Then there exists an automorphism $\delta$ of $S$ and some integers $m_1,\cdots,m_n$ such that $\gExt_A^d(S,A)\cong  \mathop{\oplus}\limits_{i=1}^n {}^1 M_{r_i}(D_i)^\delta(m_i)$ as graded $(S,S)$-bimodules if and only if $r_i=r_{\sigma(i)}$. In this case, $l_i=m_i$ for all $i$.
\end{lemma}
\begin{proof}
``$\Rightarrow$" 
On one hand, it follows from Lemma \ref{Gorenstein parameters} and
Corollary \ref{gExt_A(S,A) cong gExt_Ao(S,A)} that as graded $S^o$-modules
$$\gExt_A^d(S,A) \cong \gExt_{A^o}^d(S,A) \cong \mathop{\oplus}\limits_{i=1}^n (S'_{\sigma(i)})^{(r_i)}(l_i).$$ 
On the other hand, $\gExt_A^d(S,A) \cong \mathop{\oplus}\limits_{i=1}^n M_{r_i}(D_i)$ as ungraded $S^o$-modules.  
Hence $r_i=r_{\sigma(i)}$ for all $i$.

``$\Leftarrow$" 
Note that $r_i=r_{\sigma(i)}$ is equivalent to say $r_i=r_{\tau(i)}$.
Hence, by Corollary \ref{gExt_A(S,A) cong gExt_Ao(S,A)}, as graded $S$-modules,
$$\gExt_A^d(S,A)\cong \gExt_{A^o}^d(S,A)\cong \mathop{\oplus}\limits_{i=1}^n (S_{\tau(i)})^{(r_i)}(l_{\tau(i)})\cong \mathop{\oplus}\limits_{i=1}^n (S_i)^{(r_i)}(l_i).$$
So as ungraded $S$-modules, $\gExt_A^d(S,A)\cong S$. Since $\gExt_A^d(S,A)$ is an invertible $(S,S)$-bimodule, there exists an automorphism $\delta$ of $S$ such that $\gExt_A^d(S,A)\cong {}^1S^\delta$. By Lemma \ref{graded invertible S bimodule}, there exist integers $m_1,\cdots,m_n$ such that as graded $(S,S)$-bimodules, $\gExt_A^d(S,A)\cong \mathop{\oplus}\limits_{i=1}^n {}^1 M_{r_i}(D_i)^\delta(m_i)$. The isomorphism $\gExt_A^d(S,A)\cong \mathop{\oplus}\limits_{i=1}^n (S_i)^{(r_i)}(l_i)$ forces $l_i=m_i$.
\end{proof}

Based on Lemmas \ref{graded invertible S bimodule} and \ref{ri=r_sigma i}, the following proposition holds.

\begin{proposition}\label{notations of canonical commonly graded AS}
Keep the notations and assumptions as above. Then the following are equivalent.
\begin{itemize}
\item[(1)] $A$ is canonical.
\item[(2)] $r_i=r_{\sigma(i)}$ for all $i$.
\item[(3)] $\gExt_A^d(S,A)\cong \mathop{\oplus}\limits_{i=1}^n {}^1 M_{r_i}(D_i)^\delta(l_i)$ where $\delta$ is an automorphism of $S$.

\end{itemize}
\end{proposition}

When discussing the canonical commonly graded AS-Gorenstein algebras, we always keep the notations in Proposition \ref{notations of canonical commonly graded AS}.

\begin{remark}
By Lemma \ref{invertible bimodule of basic semisimple algebra}, every Ext-finite basic commonly graded AS-Gorenstein algebra is canonical. It follows from Propositions \ref{Morita to basic algebra} and \ref{GAS-Gorenstein is preserved by Morita} that every Ext-finite commonly graded AS-Gorenstein algebra is graded Morita equivalent to a canonical one.
\end{remark}

Suppose $A$ is a canonical commonly graded AS-Gorenstein algebra. 
With the notations as above, by Lemma \ref{ri=r_sigma i},
$$I'^d=\mathop{\oplus}\limits_{i=1}^n D(Q'_{\tau(i)})^{(r_i)}(l_{\tau(i)})=\mathop{\oplus}\limits_{i=1}^nD(Q'_i)^{(r_i)}(l_i).$$

Let $F$ be the forgetful functor sending graded $A$-modules to ungraded $A$-modules. Recall Proposition \ref{gExt(M,A) cong gHom(M,I'd) as right module} states that $\lambda_1:A\to (\gEnd_A(I'^d))^o$ gives $I'^d$ a graded $A^o$-module structure so that $I'^d$ becomes a graded $(A,A)$-bimodule. Then $F(D(I'^d))$ is an $(A,A)$-bimodule with the right $A$-module structure isomorphic to $A_A$ and the left structure induced by $\lambda_1$. Thus there is an algebra morphism $\mu:A\to A$ such that $F(D(I'^d))\cong {}^\mu A^1$ and $F(I'^d)={}^1F(D(A))^\mu$.

More precisely, 
suppose $\{\mathfrak{e}_1,\cdots,\mathfrak{e}_n\}$ is a set of complete orthogonal idempotents of $A$ such that $\bar{\mathfrak{e}}_i\in S=A/J=\mathop{\oplus}\limits_{i=1}^n M_{r_i}(D_i)$ is the identity of $M_{r_i}(D_i)$. 
Then $\mathfrak{e}_iA\cong (Q'_i)^{(r_i)}$ and
 $I'^d=\mathop{\oplus}\limits_{i=1}^nD(Q'_i)^{(r_i)}(l_i)=\mathop{\oplus}\limits_{i=1}^nD(\mathfrak{e}_iA)(l_i)$. Thus the canonical map $\theta: D(I'^d)\to\mathop{\oplus}\limits_{i=1}^n\mathfrak{e}_iA(-l_i)$ is an isomorphism of graded $A^o$-modules. 
Then $\mu$ is the composition of following maps:
\begin{center}
$
A\xrightarrow{\lambda_1} \gEnd_A(I'^d)^{o} \xrightarrow{D(-)} 
\gEnd_{A^{o}}(D(I'^d))\xrightarrow{f}$\\
$\gEnd_{A^{o}}(\mathop{\oplus}\limits_{i=1}^n\mathfrak{e}_iA(-l_i))
\xrightarrow{g} \mathop{\oplus}\limits_{i=1}^n\mathfrak{e}_iA(-l_i)\xrightarrow{F} F(\mathop{\oplus}\limits_{i=1}^n\mathfrak{e}_iA(-l_i))=A
$
\end{center}
where $f(\varphi)=\theta\circ \varphi \circ \theta^{-1}$ for any $\varphi\in \gEnd_{A^{o}}(D(I'^d))$ and $g(\psi)=\psi(\sum \mathfrak{e}_i(-l_i))$ for any $\psi\in \gEnd_{A^{o}}(\mathop{\oplus}\limits_{i=1}^n\mathfrak{e}_iA(-l_i))$, and $F(\sum \mathfrak{e}_i(-l_i)a_i)=\sum \mathfrak{e}_ia_i$ where $\mathfrak{e}_i(-l_i)$ means $\mathfrak{e}_i$ being shifted $(-l)_i$-th. 

\begin{lemma}\label{mu is automorphism}
If $A$ is a canonical commonly graded AS-Gorenstein algebra, 
with the notations as above, then $\mu$ is an automorphism of $A$. Furthermore, $\mu$ is graded if and only if all the $l_i$'s are equal.
\end{lemma}
\begin{proof}
It follows from Theorem \ref{D(I'd) is invertible} that $\lambda_1$ is an isomorphism. All other maps in the construction of $\mu$ are isomorphisms. Thus $\mu$ is an automorphism. The verification of the second statement is routine.
\end{proof}

When $A$ is $\mathbb{N}$-graded and $A_0$ is a finite product of base field, \cite{V2} also constructed the automorphism $\mu$. We generalize their construction.

\begin{definition}
The automorphism $\mu$ is called the Nakayama automorphism of $A$.
\end{definition}

\begin{proposition}\label{balanced dualizing module of canonical AS-Gorenstein}
Keep the notations and assumptions as above.
Then $A$ is canonical if and only if there is an automorphism $\mu$ of $A$ (not necessarily graded) such that $D(R^d\Gamma_A(A))\cong {}^\mu A^1$ as ungraded $(A,A)$-bimodules.
\end{proposition}
\begin{proof}
It remains to prove sufficiency.
It follows from the proof of Proposition \ref{D(I'd) is invertible} that
$$J(D(R^d\Gamma(A))\otimes_AS)=0=(S\otimes_A D(R^d\Gamma_A(A)))J.$$
So $J\cdot{}^\mu S^1={}^1S^{\mu^{-1}}\cdot J=0$. It follows that $\mu(J)=J$. Therefore $\mu$ induces an automorphism of $S$. 

By the proof of Proposition \ref{D(I'd) is invertible} again,
$D(R^d\Gamma(A))\otimes_AS \cong D(\gExt^d_A(S,A))$. So
$$\gExt_A^d(S,A)\cong {}^1S^{\delta}$$
as ungraded $(S,S)$-bimodules where $\delta$ is the automorphism of $S$ induced by $\mu$.
\end{proof}

\begin{corollary}
Suppose $A$ is a canonical commonly graded AS-Gorenstein algebra of dimension $d$.
\begin{itemize}
    \item [(1)] Let $\mu$ be the Nakayama automorphism of $A$. Then $\mu(J)=J$.
    \item [(2)] Let $\delta$ be the automorphism of $S=A/J$ induced by $\mu$. Then $\gExt_A^d(S,A)\cong {}^1S^\delta$ as ungraded $(S,S)$-bimodules.
\end{itemize}
\end{corollary}

Now Proposition \ref{gExt(M,A) cong gHom(M,I'd) as right module} becomes the following.

\begin{proposition}\label{graded Ext(M,A)=Hom(M,D(A)}
 If $A$ is a canonical commonly graded AS-Gorenstein algebra of dimension $d$, then for any finite-dimensional graded $A$-module (resp. $(A,A)$-bimodule) $M$,
$$F\gExt_A^d(M,A)\cong (FD(M))^\mu$$
as ungraded $A^o$-modules (resp. $(A,A)$-bimodules).

If moreover $l_1=\cdots=l_n=l$, then
$$\gExt_A^d(M,A)\cong \gHom_A(M,{}^1D(A)^\mu(l))\cong D({}^\mu M)(l)$$
as graded $A^o$-modules (resp. $(A,A)$-bimodules).
\end{proposition}
\begin{proof}
By Proposition \ref{gExt(M,A) cong gHom(M,I'd) as right module}, $\gExt_A^d(M,A)\cong\gHom_A(M,I'^d)$.
Since $M$ is finite-dimensional,
\begin{align*}
    \gHom_A(M,I'^d)&=\Hom_A(FM,F(I'^d))\\
    &\cong \Hom_A(FM,{}^1F(D(A))^{\mu})\\
    &\cong \Hom_A(FM,F(D(A)))^{\mu}\\
    &\cong F(\gHom_A(M,D(A)))^{\mu}\\
    &\cong (FD(M))^\mu.
\qedhere
\end{align*}
\end{proof}

Similarly, if we consider $A^o$ and $J'^d$, there is an automorphism $\mu'$ of $A$, such that $F(J'^d)\cong {}^{\mu'}F(D(A))^1$. It follows from Proposition \ref{I and J iso as bimodule} that $\mu'=\mu^{-1}$. 

\begin{proposition}\label{graded Ext(M,A)=Hom(M,D(A) in Aop}
If $A$ is a canonical commonly graded AS-Gorenstein algebra of dimension $d$, then for any finite-dimensional graded $A^o$-module (resp. $(A,A)$-bimodule) $N$,
$$F\gExt_{A^o}^d(N,A)\cong {}^{\mu'}(FD(N))$$
as ungraded $A$-modules (resp. $(A,A)$-bimodules).

If moreover $l_1=\cdots=l_n=l$,then
$$\gExt_{A^o}^d(N,A)\cong \gHom_{A^o}(N,{}^{\mu'}D(A)^1(l))\cong D(N^{\mu'})(l)$$
as graded $A$-modules (resp. $(A,A)$-bimodules).
\end{proposition}

\begin{corollary}\label{delta on e_i}
Suppose $A$ is a canonical commonly graded AS-Gorenstein algebra of dimension $d$ with Nakayama automorphism $\mu$. Let $\delta$ be the automorphism of $S=A/J$ induced by $\mu$ and $e_i$ be the idempotent of $A$ such that $S_i\cong S\bar{e}_i$. Then $\delta(\bar{e}_{\sigma(i)})=\bar{e}_i$.
\end{corollary}
\begin{proof}
With the notations as above, $\gExt^d_A(S_i,A)\cong S'_{\sigma(i)}=\bar{e}_{\sigma(i)}S$ as ungraded modules. By Proposition \ref{graded Ext(M,A)=Hom(M,D(A)}, $\gExt^d_A(S_i,A)\cong \bar{e}_iS^\delta$ as ungraded modules. It follows that $\delta(\bar{e}_{\sigma(i)})=\bar{e}_i$.
\end{proof}
When $A$ is basic, Corollary \ref{delta on e_i} can be  proved also by the structure of $D(I'^d)$ and \cite[Proposition 5.2]{RR2}.

Recall the definition of ``AS-Gorenstein over $A_0$" given in \cite{MM}. 

\begin{definition}\cite[Definition 3.1]{MM}
A locally finite $\mathbb{N}$-graded algebra $A$ is called AS-Gorenstein (resp. regular) over $A_0$ of dimension $d$, if $A$ has finite left injective (resp. global) dimension $d$, and
$$R\gHom_A(A_0,A)\cong {}^1 D(A_0)^\nu(l)[-d]$$
in $\D(\Gr A^e)$ for some automorphism $\nu$ of $A_0$. The integer $l$ is called Gorenstein parameter of $A$.
\end{definition}

\begin{remark}
\begin{itemize}
\item[(1)] The original definition \cite[Definition 3.1]{MM} ``AS-regular over $A_0$" requires further $\gldim A_0<\infty$, but this is a consequence of $\gldim A<\infty$ (see \cite[Lemma 3.23]{RR}).
\item[(2)] If $A$ (resp. $A^o$) is AS-regular over $A_0$, then $A^o$ (resp. $A$) is regular over $A_0$ and $A$ is Ext-finite (see \cite[Proposition 2.3 and Proposition 3.6]{MM}). 
\end{itemize}
\end{remark}

\begin{definition}\cite[Definition 3.9]{MM}
A locally finite $\mathbb{N}$-graded algebra $A$ is called ASF-Gorenstein (resp. ASF-regular) of dimension $d$ with Gorenstein parameter $l$, if $A$ has finite left injective (resp. global) dimension $d$ and 
$$R\Gamma_A(A)\cong {}^1D(A)^{\mu}(l)[-d]$$ 
in $\D(\Gr A^e)$, where $\mu$ is a graded automorphism of $A$.
\end{definition}

``ASF" here stands for Artin-Schelter-Frobenius.
The following theorems show the relations between these definitions. The implication (3) $\Rightarrow$ (1) is proved in \cite[Theorem 2.10]{U} under the assumption that $A$ is noetherian with finite global dimension.
\begin{theorem}\label{GAS-Gorenstein and AS over A0}
Let $A$ be an Ext-finite $\mathbb{N}$-graded algebra. Then the following are equivalent.
\begin{itemize}
\item[(1)] $A$ and $A^o$ are AS-Gorenstein over $A_0$ of dimension $d$ with Gorenstein parameter $l$.
\item[(2)] $A$ is a canonical $\mathbb{N}$-graded AS-Gorenstein algebra of dimension $d$ and $l_1=\cdots=l_n=l$. 
\item[(3)] $A$ and $A^o$ are ASF-Gorenstein algebras of dimension $d$ with Gorenstein parameter $l$.
\end{itemize} 
\end{theorem}
\begin{proof}
(1) $\Rightarrow$ (2). It follows from the proof of (3) $\Rightarrow$ (4) of Theorem \ref{supplement of equivalent definition of GAS} that
$$R\gHom_A(S,A)\cong S\otimes_{A_0} {}^1 A_0^\nu(l)[-d]$$
in $\D(\Gr A^e)$.
Note that $\nu$ induces an automorphism of $S$. So, $\gExt_A^i(S,A)=0$ for all $i\neq d$ and
$$\gExt_A^d(S,A)\cong {}^1 (\mathop{\oplus}\limits_{i=1}^n M_{r_i}(D_i))^\nu(l).$$

(2) $\Rightarrow$ (1). In this case $\mu$ and $\mu'$ is graded. 
By Propositions \ref{graded Ext(M,A)=Hom(M,D(A)} and \ref{graded Ext(M,A)=Hom(M,D(A) in Aop}, 
$$R\gHom_A(A_0,A)\cong D({}^\mu A_0)(l)[-d] \textrm{ and } R\gHom_{A^o}(A_0,A)\cong D(A_0^{\mu'})(l)[-d].$$

(2) $\Rightarrow$ (3). With the notations as above, there are isomorphisms in $\D(\Gr A^e)$:
$$R\Gamma_A(A)\cong I'^d[-d]\cong {}^1D(A)^{\mu}(l)[-d] \textrm{  
and } R\Gamma_{A^o}(A)=J'^d[-d]={}^{\mu'}D(A)^1(l)[-d].$$

(3) $\Rightarrow$ (2).  By Theorem \ref{another equivalent definition of GAS}, $A$ is an $\mathbb{N}$-graded AS-Gorenstein algebra.  
Since $R\Gamma_A(A)\cong {}^1D(A)^{\mu}(l)[-d]$, as the proof of Theorem \ref{another equivalent definition of GAS}, $\gExt_A^d(S,A)\cong D({}^\mu A\otimes_A S) \cong D(S)^\mu \cong {}^{\mu^{-1}} S^1$. 
Thus $A$ is a canonical $\mathbb{N}$-graded AS-Gorenstein algebra.
\end{proof}

\begin{theorem}\label{GAS-regular and AS over A0}
Let $A$ be a locally finite $\mathbb{N}$-graded algebra. Then the following are equivalent.
\begin{itemize}
\item[(1)\hphantom{'}] $A$ is AS-regular over $A_0$ of dimension $d$ with Gorenstein parameter $l$.
\item[$(1)^{\prime}$] $A^o$ is AS-regular over $A_0$ of dimension $d$ with Gorenstein parameter $l$.
\item[(2)\hphantom{'}] $A$ is a canonical $\mathbb{N}$-graded AS-regular algebra of dimension $d$ and $l_1=\cdots=l_n=l$.
\item[(3)\hphantom{'}] $A$ is ASF-regular algebra of dimension $d$ with Gorenstein parameter $l$.
\item[$(3)^{\prime}$] $A^o$ is ASF-regular algebra of dimension $d$ with Gorenstein parameter $l$.
\end{itemize} 
\end{theorem}
\begin{proof}
First we claim that  $A$ is Ext-finite in any case above.

If $A$ is AS-regular over $A_0$, or $A^o$ is AS-regular over $A_0$, then $A$ is Ext-finite by Proposition \ref{criterion for Ext-finite}.

If $A$ is $\mathbb{N}$-graded AS-regular, then $A$ is Ext-finite by Corollary \ref{GAS-regular is Ext-finite}. 

If $A$ is ASF-regular, then by an analogous argument as (2) $\Rightarrow$ (1) of Theorem \ref{another equivalent definition of GAS}, $\gExt_A^i(S,A)=0$ for $i\neq d$ and $\gExt_A^d(S,A)\cong \soc D(A)(l)$ which is finite-dimensional. 
By Proposition \ref{criterion for Ext-finite}, $A$ is Ext-finite. Similarly if $A^o$ is ASF-regular, then it is Ext-finite.

(1) $\Leftrightarrow (1)^{\prime}$. By \cite[Proposition 3.6]{MM}.

(1) + $(1)^{\prime} \Leftrightarrow$ (2) $\Rightarrow$ (3) or $(3)^{\prime}$. They are direct from Theorem \ref{GAS-Gorenstein and AS over A0}.

(3) $\Rightarrow$ (2). Since $\gExt_A^d(S,A)\neq 0$ and $\gldim A=d$, the injective dimension of ${}_AA$ is $d$. By Proposition \ref{dual property}, $\gExt^d_{A}(\gExt_{A^o}^d(S,A),A)\cong S\neq 0$, so the injective dimension of $A_A$ is also $d$. It follows from Theorem \ref{another equivalent definition of GAS} that $A$ is $\mathbb{N}$-graded AS-regular. Moreover, by Corollary \ref{gExt_A(S,A) cong gExt_Ao(S,A)}, $A$ is canonical.

$ (3)^{\prime}\Rightarrow$ (2). Similar to (3) $\Rightarrow$ (2).
\end{proof}

In \cite[Theorem 3.12]{MM}, the author proved that $A$ being AS-regular over $A_0$ implies $A$ being ASF-regular, and they left the converse to be open. Theorem \ref{GAS-regular and AS over A0} completes this equivalent characterization.

\begin{remark}\label{GAS-Gorenstein in RRZ}
There is another definition of generalized AS-Gorenstein algebra given in \cite[Definition 3.3]{RRZ} as in the following: a locally finite $\mathbb{N}$-graded algebra $A$ is generalized AS-Gorenstein if (1) $A$ has finite injective dimension $d$ on both sides; (2) $A$ is noetherian and satisfies the $\chi$-condition, and the functor $\Gamma_A$ has finite cohomological dimension; (3) $D(R^d\Gamma_A(A))\cong {}^{\mu}A^1(-l)$ for some graded automorphism $\mu$ of $A$. 

Note that if $A$ is Ext-finite, then (1) and (3) are equivalent to say that $A$ and $A^o$ are AS-Gorenstein over $A_0$. When $A$ is noetherian then (2) is a consequence of (1) and (3) by Corollary \ref{corollary of existence of balanced dualizing complex}.
\end{remark}

\section{Graded and Ungraded Calabi-Yau Algebras}\label{Graded and Ungraded Calabi-Yau Algebras}
As mentioned before, an $\mathbb{N}$-graded algebra $A$ is twisted Calabi-Yau if and only if $A$ is $\mathbb{N}$-graded AS-regular and $S=A/J$ is a separable algebra \cite[Theorem 5.15]{RR}. In this section, we show that when an $\mathbb{N}$-graded AS-regular algebra $A$ with $S=A/J$  separable is (skew) Calabi-Yau in either graded sense or ungraded sense.

Here is the definition of twisted Calabi-Yau algebras.
\begin{definition}\label{definition of CY of dimension d}
Let $A$ be an $\mathbb{N}$-graded algebra with $A_0$ finite-dimensional. If 
\begin{enumerate}
\item [(1)] $A$ is homologically smooth, i.e. as a graded $(A,A)$-bimodule, $A$ has a  finitely generated projective resolution of finite length; 
\item [(2)] there is an invertible graded $(A,A)$-bimodule $U$ such that
$$
\gExt_{A^e}^i(A,A^e)\cong\left\{
\begin{aligned}
&0,&i\neq d\\
&U, &i=d
\end{aligned}
\right.$$
\end{enumerate}
then  $A$ is called a graded twisted Calabi-Yau algebra of dimension $d$. 

Moreover if $U\cong {}^1A^\mu(l)$ as graded $(A,A)$-bimodules for some graded automorphism $\mu$ of $A$, then $A$ is called graded skew Calabi-Yau and $\mu$ is called a Nakayama automorphism. 
If $\mu$ is an inner automorphism, then $A$ is called graded Calabi-Yau.

Ungraded (twisted, skew) Calabi-Yau algebras can be defined similarly.
\end{definition}

As proved in \cite[Theorem 3.10]{RR}, $A$ is homologically smooth if and only if $S=A/J$ is separable, and either ${}_AS$ or $S_A$ possesses a finitely generated projective resolution of finite length.

\begin{lemma}\label{lemma about S separable}
Let $A$ be an $\mathbb{N}$-graded algebra with $A_0$ finite dimensional and $S=A/J$. Suppose $S$ is a separable algebra.
\begin{itemize}
    \item [(1)] Let $J(A^e)$ be the graded Jacobson radical of $A^e$. Then $A^e/J(A^e)\cong S^e$.
    \item [(2)] ${}_AS$ has a finitely generated graded projective resolution if and only if so does ${}_{A^e}A$.
\end{itemize} 
\end{lemma}
\begin{proof}
(1). By \cite[Lemma 3.7]{RR}.

(2). Let $P_\bullet$ be a minimal graded projective resolution of ${}_{A^e}A$ and $P'_\bullet$ be a minimal graded projective resolution of ${}_AS$. Then $P_\bullet\otimes_A S$ is a graded projective resolution of ${}_A S$. Since
$$S^e\otimes_{A^e} P_\bullet\cong S\otimes_A P_\bullet\otimes_A S=S\otimes_A (P_\bullet\otimes_A S),$$
$\Tor_i^{A^e}(S^e,A)\cong \Tor_i^A(S,S)$ for all $i \geqslant 0$. 

Observe that $\Tor_i^{A^e}(S^e,A)\cong S^e\otimes_{A^e} P_i$ and $\Tor_i^A(S,S)\cong S\otimes_A P'_i$. Hence $P_i$ is finitely generated as a graded $(A,A)$-bimodule if and only if $\Tor_i^{A^e}(S^e,A)$ is finite-dimensional, if and only if $\Tor_i^A(S,S)$ is finite-dimensional, if and only if $P'_i$ is finitely generated.
\end{proof}

The next proposition is related to Proposition \ref{Ext(A,Ae) and balanced dualizing complex}.

\begin{proposition}\label{Ext(A,Ae)/J and  balanced dualizing complex}
Let $A$ be an Ext-finite $\mathbb{N}$-graded algebra with $S=A/J$ separable. Suppose $A$ has finite injective dimension $d$ on both sides, 
$\gExt_{A^e}^i(A,A^e)=0$ for all $i\neq d$ and $\gExt^d_{A^e}(A,A^e)$ is an invertible graded $(A,A)$-bimodule. 
Set $U=\gExt_{A^e}^d(A,A^e)$ and $\omega=\gHom_{A^o}(D(R^d\Gamma_A(A)),A)$. 
Then
\begin{itemize}
    \item [(1)] $A$ is an $\mathbb{N}$-graded AS-Gorenstein algebra of dimension $d$.
    \item [(2)] $U/JU\cong \omega/J\omega\cong\gExt_A^d(S,A)\cong \gExt_{A^o}^d(S,A)$ as graded $(A,A)$-bimodules.
    \item [(3)] $U\cong \omega$  in $\Gr A$ and in $\Gr A^o$.  
\end{itemize}
\end{proposition}
\begin{proof}
(1) Let 
$\cdots\to P_{d+1} \xrightarrow[]{\partial_{d+1}} P_d \xrightarrow[]{\partial_d} \cdots \to  P_0 \xrightarrow[]{\partial_0} A \to 0$
be a minimal graded projective resolution of ${}_{A^e}A$. 
As $A$ is Ext-finite, $P_i$ is finitely generated as a graded $(A,A)$-bimodule for all $i$ by Lemma \ref{lemma about S separable}.
Note that the complex $\gHom_{A^e}(P_\bullet,A^e)$ is exact except the place $\gHom_{A^e}(P_d,A^e)$. 
Therefore for any $i>d$, the complex
$$ \gHom_{A^e}(P_{d},A^e)\to\cdots\to \gHom_{A^e}(P_{i-1},A^e)\xrightarrow[]{\partial_i^*}\gHom_{A^e}(P_i,A^e)\to \im\partial_{i+1}^*\to 0$$
is a part of a graded projective resolution of $\im \partial_{i+1}^*$ as a graded $(A,A)$-bimodule and as a graded $A$-module. 
Thus for $i\geqslant d+3$, $\Tor_1^A(S,\im \partial^*_{i+1})=H^{i-1}(S\otimes_A\gHom_{A^e}(P_\bullet,A^e))$.

By the finiteness of $P_\bullet$,
\begin{align*}
    S\otimes_A \gHom_{A^e}(P_\bullet,A^e)&\cong \gHom_{A^e}(P_\bullet,A\otimes S\otimes_A A)\\
    &\cong \gHom_{A^e}(P_\bullet,A\otimes S)\\
    &\cong \gHom_{A^e}(P_\bullet, A\otimes D(S))\\
    &\cong \gHom_{A^e}(P_\bullet,\gHom_k(S,A))\\
    &\cong \gHom_A(P_\bullet\otimes_A S,A).
\end{align*}
Note that $P_\bullet\otimes_A S$ is a graded projective resolution of ${}_AS$. 
Since the injective dimension of ${}_AA$ is $d$, $H^i(S\otimes_A \gHom_{A^e}(P_\bullet,A^e))\cong \gExt_A^i(S,A)=0$ for all $i>d$. 
Therefore for all $i\geqslant d+3$, $\Tor_1^A(S,\im \partial^*_{i+1})=0$. 
As a graded $A$-module, $\im \partial^*_{i+1}$ has a minimal graded projective resolution, denoted by $Q_\bullet$. Then $\Tor_1^A(S,\im \partial^*_{i+1})=S\otimes_A Q_1$
which implies that $Q_1=0$. 
Thus $\im \partial^*_{i+1}$ is projective as a graded $A$-module. 
Together with the fact that homology groups of $\gHom_{A^e}(P_\bullet,A^e)$ are either graded projective $A$-modules or $0$, it implies that $\gHom_{A^e}(P_\bullet, A^e)$ is split as a graded $A$-module complex. 
Hence
$$H^i(S\otimes_A \gHom_{A^e}(P_\bullet,A^e))\cong S\otimes_A H^i(\gHom_{A^e}(P_\bullet,A^e)),$$
and
$$\gExt_A^i(S,A)\cong S\otimes_A \gExt_{A^e}^i(A,A^e)\cong\left\{
\begin{aligned}
&0,&i\neq d,\\
&S\otimes_A U, &i=d.
\end{aligned}
\right.$$
By Theorem \ref{equvalent definitions of AS-Gorenstein}, $A$ is an $\mathbb{N}$-graded AS-Gorenstein algebra of dimension $d$.

(2) It follows from (1) that $\gExt_A^d(S,A)\cong S\otimes_A U\cong U/JU$.
On the other hand, Corollary \ref{gExt_A(S,A) cong gExt_Ao(S,A)} shows that $\gExt_A^d(S,A)\cong \omega/J\omega$.
So $U/JU\cong \omega/J\omega$ as graded $(S,S)$-bimodules.

(3) Since $U$ and $\omega$ are finitely generated projective on both sides, by Nakayama's Lemma, $U\cong \omega$ in $\Gr A$ and in $\Gr A^o$.
\end{proof}

Finally we discuss when an $\mathbb{N}$-graded AS-regular algebra becomes graded (skew) Calabi-Yau. 

\begin{theorem}\label{when twisted CY is skew CY}
Suppose $A$ is an $\mathbb{N}$-graded AS-regular algebra of dimension $d$ with $S=A/J$ being a separable algebra, or equivalently $A$ is a twisted Calabi-Yau algebra of dimension $d$. Then
\begin{itemize}
\item[(1)] $A$ is an ungraded skew Calabi-Yau algebra if and only if $A$ is canonical;
\item[(2)] $A$ is a graded skew Calabi-Yau algebra if and only if $A$ is canonical and $l_1=\cdots=l_n=l$. 
\end{itemize}
\end{theorem}
\begin{proof}
(1) With the notations as above, let $U=\gExt^d_{A^e}(A,A^e)$. By Proposition \ref{Ext(A,Ae)/J and  balanced dualizing complex}, $U/JU\cong\gExt_A^d(S,A)$. 

``$\Rightarrow$" Suppose $U\cong {}^1A^{\varphi}$ as ungraded $(A,A)$-bimodules where $\varphi$ is an automorphism of $A$. 
By \cite[Lemma 2.11]{RR}, $U\otimes_A S=S\otimes_A U$ which implies that $\varphi(J)=J$. 
Then $\varphi$ induces an automorphism of $S$, denoted by $\delta$, and $\gExt_A^d(S,A)\cong U/JU\cong {}^1S^\delta$ as ungraded $(A,A)$-bimodules. So $A$ is canonical.

``$\Leftarrow$" If $A$ is canonical, then by Proposition \ref{notations of canonical commonly graded AS}, $\gExt_A^d(S,A)\cong\mathop{\oplus}\limits_{i=1}^n {}^1 M_{r_i}(D_i)^\delta(l_i)$. 
Hence $U\cong \mathop{\oplus}\limits_{i=1}^n Ae_i(l_i)$ as graded $A$-modules and $U\cong A$ as ungraded $A$-modules. Because $U$ is invertible, there is an automorphism $\varphi$ of $A$ such that $U\cong {}^1 A^\varphi$ as ungraded $(A,A)$-bimodules.

Since $A$ has a finitely generated projective resolution with finite length as a graded $(A,A)$-bimodule, 
$\gExt^d_{A^e}(A,A^e)= \Ext^d_{A^e}(A,A^e)$ and $\Ext_{A^e}^i(A,A^e)=0$ for all $i\neq d$. 
Hence
$$\Ext^d_{A^e}(A,A^e)\cong U\cong {}^1 A^\varphi$$
as ungraded $(A,A)$-bimodules. Finally by \cite[Lemma 2.4]{RR}, if $A$ is graded homologically smooth, then $A$ is ungraded homologically smooth. In conclusion, $A$ is an ungraded skew Calabi-Yau algebra.

(2) ``$\Rightarrow$" By Proposition \ref{Ext(A,Ae)/J and  balanced dualizing complex}.

``$\Leftarrow$" By (1), $U\cong {}^1A^{\varphi}$ where $\varphi$ is an automorphism of $A$. Moreover $\varphi$ is graded as $l_1=\cdots=l_n=l$. It follows that $U\cong {}^1A^\varphi(l)$ as graded $(A,A)$-bimodules. Thus $A$ is skew Calabi-Yau.
\end{proof}

In conclusion, $A$ is AS-regular over $A_0$ with $S=A/J$ being a separable algebra if and only if $A$ is a graded skew Calabi-Yau algebra. 

The following is an example of graded twisted Calabi-Yau algebra which is ungraded skew Calabi-Yau but not graded skew Calabi-Yau.

\begin{example}\cite[Example 6.7]{RR2}\label{example: graded twisted CY is not graded skew CY}
Let $Q$ be the quiver with vertices $\{1,\cdots, n\}$ and  arrows  $\alpha_i:i-1\to i$ for $2\leqslant i\leqslant n$ and $\alpha_1:n\to 1$, where $\alpha_i$ has weight $l_i$. Let $A=kQ$ and assume at least one $l_i\neq 0$. Then $A$ is graded twisted Calabi-Yau of dimension $1$ and the Nakayama's automorphism $\mu$ satisfies $\mu(e_i)=e_{i+1}$. Moreover $\mu$ is graded if and only if all $l_i$ are equal. Thus $A$ is graded skew Calabi-Yau if and only if all $l_i$ are equal. To construct a graded twisted Calabi-Yau algebra but not graded skew Calabi-Yau, we only need to choose different $l_i$. 
\end{example}

\section{Twisted Derived-Calabi-Yau Algebras}\label{commonly graded AS-Regular Algebras and Twisted Derived-Calabi-Yau Algebras}
In this section we first introduce the notion of derived-CY algebras, then we 
prove that $A$ is a noetherian commonly graded AS-regular algebra of dimension $d$ if and only if $A$ is a twisted derived-CY algebra of dimension $d$. 
Consequently, for any $\mathbb{N}$-graded algebra $A$, $A$ is a twisted Calabi-Yau algebra of dimension $d$ if and only if $A$ is a twisted derived-CY algebra of dimension $d$ and $S=A/J$ is a separable algebra. 

Let $\proj A$ be the full subcategory of $\Gr A$ consisting of finitely generated projective $A$-modules. Next definition is motivated by \cite{IR}.

\begin{definition}\label{definition of dCY algebras}
Let $A$ be a noetherian commonly graded algebra and $U$ an invertible graded $(A,A)$-bimodule. If for any $X, Y\in \D^b(\fd A)$, there exists a functorial isomorphism 
\begin{equation}\label{equation}
    \Hom_{\D(\Gr A)}(X, Y)\cong D\Hom_{\D(\Gr A)}(Y,U\otimes_AX[d]),
\end{equation}
then $A$ is called twisted derived-Calabi-Yau (twisted derived-CY for short) of dimension $d$.

If the functorial isomorphism \eqref{equation} holds for any $X\in \K^b(\proj A)$ and $Y\in \D^b(\fd A)$, then $A$ is called twisted derived-singular-Calabi-Yau (twisted derived-sCY for short). 

If further $U\cong {}^\mu A^1(-l)$ as $(A,A)$-bimodules, where $\mu$ is a graded automorphism of $A$, then $A$ is called skew derived-CY (resp. skew derived-CY). Moreover if $\mu$ is an inner automorphism then $A$ is called derived-CY (resp. derived-sCY).
\end{definition}

Let $A$ be a noetherian commonly graded algebra. Let $\D^b_{\fd}(\Gr A)$ be the full subcategory of $\D^b(\Gr A)$ consisting of the complexes whose cohomologies are finite-dimensional.

\begin{lemma}\label{Db(fd A) Hom finite}
Suppose $A$ is a noetherian commonly graded algebra. Then the triangulated category $\D^b(\fd A)$ is equivalent to $\D^b_{\fd}(\Gr A)$, and $\D^b(\fd A)$ is Hom-finite.
\end{lemma}
\begin{proof}
Recall that $\D^b_{\fg}(\Gr A)$ is the full subcategory of $\D^b(\Gr A)$ consisting of the complexes whose cohomologies are finitely generated. Let $\D^b_{\fd}(\gr A)$ be the full subcategory of $\D^b(\gr A)$ consisting of complexes whose cohomologies are finite-dimensional. By \cite[Proposition 1.7.11]{KS} and its dual version, $\D^b(\fd A)$ is equivalent to $\D^b_{\fd}(\gr A)$ and  $\D^b(\gr A)$ is equivalent to $\D^b_{\fg}(\Gr A)$. Let $F$ be the composition:
$$\D^b(\fd A)\xrightarrow{\cong} \D^b_{\fd}(\gr A)\hookrightarrow \D^b(\gr A)\xrightarrow{\cong} \D^b_{\fg}(\Gr A)\hookrightarrow \D^b(\Gr A).$$
Then $F$ is fully faithful. Since $F(\D^b(\fd A))\subseteq \D^b_{\fd}(\Gr A)$, $F$ can be regarded as a fully faithful functor from $\D^b(\fd A)$ to $\D^b_{\fd}(\Gr A)$. 

For any $X\in \D^b_{\fd}(\Gr A)$, $X$ is quasi-isomorphic to a complex $X'\in \D^b(\gr A)$. Since the cohomologies of $X$ are finite-dimensional, $X'\in  \D^b_{\fd}(\gr A)$. Hence $X'$ is quasi-isomorphic to a complex in $\D^b(\fd A)$. Therefore $F:\D^b(\fd A)\to \D^b_{\fd}(\Gr A)$ is an equivalent functor. 

Let $X,Y\in \D^b(\fd A)$. Since $A$ is noetherian, $X$ is quasi-isomorphic to a complex $P \in \K^-(\proj A)$. So $\Hom_{\D^b(\fd A)}(X,Y)\cong \Hom_{\K(\Gr A)}(P,Y) =\textrm{H}^0(\Hom_{\Gr A}^\bullet(P, Y))$. Since $\textrm{H}^0(\Hom_{\Gr A}^\bullet(P, Y))$ is finite-dimensional, it follows that $\D^b(\fd A)$ is Hom-finite.
\end{proof}

The following proposition is a commonly graded version of \cite[Proposition 2.6]{IR}. Note that if $A$ is finitely generated as a graded $k$-algebra, then $\underleftarrow{\Lim} M/J^iM\cong M$ for any finitely generated graded module ${}_AM$. For $P\in \K^b(\proj A)$, let $P^{(i)}$ be the complex:
$$\cdots\to P^{-1}/J^iP^{-1}\to P^0/J^iP^0\to P^1/J^iP^1\to \cdots.$$

\begin{proposition}\label{limit of derived hom set}
 Let $P\in \K^b(\proj A)$ and $U$ be an invertible graded $(A,A)$-bimodule. With the above notations, then 
\begin{itemize}
\item[(1)] $\underleftarrow{\Lim}\gHom_{\D(\Gr A)}(X,U\otimes_A P^{(i)})\cong \gHom_{\D(\Gr A)}(X,U\otimes_A P)$ for any $X\in \D^-(\gr A)$.
\item[(2)] $\underrightarrow{\Lim}\gHom_{\D(\Gr A)}(P^{(i)},X) \cong \gHom_{\D(\Gr A)}(P,X)$ for any $X\in \D^b(\fd A)$.
\end{itemize}
\end{proposition}
\begin{proof}
To convince the readers, we give a proof here. 

(1) Suppose $Q$ is a finitely generated graded projective resolution of $X$. Since $U$ is an invertible bimodule, 
$U\otimes_A J^iP^t\cong J^i(U\otimes_A P^t)$ for any $t$ and $i$.
Thus  
$$\underleftarrow{\Lim}U\otimes_A P^t/U \otimes J^iP^t\cong \underleftarrow{\Lim} (U\otimes_A P^t)/J^i(U\otimes_A P^t)\cong U\otimes_A P^t$$
 and 
$$\underleftarrow{\Lim}\gHom_A(Q^s,(U\otimes_A P^{(i)})^t)\cong \gHom_A(Q^s,\underleftarrow{\Lim} U\otimes_A P^t/J^iP^t)\cong \gHom_A(Q^s,U\otimes_A P^t).$$

So there is an isomorphism of complexes

$$\xymatrix{
	\cdots     \prod_s\gHom_A(Q^s,U\otimes_A P^{s+n}) \ar[d]^{\cong} \ar[r]^{\partial_{\Hom}} & \prod_s\gHom_A(Q^s,U\otimes_A P^{s+n+1}) \ar[d]^{\cong} \cdots \\
\cdots   \prod_s\gHom_A(Q^s, \underleftarrow{\Lim}  \frac{U\otimes_AP^{s+n}}{J^i(U\otimes_A P^{s+n})})\ar[d]^{\cong}  \ar[r]^{\prod_s \partial_{\Hom}} & \prod_s\gHom_A(Q^s, \underleftarrow{\Lim}  \frac{U\otimes_AP^{s+n+1}}{J^i(U\otimes_A P^{s+n+1})}) \ar[d]^{\cong} \cdots \\
\cdots   \prod_s\gHom_A(Q^s, \underleftarrow{\Lim}  \frac{U\otimes_AP^{s+n}}{U\otimes_A J^iP^{s+n}})\ar[d]^{\cong}  \ar[r]^{\prod_s\partial_{\Hom}} & \prod_s\gHom_A(Q^s, \underleftarrow{\Lim}  \frac{U\otimes_AP^{s+n+1}}{U\otimes_A J^iP^{s+n+1}}) \ar[d]^{\cong}  \cdots \\
\cdots   \underleftarrow{\Lim} \prod_s\gHom_A(Q^s,  U\otimes_A\frac{ P^{s+n}}{J^iP^{s+n}}) \ar[r]^{\underleftarrow{\Lim} \prod_s\partial_{\Hom}} & \underleftarrow{\Lim}\prod_s\gHom_A(Q^s,  U\otimes_A \frac{P^{s+n+1}}{J^iP^{s+n+1}})  \cdots  }
$$

Since $P^{(i)}$ is a bounded complex consisting of finite-dimensional modules, every item in complex $\prod_s\gHom_A(Q^s,(U\otimes_A P^{(i)})^{s+\bullet})$ is finite-dimensional.  Thus, by Mittag-Leffler condition, the inverse limit commutes with cohomology. Therefore 
$$\gHom_{\D(\Gr A)}(X,U\otimes_A P)=H^0\gHom_A(Q,U\otimes_A P)\cong \underleftarrow{\Lim}\gHom_{\D(\Gr A)}(X,U\otimes_A P^{(i)}). $$ 

(2) Let $p:P\to P^{(i)}$ and  $\varphi$ : $\underrightarrow{\Lim}\gHom_{\D(\Gr A)}(P^{(i)},X)\to\gHom_{\D(\Gr A)}(P,X)$ be the canonical maps. 
For any $s^{-1}g\in \gHom_{\D(\Gr A)}(P^{(i)},X)$ where $s:X\to X'$ is a quasi-isomorphism and $g:P^{(i)}\to X'$ is a complex 
cochain map, $\varphi(s^{-1}g)=s^{-1}gp$. If $\varphi(s^{-1}g)=0$, then $gp$ is null-homotopic. By Lemma \ref{Db(fd A) Hom finite}, $X'$ is in $\D^b(\fd A)$. Hence there is an integer $n$ such that $J^nX'=0$. 
Then for any $t$ and any map $P^t\to X'^{t-1}$, this map factors through $P^t/J^nP^t$. 
Therefore the composition $P^{(n)}\to P^{(i)}\xrightarrow{g} X'$ is null-homotopic, which implies $s^{-1}g$ is zero in $\underrightarrow{\Lim}\gHom_{\D(\Gr A)}(P^{(i)},X)$, so $\varphi$ is injective.

For any $f\in \gHom_{\D(\Gr A)}(P,X)=\gHom_{\K(\Gr A)}(P,X)$, since $X\in \K^b(\fd A)$, $f$ factors through $P^{(i)}$ for some $i$, which means there exists $g:P^{(i)}\to X$ such that $f$ is the composition of $P\xrightarrow{p} P^{(i)}\xrightarrow{g} X$. Then $\varphi(g)=f$. So $\varphi$ is surjective.
\end{proof}

Next proposition is proved in \cite{IR} for module-finite algebras.
\begin{proposition}\label{d-CY and d-sCY}
\begin{itemize}
\item[(1)] Derived-CY (resp. derived-sCY) algebras are closed under derived equivalence.
\item[(2)] $A$ is twisted (resp. skew, none) derived-CY of dimension $d$ if and only if so is $A^o$.
\item[(3)] Any twisted derived-CY algebra of dimension $d$ has global dimension $d$.
\item[(4)] $A$ is twisted (resp. skew) derived-CY of dimension $d$ if and only if $A$ is twisted (resp. skew) derived-sCY of dimension $d$ and $A$ has finite global dimension.
\item[(5)] If $A$ is twisted (resp. skew) derived-CY of dimension $d$, then there exists a functorial isomorphism $$\gHom_{\D(\Gr A)}(X, U\otimes_A Y[d])\cong D\gHom_{\D(\Gr A)}(Y,X)$$ for any $X\in \D^b(\fd A)$ and $Y\in \D^b(\gr A)$.
\end{itemize}
\end{proposition}
\begin{proof}
(1) Let $A$ and $B$ be two noetherian commonly graded algebras. Note that $\K^b(\proj A)$ consists of the compact objects of $\D^b(\Gr A)$, and $X$ is in $\D^b(\fd A)$ if and only if $\oplus_i\gHom_{\D(\Gr A)}(Y,X[i])$ is finite-dimensional for any $Y\in \K^b(\proj A)$. 
So any triangulated equivalence $\D^b(\Gr A)\to \D^b(\Gr B)$ induces triangulated equivalences $\K^b(\proj A)\to \K^b(\proj B)$ and $\D^b(\fd A)\to \D^b(\fd B)$. 
It follows that derived-CY (resp. derived-sCY) algebras are closed under derived equivalence.

(2) Note that $D=\gHom_k(-,k)$ induces a duality $\D^b(\fd A)\to \D^b(\fd A^o)$. 
For any invertible graded $(A,A)$-bimodule $U$ and $X\in \D^b(\fd A)$, $D(U\otimes_A X)\cong D(X)\otimes_A U^{-1}$  where $U^{-1}=\gHom_{A^o}(U,A)$. 
So $A$ is twisted (resp. skew, none) derived-CY of dimension $d$ if and only if so is $A^o$.

(3) Let $A$ be a twisted derived-CY algebra of dimension $d$ and $S=A/J$. By \cite[Proposition 3.18]{RR}, the global dimension of $A$ is $\max\{i\mid \gExt_A^i(S,S)\neq 0\}$. Since $S=A/J$ belongs to $\fd A$, $\gExt_A^i(S,S)\cong D\gExt^{d-i}_A(S,U\otimes_A S)=0$ for $i>d$,
and $\gExt_A^d(S,S)\cong D\gHom_A(S,U\otimes_A S)\cong D(U\otimes_A S)\neq 0$. It follows that the global dimension of $A$ is $d$. 

(4) We only prove the twisted case. If $A$ has finite global dimension, then $\D^b(\fd A)$ is a full subcategory of $\K^b(\proj A)$. Hence in this case,  $A$ being twisted derived-sCY implies $A$ being twisted derived-CY.

If $A$ is twisted derived-CY of dimension $d$, taking $X\in \K^b(\proj A)$ and $Y\in\D^b(\fd A) $, then $X^{(i)}\in \D^b(\fd A)$ and 
$$\gHom_{\D(\Gr A)}(X^{(i)},Y)\cong D\gHom_{\D(\Gr A)}(Y,U\otimes_A X^{(i)}[d]).$$ 
It follows from Proposition \ref{limit of derived hom set} that
$$\gHom_{\D(\Gr A)}(X,Y)\cong D\gHom_{\D(\Gr A)}(Y,U\otimes_A X[d]).$$
So $A$ is derived-sCY of dimension $d$.

(5) Since $A$ has finite global dimension, $\D^b(\gr A)=\K^b(\proj A)$. Thus (5) holds.
\end{proof}

Next proposition is a generalization of the Serre duality.
\begin{proposition}\label{generalized Serre duality}
Suppose $\Gamma_A$ has finite cohomological dimension. Let $W=D(R\Gamma_A(A))$. Then for any $X\in \K^b(\proj A)$ and $Y\in\D^b(\fd A) $, we have functorial isomorphism
$$R\!\gHom_A(X,Y)\cong DR\!\gHom_A(Y,W\otimes_A^L X).$$
\end{proposition}
\begin{proof}
By Theorem \ref{local duality} and \cite[Lemma 4.4]{VdB}, $R\!\gHom_A(Y,W)\cong D(R\Gamma_A(Y))\cong Y$. 
Since $W \in \D^b(\Gr A^e)$ and $X\in \K^b(\proj A)$, $W\otimes_A^L X=W\otimes_A X$ is bounded. Let $Q$ be a finitely generated graded projective resolution of $Y$. Then  $\gHom_A(Q,W\otimes_A X)\cong \gHom_A(Q,W)\otimes_A X$. 
Hence \begin{align*}
R\!\gHom_A(Y,W\otimes_A^L X)&=\gHom_A(Q,W\otimes_A X)\\
&\cong \gHom_A(Q,W)\otimes_A X\\
&=R\!\gHom_A(Y,W)\otimes_A X\\
&\cong D(R\Gamma_A(Y))\otimes_A X\\
&\cong D(Y)\otimes_A X.
\end{align*} 

On the other hand,
\begin{align*}
DR\!\gHom_A(X,Y)&\cong D(\gHom_A(X,A)\otimes_A Y)\\
&\cong \gHom_{A^o}(\gHom_A(X,A),D(Y))\\
&\cong D(Y)\otimes_A \gHom_{A^o}(\gHom_A(X,A),A)\\
&\cong D(Y)\otimes_A X.
\end{align*}
It follows that $R\!\gHom_A(X,Y)\cong DR\!\gHom_A(Y,W\otimes_A^L X)$.
\end{proof}

For the invertible bimodule in (\ref{equation}), we have the following characterization.

\begin{proposition}\label{invertible bimodule of der-sCY}
    Suppose $A$ is twisted derived-sCY of dimension $d$. Then $U\cong D(R^d\Gamma_A(A))$ as graded $(A,A)$-bimodules.
\end{proposition}
\begin{proof}
    Taking $X=U^{-1}:=\gHom_A(U,A)$ and $Y=A/J^i$ in (\ref{equation}) yields that
    \[D(\gHom_A(U^{-1},A/J^i))\cong \gExt_A^d(A/J^i,A).\]
    Applying direct limit to this isomorphism and by Lemma \ref{Hom iso}, we have
    \[
        D(A)\otimes_A U^{-1}\cong R^d\Gamma_A(A)
    \]
    which implies that $U\cong D(R^d\Gamma_A(A))$. Since (\ref{equation}) is functorial, all the isomorphisms are $(A,A)$-bimodule morphisms.
\end{proof}

Now we describe the relations between the commonly graded AS-Gorenstein algebras and the twisted derived-sCY algebras. 
\begin{theorem}\label{GAS-Gorenstein and d-sCY}
Suppose $A$ is a noetherian commonly graded algebra. Then the following are equivalent.
\begin{itemize}
\item[(1)] $A$ is a commonly graded AS-Gorenstein algebra of dimension $d$;
\item[(2)] $A$ is twisted derived-sCY of dimension $d$ and admits a balanced dualizing complex;
\item[(3)] $A$ is twisted derived-sCY of dimension $d$ and $A$ has finite injective dimension $d$ on both sides;
\item[(4)] $A$ admits a balanced dualizing complex, and there exits an invertible graded $(A,A)$-bimodule $U$ such that for any $X\in \fd A$, $\gExt^d_A(X,U)\cong D(X)$ as graded $A^o$-modules and $\gExt^i_A(X,U)=0$ for all $i\neq d$;
\item[(5)] $A$ has finite injective dimension $d$ on both sides, and there exits an invertible graded $(A,A)$-bimodule $U$ such that for any $X\in \fd A$, $\gExt^d_A(X,U)\cong D(X)$ as graded $A^o$-modules and $\gExt^i_A(X,U)=0$ for all $i\neq d$;
\item[(6)] $A$ admits a balanced dualizing complex $R$ such that $R\cong U[d] \in \D(\Gr A^e)$ for some invertible graded $(A,A)$-bimodule $U$.
\end{itemize}
\end{theorem}
\begin{proof}
(1) $\Rightarrow$ (2).
By Theorem \ref{existence of baalanced dualizing complex of GAS},
$A$ admits a balanced dualizing complex $R$ given by $R=D(R\Gamma_A(A))\cong D(R^d\Gamma_A(A))[d]$ and $D(R^d\Gamma_A(A))$ is an invertible graded $(A,A)$-bimodule. 
Corollary \ref{corollary of existence of balanced dualizing complex} deduces that $\Gamma_A$ has finite cohomological dimension. 
It follows from Proposition \ref{generalized Serre duality} that $A$ is twisted derived-sCY of dimension $d$.

(2) $\Rightarrow$ (3). Proposition \ref{invertible bimodule of der-sCY} shows that $U[d]$ is a balanced dualizing complex of $A$. By Theorem \ref{local duality}, the injective dimension of $U$ is $d$ on both sides. Since $U$ is invertible, $A$ shares the same injective dimension as $U$.

(2) $\Rightarrow$ (4) and 
(3) $\Rightarrow$ (5) are obvious.

(5) $\Rightarrow$ (1). Let $U^{-1}=\gHom_A(U,A)$. Then $\gExt^i_A(X,U)\cong \gExt_A^i(X,A)\otimes_A U$. Thus 
$$\gExt_A^d(S,A)\cong\gExt_A^d(S,U)\otimes_A U^{-1}\cong D(S)\otimes_A U^{-1}\cong S\otimes_A U^{-1}$$
 as graded $A^o$-modules. By assumption, $\gExt_A^i(S,A)=0$ for all $i\neq d$. 
 Since $S\otimes_A U^{-1}$ is an invertible graded $(S,S)$-bimodule, $A$ is a commonly graded AS-Gorenstein algebra of dimension $d$ by Theorem \ref{equvalent definitions of AS-Gorenstein}.

(4) $\Rightarrow$ (6). Similarly to the proof in (5) $\Rightarrow$ (1), $\gExt_A^d(S,A)\cong S\otimes_A U^{-1}$ as graded $A^o$-modules and $\gExt_A^i(S,A)=0$ for all $i\neq d$. It follows that $R\Gamma_A(A)\cong R^d\Gamma_A(A)[-d]$.

By a similar argument as previous sections, we see that $D(R^d\Gamma_A(A))$ is a finitely generated graded projective generator in $\Gr A^o$. 
Since $R:=D(R\Gamma(A))$ is a balanced dualizing complex of $A$, $\gHom_{A^o}(D(R^d\Gamma_A(A),D(R^d\Gamma_A(A)))\cong R\gHom_{A^o}(R,R)\cong A$ in $\D(\Gr A^e)$.  
It follows that $D(R^d\Gamma_A(A))$ is an invertible graded $(A,A)$-bimodule.

(6) $\Rightarrow$ (1). By Theorem \ref{another equivalent definition of GAS}, it is left to show that $A$ has finite injective dimension $d$ on both sides.
It follows from Theorem \ref{existence of balanced dualizing complex} that $U[d]\cong D(R^d\Gamma_A(A))[d]\cong D(R\Gamma_A(A))$.
So,
by Theorem \ref{local duality}, the injective dimension of ${}_AU$ is $d$. Dually, the injective dimension of $U_A$ is also $d$. 
Since $U$ is invertible, $A$ has the same injective dimension as $U$ on both sides.
\end{proof}

Now we come to the main result of this section.
\begin{theorem}\label{GAS-regular and d-CY}
Suppose $A$ is a noetherian commonly graded algebra. Then
\begin{itemize}
\item[(1)] $A$ is a commonly graded AS-regular algebra of dimension $d$ if and only if $A$ is a twisted derived-CY algebra of dimension $d$.
\item[(2)] When $A$ is $\mathbb{N}$-graded, $A$ is a twisted (resp. skew) Calabi-Yau algebra of dimension $d$ if and only if $A$ is a twisted (resp. skew) derived-CY algebra of dimension $d$ with $S=A/J$ being a separable algebra. 
\end{itemize}
\end{theorem}
\begin{proof}
(1) It follows from Proposition \ref{d-CY and d-sCY} and Theorem \ref{GAS-Gorenstein and d-sCY}.

(2) The ``twisted" case follows from (1) and \cite[Theorem 5.15]{RR}. The ``skew" case follows from (1) and Proposition \ref{Ext(A,Ae) and balanced dualizing complex}. 
\end{proof}

\section{Auslander-Buchsbaum Formula and Bass Theorem}\label{Auslander-Buchsbaum Formula and Bass Theorem}
In this section, we examine the conditions under which the Auslander-Buchsbaum formula, the Bass theorem, and the No-Hole theorem hold for a noetherian commonly graded algebra $A$ possessing a balanced dualizing complex. 
We give a series of equivalent conditions, denoted as $\mathbf{(C)}$ and $\mathbf{(C^o)}$, and subsequently employ these conditions to investigate noetherian commonly graded AS-Gorenstein algebras.

For any complex $X\in \D(\Gr A)$, let
$$\sup X:=\sup\{i\mid H^i(X)\neq 0\} \text{ and } \inf X:=\inf\{i\mid H^i(X)\neq 0\}.$$

    Next lemma is proved in \cite{Jo1} and \cite{WZ} in the case that $A$ is connected graded and ungraded local respectively.

\begin{lemma}\label{depth, sup, projective and injective dim}
Let $A$ be a noetherian commonly graded algebra with a balanced dualizing complex $R$ and $S=A/J$. Then, for any $0 \neq X\in \D^+_{\fg}(\Gr A)$,
\begin{itemize}
    \item [(1)] $\depth_A X=-\sup R\gHom_A(X,R)$.
    \item [(2)] $\depth_A X<\infty$.
    \item [(3)] $\idim_A X=\sup\{i\mid \gExt_A^i(S,X)\neq 0\}=\pdim_{A^o}R\gHom_A(X,R).$
\end{itemize}
\end{lemma}
\begin{proof}
     (1) and (2) By Lemma \ref{dualizing complex induces a duality}, $R\gHom_A(X,R)\in \D^-_{\fg}(\Gr A^o)$. 
     Let $P$ be a minimal graded projective resolution of $R\gHom_A(X,R)$.
     Then it follows from Lemma \ref{pdim and idim and resolution} that 
     $$\pdim_{A^o}R\gHom_A(X,R)=-\inf \{i\mid P^i\neq 0\}.$$
    Lemma \ref{dualizing complex induces a duality}, Theorems \ref{local duality} and \ref{existence of balanced dualizing complex} yield that
    \begin{align*}
        R\gHom_A(S,X)&\cong \gHom_{A^o}(R\gHom_A(X,R),R\gHom_A(S,R))\\
        &\cong R\gHom_{A^o}(R\gHom_A(X,R),D(R\Gamma_A(S)))\\
        &\cong \gHom_{A^o}(P,S).
    \end{align*}
    Since $P$ is minimal, the differentials in $\gHom_{A^o}(P,S)$ are zero. Hence
    $$\depth_A X=\inf R\gHom_A(S,X)=\inf \gHom_{A^o}(P,S)=-\sup \{i\mid P^i\neq 0\}=-\sup R\gHom_A(X,R)$$
    and
    $$\sup R\gHom_A(S,X)=\sup \gHom_{A^o}(P,S)=-\inf \{i\mid P^i\neq 0\}=\pdim_{A^o}R\gHom_A(X,R).$$

    Since $R\gHom_A(X,R)\in \D^-_{\fg}(\Gr A^o)$, $-\sup R\gHom_A(X,R)$ is finite. Thus $\depth_A X<\infty$.

    (3) By Lemma \ref{dualizing complex induces a duality}, $\idim_A X=\infty$ if and only if $\pdim_{A^o}R\gHom_A(X,R)=\infty$.

    If $\idim_AX<\infty$, then $\pdim_{A^o}R\gHom_A(X,R)<\infty$. Note that
    \begin{align*}
        \idim_AX&=\sup\limits_{M\in\gr A}\{\sup R\gHom_A(M,X)\}\\
        &=\sup\limits_{M\in\gr A}\{\sup R\gHom_{A^o}(R\gHom_A(X,R),R\gHom_A(M,R))\}\\
        &=\sup\limits_{M\in\gr A}\{\sup \gHom_{A^o}(P,R\gHom_A(M,R))\}.
    \end{align*}
    Suppose $\pdim_{A^o}R\gHom_A(X,R)=n$. Then $\inf P=-n$. Since $\idim_AR=0$, $\sup R\gHom_A(M,R)\leqslant 0$. It follows that
    $$\sup\limits_{M\in\gr A}\{\sup \gHom_{A^o}(P,R\gHom_A(M,R))\}\leqslant -\inf P=n.$$
    Hence
    $$\idim_AX\leqslant \pdim_{A^o} R\gHom_A(X,R)=\sup\{i\mid \gExt_A^i(S,X)\neq 0\}\leqslant \idim_A X.$$
    So $\idim_AX= \pdim_{A^o} R\gHom_A(X,R)=\sup\{i\mid \gExt_A^i(S,X)\neq 0\}$.
\end{proof}

For any Noetherian commonly graded algebra that possesses a balanced dualizing complex, we show next that the Auslander-Buchsbaum formula holds if and only if the Bass Theorem holds.

Let $\tilde{I}(A)$ (resp. $\tilde{P}(A)$) be the subcategory of finitely generated graded $A$-modules with finite injective (resp. projective) dimension.

\begin{theorem}\label{AB formula and Bass thm for locally finite alg}
    Let $A$ be a commonly graded algebra with a balanced dualizing complex $R$. 
    Suppose $0\to A\to I^0\to \cdots \to I^d\to \cdots$ is a minimal graded injective resolution of ${}_AA$. 
    If $\depth_A A=d$, then the following are equivalent.
    \begin{itemize}
        \item [C1\hphantom{'}] (Bass theorem) If $0\neq X\in \D^b_{\fg}(\Gr A)$ and $\idim_AX<\infty$, then $\idim_AX=d+\sup X$.
        \item [C2\hphantom{'}] (Auslander-Buchsbaum Formula) If $0\neq Y\in \D^b_{\fg}(\Gr A^o)$ and $\pdim_{A^o}Y<\infty$, then $$\pdim_{A^o}Y+\depth_{A^o} Y=d.$$
        \item [C$2^{\prime}$] $\pdim_{A^o}N+\depth_{A^o}N=d$, for any $N\in \tilde{P}(A^o)$.
        \item [C3\hphantom{'}] $\min\{i\mid \gExt_A^i(M,A)\neq 0\}\leqslant d$ for any $M\in \gr A$.
        \item [C4\hphantom{'}] $\gExt_A^d(M,A)\neq 0$ for every graded simple $A$-module $M$.
        \item [C5\hphantom{'}] $\soc_A I^d$ is a projective generator in $\gr S$.
        \item [C6\hphantom{'}] $\soc_A R^d\Gamma_A(A)$ is a projective generator in $\gr S$.
        \item [C7\hphantom{'}] $\Omega'/\Omega' J$ is a projective generator in $\gr S^o$ where $\Omega'=D(R^d\Gamma_A(A))$.
    \end{itemize}
\end{theorem}
\begin{proof}
    C1 $\Rightarrow$ C2. Let $X=R\gHom_{A^o}(Y,R)$. Then $\idim_AX<\infty$. Thus
    $\idim_AX=d+\sup X$.
    By Lemma \ref{depth, sup, projective and injective dim}, $\sup X=-\depth_{A^o} Y$ and $\idim_AX=\pdim_{A^o}R\gHom_A(X,R)=\pdim_{A^o} Y$.
    So $\pdim_{A^o}Y+\depth_{A^o} Y=d$.

    C2 $\Rightarrow$ C$2^{\prime}$. Clearly.

    C$2^{\prime}$ $\Rightarrow$ C3. Suppose on the contrary that  $\min\{i\mid \gExt_A^i(M,A)\neq 0\}>d$ for some $M\in\gr A$. Let $\cdots \to P_{d+1}\xrightarrow[]{\partial} P_d\to \cdots \to P_0\to M\to 0$ be a minimal graded projective resolution of $M$.
    Then
    $$0\to \gHom_A(P_0,A)\to \cdots\to \gHom_A(P_d,A)\xrightarrow[]{\partial^*}\gHom_A(P_{d+1},A)\to \Coker\partial^*\to 0$$
    is exact, and so it is a minimal graded projective resolution of $\Coker \partial^*$. Hence $\pdim_{A^o}\Coker\partial^*=d+1$, which contradicts to C2
    as $\depth_{A^o} \Coker\partial^*\geqslant 0$.

    C3 $\Rightarrow$ C4. It follows from $\depth_AA=d$.

    C4 $\Rightarrow$ C5. For any graded simple $A$-module $M$, $\gExt^d_A(M,A)=\gHom_A(M,I^d)$. Note that $\gHom_A(M,I^d)\neq 0$ if and only if $\soc_A I^d$ contains some shifts of $M$. It follows that $\soc_A I^d$ contains every graded simple module up to shift. Hence $\soc_A I^d$ is a projective generator in $\gr S$.

    C5 $\Rightarrow$ C6. Since $\depth_A A=d$, $R^d\Gamma_A(A)$ is the kernel of the canonical map $\Gamma_A(I^d)\to \Gamma_A(I^{d+1})$.
    Since $\Ker (I^d\to I^{d+1})$ is an essential submodule of $I^d$,
    it follows from $\soc_A I^d=\soc_A \Gamma_A(I^d)$ that $\soc_A \Gamma_A(I^d)\subseteq \Ker(\Gamma_A(I^d)\to \Gamma_A(I^{d+1}))=R^d\Gamma_A(A)$. Hence $\soc_A R^d\Gamma_A(A)$ is a projective generator in $\gr S$.

    C6 $\Rightarrow$ C7. The Matlis duality induces a surjective morphism of graded $A^o$-modules $\pi:\Omega'=D(R^d\Gamma_A(A))\to D(\soc_AR^d\Gamma_A(A))$  such that $\Omega' J\subseteq \Ker\pi$. 
    It follows from (6) that $D(\soc_A R^d\Gamma_A(A))$ is a projective generator in $\gr S^o$. Hence $\Omega'/\Omega' J$ is a projective generator in $\gr S^o$.

    C7 $\Rightarrow$ C1. Assume $\idim_AX=m$ and let $Y=R\gHom_A(X,R)$. Then $\pdim_{A^o} Y=\idim_AX=m$ by Lemma \ref{depth, sup, projective and injective dim}. 
    It is clear that
    $$\sup R\gHom_{A^o}(Y,R)\leqslant \pdim_{A^o}Y+\sup R.$$

    Since $R=D(R\Gamma_A(A))$ and $\depth_A A=d$, $R$ is isomorphic to the following complex in $\D^-_{\fg}(A^o)$:
    $$\cdots\to D(\Gamma_A(I^{d+1})\to D(\Gamma_A(I^d))\to 0.$$
    Hence $\sup R=-d$ and $H^{-d}(R)=D(R^d\Gamma_A(A))=\Omega'$.
    Since $\Omega'$ is a quotient module of $D(\Gamma_A(I^d))$, there is an exact sequence of complexes
    $$0\to Z\to R\to \Omega'/\Omega' J[d]\to 0$$
    with $\sup Z\leqslant -d$. 
    It induces a long exact sequence
    $$\cdots \to \gExt_{A^o}^{m-d}(Y,R)\to \gExt_{A^o}^{m-d}(Y,\Omega'/\Omega' J [d])\to \gExt_{A^o}^{m-d+1}(Y,Z)\to \cdots.$$
    
Let $P$ be a minimal graded projective resolution of $Y$. Then $\inf \{i\mid P^i\neq 0\}=-\pdim_{A^o}Y=-m$ which implies that
$\sup R\gHom_{A^o}(Y,Z)\leqslant m-d$.
Hence $\gExt_{A^o}^{m-d+1}(Y,Z)=0$.
     
     On the other hand, $H^i(R\gHom_{A^o}(Y,S))=H^i(\gHom_{A^o}(P,S))=\gHom_{A^o}(P^{-i},S)$ for any $i$. It follows that $\gExt_{A^o}^m(Y,S)\neq 0$. Since $\Omega'/\Omega' J$ is a projective generator in $\gr S^o$, 
    $$\gExt_{A^o}^{m-d}(Y,\Omega'/\Omega' J[d])=\gExt_{A^o}^m(Y,\Omega'/\Omega' J)\neq 0.$$
    It follows that $\gExt_{A^o}^{m-d}(Y,R)\neq 0$. 
    
    As a result,
    $$\pdim_{A^o}Y+\sup R=m-d\leqslant \sup R\gHom_{A^o}(Y,R)\leqslant \pdim_{A^o}Y+\sup R.$$
    So, $\sup R\gHom_{A^o}(Y,R)= \pdim_{A^o}Y+\sup R$. Since
        $$\sup X=\sup R\gHom_A(A,X)=\sup R\gHom_{A^o}(R\gHom_A(X,R),R)=\sup R\gHom_{A^o}(Y,R),$$ it follows that
        $\sup X=\pdim_{A^o}Y+\sup R=m-d=\idim_AX-\depth_A A$.
\end{proof}

Balanced Cohen-Macaulay algebras are well studied in \cite{Mo} in the connected graded case. We give a definition of balanced Cohen-Macaulay algebra for commonly graded algebras. 

\begin{definition}\label{def of balanced CM}
    A noetherian commonly graded algebra $A$ is called balanced Cohen-Macaulay (balanced CM, for short) of dimension $d$, if there is a graded $(A,A)$-bimodule $\Omega$ such that $\Omega[d]$ is a balanced dualizing complex of $A$. In this case, 
    $\Omega$ is called a balanced dualizing module of $A$.
\end{definition}

In the rest of this section, we always assume that $A$ is a balanced CM algebra of dimension $d$, and $\Omega$ is the balanced dualizing module of $A$. By Theorem \ref{local duality}, the injective dimension of a balanced dualizing complex is $0$ on both sides, so the injective dimension of $\Omega$ is $d$ on both sides. By Theorem \ref{local duality} again, the cohomological dimension of $\Gamma_A$ (or $\Gamma_{A^o}$) is $d$.

\begin{definition}
    Let $A$ be a balanced CM algebra of dimension $d$. A finitely generated graded $A$-module $M$ is called maximal Cohen-Macaulay (MCM, for short) if $\depth_AM=d$.
\end{definition}

The full subcategory consisting of all MCM $A$-modules is denoted by $\MCM A$. 
For convenience, let $(-)^\vee=\gHom_A(-,\Omega)$ (or $\gHom_{A^o}(-,\Omega)$).
\begin{lemma}\label{property of MCM}
    Suppose $A$ is a balanced CM algebra of dimension $d$ with balanced dualizing module $\Omega$. The following are equivalent for $M\in \gr A$.
    \begin{itemize}
        \item [(1)] $M$ is an MCM $A$-module.
        \item [(2)] $\gExt_A^i(M,\Omega)=0$ for all $i>0$.
        \item [(3)] $M\cong M^{\vee\vee}$ and $\gExt_A^i(M,\Omega)=\gExt_{A^o}^i(M^\vee,\Omega)=0$ for all $i>0$.
    \end{itemize}
    In particular, $(-)^\vee$ induces a duality between $\MCM A$ and $\MCM A^o$.
\end{lemma}
\begin{proof}
    The proof is same as the connected graded case. See \cite[Lemma 4.6]{Mo}.
\end{proof}

By Lemma \ref{property of MCM} and the assumption that $\Omega[d]$ is a dualizing complex, the endomorphism ring $\gEnd_A(\Omega)$ is isomorphic to $A$.

Next, we recall some facts about Auslander-Buchweitz approximation.

Let $\cA$ be an abelian category, $\cB$ be a full subcategory of $\cA$ and $\cC$ be a full subcategory of $\cB$. Let $\hat{\cB}$ (resp. $\hat{\cC}$) be the full subcategory of $\cA$ consisting of objects admitting a $\cB$-resolution (resp. $\cC$-resolution) of finite length, that is, $X\in \hat{\cB}$ (resp. $X\in \hat{\cC}$) if and only if there is an exact sequence
$$0\to B_n\to \cdots\to B_0\to X\to 0$$
for some $n\geqslant 0$ and $B_i\in \cB$ (resp. $B_i\in \cC$) for all $0 \leqslant i \leqslant n$. 

If for any $B\in \cB$, there is an exact sequence
$$0\to B\to C\to B'\to 0$$
in $\cA$, where $C\in \cC$,  
then $\cC$ is called a cogenerator of $\cB$.

\begin{theorem}[Auslander-Buchweitz Approximation]\cite{AB}\label{AB approximation}
 Let $\cA$ be an abelian category and $\cC\subseteq \cB\subseteq \cA$ be full subcategories. Suppose $\cB$ is closed under finite direct sum, direct summands and extensions in $\cA$. Suppose $\cC$ is a cogenerator of $\cB$. Then for any $X\in \hat{\cB}$, there are exact sequences
 $$0\to Z\to Y\to X\to 0 \text{ and } 0\to X\to Z'\to Y'\to 0$$
 such that $Y,Y'\in \cB$ and $Z,Z'\in \hat{\cC}$.
\end{theorem}

For a graded $A$-module $M$, let $\add_A M$ be the category of graded $A$-modules which are isomorphic to the direct summands of the finite direct sums of the shifts of $M$.

Take $\cA=\gr A$, $\cB=\MCM A$ and $\cC=\add_A\Omega$. Clearly, $\cB$ is closed under finite direct sum, direct summands and extensions in $\gr A$. For any $M\in\gr A$, let $0\to K_d\to P_{d-1}\to \cdots \to P_0\to M\to 0$ be a resolution of $M$ such that $P_i$ is projective for $0\leqslant i\leqslant d-1$. 
Then $\depth_A K_d=d$ and $K_d$ is an MCM $A$-module. So $M$ has an ($\MCM A$)-resolution of length at most $d$ and $\gr A=\widehat{\MCM A}$.

For any $M\in \MCM A$, there is an exact sequence $0\to K\to P\to M^\vee\to 0$ where $P$ is a finitely generated projective $A^o$-module.
Applying $(-)^\vee$ to this exact sequence, we get an exact sequence $0\to M\to P^\vee\to K^\vee\to 0$ such that $P^\vee\in \add_A \Omega$.
It follows that $\add_A \Omega$ is a cogenerator of $\MCM A$. By Theorem \ref{AB approximation}, the next proposition holds.

\begin{proposition}\label{AB approximation for balanced CM}
    Let $A$ be a balanced CM algebra of dimension $d$ with balanced dualizing module $\Omega$. Then for any $X\in \gr A$, there are exact sequences
    $$0\to Z\to Y\to X\to 0 \text{ and } 0\to X\to Z'\to Y'\to 0$$
 such that $Y,Y'\in \MCM A$ and $Z,Z'\in \widehat{\add_A\Omega}$.
\end{proposition} 

Let $\overline{\add}_A\Omega$ be the full subcategory of $\gr A$ consisting of the modules $M$ having a finite coresolution $0\to M\to W^0\to \cdots\to W^n\to 0$ such that $W^i\in \widehat{\add_A\Omega}$.

The proof of the following lemma is similar to the one in \cite[Proposition 2.6(1)]{GN} for module-finite algebras.

\begin{lemma}\label{I(A) and Omega-coresolution}
    Let $A$ be a balanced CM algebra of dimension $d$ with  balanced dualizing module $\Omega$. Keep the notations as above. Then $\widehat{\add_A\Omega}\subseteq \overline{\add}_A\Omega=\tilde{I}(A)$.
\end{lemma}
\begin{proof}
    Obviously, $\widehat{\add_A\Omega}\subseteq \overline{\add}_A\Omega.$
    It follows from $\idim_A\Omega=d$ that $\overline{\add}_A\Omega\subseteq \tilde{I}(A)$.
    
    For any $X\in \tilde{I}(A)$, by Theorem \ref{AB approximation for balanced CM}, there is an exact sequence $0\to X\to Z\to Y\to 0$ where $Z\in \widehat{\add}_A\Omega$ and $Y\in \MCM A$. Then it follows from $\idim_AX<\infty$ that $\idim_AY<\infty$. Thus $Y\in \MCM A\cap \tilde{I}(A)$. 

    Since $\add_A\Omega$ is a cogenerator for $\MCM A$, there is a coresolution of $Y$:
    $$0\to Y\to W^0\xrightarrow[]{\partial^0} W^1\xrightarrow[]{\partial^1} \cdots$$
    where $W^i\in \add_A\Omega$. 
    On one hand, $\idim_AY<\infty$ implies that $\gExt_A^{n+1}(\im \partial^n,Y)=0$ for $n\gg 0$. 
    On the other hand, $\im\partial^n$ is MCM which is equivalent to say that
    for any $i>0$, $\gExt^i(\im\partial^n,\Omega)=0$. So
    $$\gHom_A(\im \partial^n,W^n)\to \gHom_A(\im \partial^n,W^{n+1})\to \gHom_A(\im\partial^n,W^{n+2})$$
    is exact. It follows that $\partial^n$ splits. 
    Then the coresolution of $Y$ is finite:
    $$0\to Y\to W^0\to \cdots \to W^m\to 0.$$
    Combining this coresolution with $0\to X\to Z\to Y\to 0$ shows that $X\in \overline{\add}_A\Omega$.
\end{proof}

Now we prove the main result in this section.
\begin{theorem}\label{AB formula and Bass thm for balanced CM alg}
    Let $A$ be a balanced CM algebra of dimension $d$ with  balanced dualizing module $\Omega$ and $0\to A\to I^0\to \cdots I^d\to \cdots$ be a minimal graded injective resolution of ${}_AA$. 
    Then C1 to C7 are also equivalent to any of the following conditions. 
    \begin{itemize}
    \item [C$1^{\prime}$] $\idim_AM=d$ for any $M\in \tilde{I}(A)$.
    \item [C8\hphantom{'}] $\idim_AM\leqslant d$, for any $M\in \tilde{I}(A)$.
    \item [C9\hphantom{'}] $\idim_A M=d$ for any $M\in \MCM A\cap \tilde{I}(A)$.
    \item [C10] 
    $\gExt^i_A(M,X)=0$ for any $i>0$, $M\in \MCM A$ and $X\in \tilde{I}(A)$.
    \item [C11] 
    $\gExt^i_A(\Omega,X)=0$ for any $i>0$ and $X\in \tilde{I}(A)$.
    \item [C12] $\widehat{\add_A\Omega}=\tilde{I}(A)=\overline{\add}_A\Omega$.
    \item [C13] If $0\to M\to W^0\to W^1\to 0$ is exact with $W^i\in \add_A\Omega$, then $M\in \add_A\Omega$.
    \item [C14] If $0\to M\to W^0\to\cdots \to W^n\to 0$ is exact with $W^i\in \add_A\Omega$, then $M\in \add_A\Omega$.
    \item [C15] If $M\in \MCM A$ and $0\to M\to W^0\to\cdots \to W^n\to 0$ is exact with $W^i\in \add_A\Omega$, then $M\in \add_A\Omega$.
    \item [C16] $\add_{A^o}A=\MCM A^o\cap \tilde{P}(A^o)$.
    \item [C17] $\add_A\Omega=\MCM A\cap \tilde{I}(A)$.
    \end{itemize}
\end{theorem}
\begin{proof}
    C1 $\Rightarrow$ C$1^{\prime}$ $\Rightarrow$ C8. Clearly.

    C$2^{\prime}$ $\Rightarrow$ C17. Obviously, $\add_A\Omega\subseteq \MCM A\cap \tilde{I}(A)$. Let $M\in \MCM A\cap \tilde{I}(A)$. By Lemma \ref{depth, sup, projective and injective dim} (3), $\pdim_{A^o} M^\vee$ is finite. 
    Then C$2^{\prime}$ yields that $M^\vee$ is a projective $A^o$-module as it is MCM. 
    Hence $M\cong M^{\vee\vee}\in \add_A \Omega$.

    C17 $\Rightarrow$ C$2^{\prime}$. Suppose $N\in \tilde{P}(A^o)$ and  $\depth_{A^o} N=n$. 
    Let $0\to P_s\to\cdots\to P_0\to N\to 0$ be a graded projective resolution of $N$. 
    If $n=d$, then $N$ is MCM. 
    Thus $0\to N^\vee\to P_0^\vee\to \cdots P_s^\vee\to 0$ is exact. 
    As $P_i^\vee\in \add_A \Omega$, $N^\vee$ has finite injective dimension. 
    By assumption, $N^\vee\in\add_A\Omega$. As a result, $N\cong N^{\vee\vee}$ is projective.

    If $n<d$, consider the exact sequence $0\to K\to P_0\to N\to 0$. 
    Then $\depth_{A^o} K=n+1$ and $\pdim_{A^o}K=\pdim_{A^o}N-1$. By a backward induction on the depth, $\depth_{A^o} K+\pdim_{A^o}K=d$ and thus $\depth_{A^o} N+\pdim_{A^o}N=d$.

    C17 $\Rightarrow$ C8. Suppose $M\in \tilde{I}(A)$. By Proposition \ref{AB approximation for balanced CM}, there is an exact sequence $0\to Z\to Y\to M\to 0$ where $Z\in \widehat{\add_A \Omega}$ and $Y\in \MCM A$. 
    If $\depth_A M=d$, that is, $M$ is MCM, then by (C17) $M\in \add_A\Omega$, which implies that $\idim_AM=d$. 
    If $\depth_A M=m<d$, then $\depth_A Z=m+1$. 
    Since $Z\in \widehat{\add_A\Omega}$, $Z$ has finite injective dimension. 
    By a backward induction on the depth, $\idim_A Z\leqslant d$. 
    The exact sequence $0\to Z\to Y\to M\to 0$ also yields that $\idim_A Y$ is finite, 
    which implies that
    $Y\in \MCM A\cap \tilde{I}(A)=\add_A\Omega$. So, $\idim_AY=d$ and $\idim_AM\leqslant d$.

    C8 $\Rightarrow$ C9. The fact $R^d\Gamma(M)\neq 0$ and  $\idim_AM\leqslant d$ imply that $\idim_AM=d$.

    C9 $\Rightarrow$ C17. Suppose $M \in \MCM A\cap \tilde{I}(A)$. By Lemma \ref{depth, sup, projective and injective dim} (3), $\pdim_{A^o}M^\vee=0$, that is, $M^\vee$ is a finitely generated projective $A^o$-module. So $M\cong M^{\vee\vee} \in \add_A \Omega$.

    C17 $\Rightarrow$ C10. 
    Suppose $X\in \tilde{I}(A)$ and  $\depth_A X=m$.
    If $m=d$, then $X\in \add_A\Omega$. It follows that C10 holds for $X$ and any $M\in \MCM A$. 
    If $m<d$, then by Proposition \ref{AB approximation for balanced CM}, there is an exact sequence $0\to Z\to Y\to X\to 0$ where $Y\in \MCM A$ and $Z\in \widehat{\add_A\Omega}$. 
    It implies that
    $\depth_A Z=m+1$, $\idim_AZ$ is finite. So, $\idim_AY$ is finite. 
    Then $Y\in\add_A\Omega$. By a backward induction on the depth we may assume that C10 holds for $Z$ and any $M\in \MCM A$. 
   It follows from the exact sequence $0\to Z\to Y\to X\to 0$ that C10 holds for $X$ and any $M\in \MCM A$.
    
    C10 $\Rightarrow$ C11. Clearly.

    C11 $\Rightarrow$ C12. By Lemma \ref{I(A) and Omega-coresolution}, it suffices to prove $\tilde{I}(A)\subseteq\widehat{\add_A\Omega}$. Note that the proofs and the results in appendix of \cite{NV} still work for commonly graded algebras. 
    So by \cite[Proposition A.2]{NV}, for any $X\in \tilde{I}(A)$, $R\gHom_A(\Omega,X)$ has finitely generated cohomologies and finite projective dimension.
    By C11, $R\gHom_A(\Omega,X)\cong \gHom_A(\Omega,X)$. 
    Then $R\gHom_A(\Omega,X)$ and $\gHom_A(\Omega,X)$ share a same finitely generated projective resolution $P_\bullet$ of finite length.
    
    Since $\idim_A X<\infty$,
    $$\Omega\otimes_A^L R\gHom_A(\Omega,X)\cong R\gHom_A(R\gHom_{A^o}(\Omega,\Omega),X)\cong X.$$
    Then
    $$\Omega\otimes_A P_\bullet\cong \Omega\otimes_A^L R\gHom_A(\Omega,X)\cong X,$$
    that is, $\Omega\otimes_A P_\bullet$ is a resolution of $X$. For any $i$, since $P_i$ is finitely generated projective, $\Omega\otimes_A P_i$ is in $\add_A\Omega$. So C12 holds.

    C12 $\Rightarrow$ C17. Let $M\in \MCM A\cap \tilde{I}(A)$. Then  $M \in \widehat{\add_A\Omega}$, and there is a resolution $0\to \Omega^n\to\cdots\to \Omega^0\to M\to 0$ where $\Omega^i\in \add_A\Omega$.

    If $n=0$, then $M\in \add_A\Omega$. 
    Now suppose $n\geqslant 1$. Consider 
    the exact sequences $0\to \Omega^n\to \Omega^{n-1}\to X\to 0$ and 
    $0\to X\to \Omega^{n-2}\to \cdots\to \Omega^0\to M\to 0.$
    Since $M$ is MCM, $X$ is MCM. So $\gExt_A^i(X,\Omega)=0$ for all $i>0$, which implies that
    $$0\to \gHom_A(X,\Omega^n)\to \gHom_A(X,\Omega^{n-1})\to \gHom_A(X,X)\to 0$$
    is exact. It follows that the exact sequence $0\to \Omega^n\to \Omega^{n-1}\to X\to 0$ splits. Hence $X\in \add_A\Omega$. It follows an induction that $M \in \add_A\Omega$.

    C12 $\Leftrightarrow$ C13 $\Leftrightarrow$ C14 $\Leftrightarrow$ C15. See \cite[4.3, 4.6]{AB}.

    C16 $\Leftrightarrow$ C17. It follows from Lemma \ref{property of MCM} that $M\in \MCM A\cap \tilde{I}(A)$ if and only if $M^\vee\in \MCM A^o\cap \tilde{P}(A^o)$. Since $\gEnd_A(\Omega)\cong A$, $M\in\add_A\Omega$ if and only if $M^\vee\in\add_{A^o}A$.
\end{proof}

As mentioned in Introduction, we simply call these equivalent conditions $\mathbf{(C)}$. When we consider right modules, these conditions are denoted by $\mathbf{(C^o)}$.

Next corollary is proved in \cite[Theorem 5.8]{Mo} for connected graded algebras.
\begin{corollary}
    Keep the notations as above. If $A$ satisfies $\mathbb{(C)}$, then for any $M\in \tilde{I}(A)$, $M$ has a resolution
$$0\to W_m\to \cdots\to W_0\to M\to 0$$
where $W_i\in \add_A\Omega$ and $m=d-\depth_AM$.
\end{corollary}
\begin{proof}
    It follows from C12 and C17.
\end{proof}

Now, we prove that No-Hole theorem holds for any balanced CM algebra satisfying $\mathbf{(C)}$.

\begin{lemma}\label{prepare for no-hole}
    Let $A$ be a noetherian commonly graded algebra with a balanced dualizing complex and $S=A/J$. For any $X\in \D^b_{\fg}(\Gr A)$, if $\gExt_A^i(S,X)=0$ for some $i$, then there is a decomposition $X=Y\oplus Z$ in $\D^b_{\fg}(\Gr A)$ such that
    $$\idim_A Y<i<\depth_A Z.$$
\end{lemma}
\begin{proof}
    The proof in \cite[Lemma 4.7]{Jo1} for connected graded algebras still works in the setting of the lemma.
\end{proof}

\begin{theorem}\label{no-hole thm}
    Let $A$ be a balanced CM algebra of dimension $d$ and $S=A/J$. If $A$ satisfies $\mathbf{(C)}$, then for any $M\in \gr A$,
    $$\gExt_A^i(S,M)\neq 0 \text{ if and only if } \depth_A M\leqslant i\leqslant \idim_AM.$$
\end{theorem}
\begin{proof}
    ``$\Rightarrow$" It is clear.

    ``$\Leftarrow$" On one hand, Theorem \ref{AB formula and Bass thm for balanced CM alg} assures that $\idim_AN\geqslant d$ for any $0 \neq N\in\gr A$. 
    Since the cohomological dimension of $\Gamma_A$ is $d$, by Lemma \ref{depth and local dimension}, $\depth_A N\leqslant d$. So for any non-zero $N_1$ and $N_2\in \gr A$, $\depth_AN_1\leqslant d\leqslant \idim_A N_2$.

    On the other hand, by Lemma \ref{prepare for no-hole}, if $\gExt_A^i(S,M)=0$ for some $i$, then there is a decomposition $M=M_1\oplus M_2$ such that
    $$\idim_AM_1<i<\depth_AM_2.$$
    It follows that $M_1=0$ or $M_2=0$. 
    
    If $M_1=0$, then $i<\depth_AM_2=\depth_AM$. If $M_2=0$, then $i>\idim_AM_1=\idim_AM$.
\end{proof}

A connected graded algebra with a balanced dualizing complex automatically satisfies $\mathbf{(C)}$ and $\mathbf{(C^o)}$. 
When $A$ is a noetherian commonly graded AS-Gorenstein algebra of dimension $d$, then $A$ is balanced CM of dimension $d$ and C4 holds for $A$, so $A$ satisfies $\mathbf{(C)}$ and $\mathbf{(C^o)}$.

\begin{corollary}
    Let $A$ be a noetherian commonly graded AS-Gorenstein algebra. Then $A$ satisfies $\mathbf{(C)}$ and $\mathbf{(C^o)}$.  In particular,
    \begin{itemize}
        \item [(1)] (Auslander-Buchsbaum formula) If $M\in \gr A$ has finite projective dimension, then
        $$\depth_AM+\pdim_AM=\depth_AA.$$
        \item [(2)] (Bass Theorem) If $M\in \gr A$ has finite injective dimension, then
        $$\idim_AM=\depth_AA.$$
        \item [(3)] (No-Hole Theorem) If $M\in \gr A$, then
        $$\gExt_A^i(S,M)\neq 0 \text{ if and only if } \depth_AM<i<\idim_AM.$$
    \end{itemize}
\end{corollary}

The fact that the Auslander-Buchsbaum formula holds for noetherian commonly graded AS-Gorenstein algebras is also proved in \cite{HY}.

Next theorem is a generalization of \cite[Corollary 5.9]{Mo} and \cite[Theorem 3.6]{DW}.
\begin{theorem}
    Let $A$ be a noetherian commonly graded algebra with a balanced dualizing complex $R$. If $A$ satisfies  $\mathbf{(C)}$ and $\mathbf{(C^o)}$, then the following are equivalent.
    \begin{itemize}
        \item [(1)] $A$ is commonly graded AS-Gorenstein.
        \item [(2)] $\idim_AA<\infty$.
        \item [(3)] $\pdim_AR<\infty$.
        \item [(4)] For any $X\in\D_{\fg}^b(\Gr A)$, $\pdim_AX<\infty$ if and only if $\idim_A X<\infty$.
    \end{itemize}
\end{theorem}
\begin{proof}
    By Corollary \ref{depth A and A^o}, we may assume $\depth_AA=\depth_{A^o}A=d$.

    (1) $\Rightarrow$ (4). Suppose $\pdim_AX<\infty$. Let $P$ be a minimal projective resolution of $X$. 
    Note that $P$ is bounded and every $P^i$ has finite injective dimension.
    It follows that $\idim_AP$ is finite and so is $X$.

    Suppose $\idim_AX<\infty$. Then $\pdim_{A^o}R\gHom_A(X,R)<\infty$. By a dual argument of the previous paragraph, $\idim_{A^o}R\gHom_A(X,R)$ is finite, which implies that 
    $\pdim_AX$ is finite by Lemma \ref{dualizing complex induces a duality}.

    (4) $\Rightarrow$ (3). Obviously.

    (3) $\Rightarrow$ (2). By Lemma \ref{dualizing complex induces a duality}, $\idim_{A^o}A<\infty$. Then C$1^o$ ( the Bass Theorem) yields that $\idim_{A^o}A=d=\depth_{A^o}A$. 
    Hence $R\cong D(R^d\Gamma_A(A))[d]$ and $A$ is balanced CM. 

    By the local duality Theorem \ref{local duality},  $R\gHom_A(S,R)\cong D(R\Gamma_A(S))\cong S$, which implies that $\depth_AR=0$. By C2 (AB-formula), $\pdim_AR=d$. 
    Hence $\pdim_A D(R^d\Gamma_A(A))=0$ and $\Omega:=D(R^d\Gamma_A(A))$ is finitely generated projective. 
    Then C$7^o$ makes $\Omega$ a finitely generated projective generator in $\gr A$. So it follows from $\idim_A\Omega=\idim_AR+d=d$ that $\idim_AA=d$.

    (2) $\Rightarrow$ (1). By Lemma \ref{dualizing complex induces a duality}, $\pdim_{A^o}R\gHom_A(A,R)=R<\infty$. Similar to (3) $\Rightarrow$ (2), $\idim_{A^o}A<\infty$. 
    Then C1 and C$1^o$ show that $\idim_AA=\idim_{A^o}A=\depth_AA=d$. Hence $A$ is balanced CM. The proof of $\pdim_AR=\pdim_{A^o}R=d$ is also similar as in (3) $\Rightarrow$ (2).

    Since $A$ is balanced CM, $R=\Omega[d]$ where $\Omega=D(R^d\Gamma_A(A))$. Then $\pdim_A\Omega=\pdim_{A^o}\Omega=0$. 
    It follows that $R\gHom_A(R,R)=\gHom_A(\Omega,\Omega)\cong A$  that $\Omega$ is an invertible graded $(A,A)$-bimodule. 
    By Theorem \ref{another equivalent definition of GAS}, $A$ is commonly graded AS-Gorenstein.
\end{proof}

\end{document}